\documentclass[a4paper,16pt]{article}
\usepackage{authblk} % need it for affiliation
\usepackage[english]{babel}
\usepackage[utf8x]{inputenc}
\usepackage[T1]{fontenc}
\usepackage{algorithm}
\usepackage[noend]{algpseudocode}
\usepackage{url}
\usepackage{xspace}
\usepackage{multirow}
\usepackage[a4paper,top=2cm,bottom=2cm,left=2cm,right=2cm,marginparwidth=1.75cm]{geometry}

\usepackage{amsmath}
\usepackage{graphicx}
\usepackage[colorinlistoftodos]{todonotes}

\usepackage{bm}
\usepackage{amssymb}
\usepackage{amsfonts}
\usepackage{tikz} \usetikzlibrary{positioning,shapes}
\usepackage{pgfplots} \usetikzlibrary{plotmarks}
\usetikzlibrary{calc,backgrounds}

\DeclareMathAlphabet{\mathitbf}{OML}{cmm}{b}{it}
\def\jmath{j}

\def\jmath{j}

\usepackage[labelformat=simple]{subcaption}

\newcommand\myeq{\mathrel{\overset{\makebox[0pt]{\mbox{\normalfont\tiny\sffamily def}}}{=}}}

\newcommand{\vTheta}{\varTheta}

\newcommand{\Psii}{\mathcal{F}}

\newcommand{\EXP}[1]{\mathbb{E}\left[#1\right ]}
\newcommand{\Var}[1]{\mathbb{V}\left[#1\right ]}
\newcommand{\di}{\mathrm{d}}
\newcommand{\btheta}{\mbox{\boldmath$\theta$}}

\def\min{\mbox{min}}

\def\sol{c}

\def\D{\mathcal{D}}

\newcommand{\var}[1]{{\ensuremath{\mathrm{Var}}\mspace{-2mu}\left[#1\right]}}

\newcommand{\MSE}{\text{MSE}}

\def\xib{\bm{\xi}}

\def\bx{\mathbf{x}}

\providecommand{\Order}[1]{ {\ensuremath{ \mathcal O\left( #1 \right)}} }

\newcommand{\conc}{c} % mass fraction
\newcommand{\pres}{p} % hydrostatic pressure
\newcommand{\poro}{\phi} % porosity
\newcommand{\perm}{{\mathbf{K}}} % permeability
\newcommand{\dens}{\rho} % density
\newcommand{\visc}{\mu} % viscosity
\newcommand{\dvel}{{\mathbf{q}}} % Darcy velocity
\newcommand{\grav}{{\mathbf{g}}} % gravity
\newcommand{\disp}{{\mathbf{D}}} % diffusion-disperion
\newcounter{theorem}

\newtheorem{defn}[theorem]{Definition}
\newtheorem{theorem}{Theorem}

\newcommand{\Tau}{\mathcal{T}}

\title{Uncertainty quantification in coastal aquifers using the multilevel Monte Carlo method}
\author[1]{Alexander Litvinenko}
\author[2]{Dmitry Logashenko}
\author[1,2]{Raul Tempone}
\author[3]{Ekaterina Vasilyeva}
\author[2]{Gabriel Wittum}
\affil[1]{RWTH Aachen, Aachen, Germany}
\affil[2]{KAUST, Thuwal-Jeddah, Saudi Arabia}
\affil[3]{Goethe-Universit\"at Frankfurt am Main, Germany}
\affil[ ]{\textit{litvinenko@uq.rwth-aachen.de,ekaterina.vasilyeva@gcsc.uni-frankfurt.de,
\{raul.tempone,dmitry.logashenko,gabriel.wittum\}@kaust.edu.sa}}
\makeatletter
\def\BState{\State\hskip-\ALG@thistlm}
\makeatother

\begin{document}
\maketitle

%\renewcommand{\thefootnote}{\fnsymbol{footnote}}

%\footnotetext[2]{E-mail: \email{david.keyes@kaust.edu.sa}, %\email{alexander.litvinenko@kaust.edu.sa}. King Abdullah University of Science and %Technology (KAUST), Thuwal, Saudi Arabia.}
%\renewcommand{\thefootnote}{\arabic{footnote}}

\begin{abstract}
We consider a class of density-driven flow problems. We are particularly interested in the problem of the salinization of coastal aquifers. We consider the Henry saltwater intrusion problem with uncertain porosity, permeability, and recharge parameters as a test case. 
The reason for the presence of uncertainties is the lack of knowledge, inaccurate measurements,
and inability to measure parameters at each spatial or time location. This problem is nonlinear and time-dependent. The solution is the salt mass fraction, which is uncertain and changes in time. Uncertainties in porosity, permeability, recharge, and mass fraction are modeled using random fields. This work investigates the applicability of the well-known multilevel Monte Carlo (MLMC) method for such problems. The MLMC method can reduce the total computational and storage costs. Moreover, the MLMC method runs multiple scenarios on different spatial and time meshes and then estimates the mean value of the mass fraction. 
The parallelization is performed in both the physical space and stochastic space. To solve every deterministic scenario, we run the parallel multigrid solver ug4 in a black-box fashion. 
We use the solution obtained from the quasi-Monte Carlo method as a reference solution.
\end{abstract}

%\begin{keywords}
%\end{keywords}

%\begin{AMS}
%15A69,  %       Multilinear algebra, tensor products
%65F10,  % Iterative methods for linear systems
%60H15,  %  Stochastic partial differential equations
%60H35,  %  Computational methods for stochastic equations 
%65C30  %  Stochastic differential and integral equations
%\end{AMS}

%

\section{Introduction}
\label{sec:1}

\begin{table}[htbp!]
\begin{tabular}{|c|l|}
\hline
\multicolumn{2}{|c|}{\textbf{Notation}} \\ \hline
%${\Nset}$        &Natural numbers \\ \hline
QoI $g$ & quantity of interest $g$ \\ \hline
$\mathcal{D}$ & computational spatial domain \\ \hline
$\mathcal{D}_0,\mathcal{D}_1,\ldots,\mathcal{D}_L$ & hierarchy of spatial meshes \\ \hline
$\mathcal{T}_0,\mathcal{T}_1,\ldots,\mathcal{T}_L$ & hierarchy of temporal meshes \\ \hline
$L$ & number of levels \\ \hline
$s$ & complexity \\ \hline
$h_{\ell}$ (or $h$), $n_{\ell}$ & spatial step size and number of spatial degrees of freedom on level $\ell$ \\ \hline
$\tau_{\ell}$ (or $\tau$), $r_{\ell}$ & time step size and number of time steps on level $\ell$ \\ \hline
$m_{\ell}$ & number of samples (scenarios) on level $\ell$ \\ \hline
$\EXP{\cdot}$, $\Var{\cdot}$           & expectation and variance\\ \hline
$\vTheta$ & multidimensional domain of integration in parametric space \\ \hline
$\omega$, $\xib(\omega)$ & random event and random vector \\ \hline
$\poro(\bx,\omega)$ & porosity random field  \\ \hline
$\perm(\bx,\omega)$ & permeability random field \\ \hline
$\dens(\bx, \omega)$ & density random field  \\ \hline
$\dvel(t,\bx, \omega)$ & volumetric velocity \\ \hline
%$\conc(t,x,\omega)$ & mass fraction\\ \hline
$\disp$ &  tensor field $\disp = \disp (\dvel)$: molecular diffusion and dispersion of salt \\ \hline
$\bar{\kappa}(\bx)$ & expectation of $\kappa(\bx,\omega)$ \\ \hline
$d$ & physical (spatial) dimension \\ \hline
% $\param$ & high-dimensional parameter \\ \hline
%$\f$            & the right-hand side \\ \hline
$\sol=\sol(t,\bx,\omega)$          & mass fraction of salt (solution of the problem)\\ \hline
\end{tabular}
\label{tab:notations}
\end{table}

% An accurate modeling and predictions of a groundwater flow is vital.
Saltwater intrusion occurs when sea levels rise and saltwater moves onto the land. Usually, this occurs during storms, high tides, droughts, or when saltwater penetrates freshwater aquifers and raises the groundwater table. Since groundwater is an essential nutrition and irrigation resource, its salinization may lead to catastrophic consequences. Many acres of farmland may be lost because they can become too wet or salty to grow crops. Therefore, accurate modeling of different scenarios of saline flow is essential \cite{Abarca07,SWLRHEGW-SaltwaterInNorthSea2018} to help farmers and researchers develop strategies to improve the soil quality and decrease saltwater intrusion effects.

Saline flow is density-driven and described by a system of time-dependent nonlinear partial differential equations (PDEs). It features convection dominance and can demonstrate very complicated behavior \cite{Voss_Souza}.
%In time, the flow and concentration fields may converge to several steady states depending on the initial conditions (cf.\ \cite{Johannsen2003}).

As a specific model, we consider a Henry-like problem with uncertain permeability and porosity. 
These parameters may strongly affect the flow and transport of salt. The original Henry saltwater intrusion problem was introduced by H.R. Henry in the 1960s (cf. \cite{henry1964effects}). The Henry problem became a benchmark for numerical solvers for groundwater flow (cf. \cite{Voss_Souza,Simpson04_Henry,Simpson2003,Dhal_review22}. In \cite{Riva2015}, the authors use the generalized polynomial chaos expansion approximation to investigate how incomplete knowledge of the system properties influences the assessment of global quantities. Particularly, they estimated the propagation of input uncertainties into a few dimensionless scalar parameters.

The hydrogeological formations typically have complicated and heterogeneous structures.
These formations may consist of a few layers of porous media with various porosity and permeability coefficients
(cf.\ \cite{ReiterLogashenkoVogelWittum2017, ScheiderKroehnPueschel2012}).
Measurements of the layer positions and their thicknesses are only possible up to some error, and
for the materials inside the layers, the average parameters are typically assumed. Thus, these layers are excellent candidates to be modeled by random fields. Further, due to the nonlinearities in the problem, averaging the parameters does not necessarily lead to the correct mathematical expectation of the solution.

To model uncertainties, we use random fields. Uncertainties in the input data propagate through the model and make the solution (e.g., the mass fraction) uncertain. An accurate estimation of the output uncertainties can facilitate a better understanding of the problem, better decisions, and improved control and design of the experiment. 

The following questions can be answered:
\begin{enumerate}
\item How long can a particular drinking water well be used (i.e., when will the mass fraction of the salt exceed a critical threshold)?
\item What regions have especially high uncertainty?
%\item After how many days the salt concentration exceeds some dangerous threshold? 
\item What is the probability that the salt concentration is higher than a threshold at a certain spatial location and time point?
\item What is the average scenario (and its variations)?
\item What are the extreme scenarios?
\item How do the uncertainties change over time?
\end{enumerate}

Many techniques can quantify uncertainties. A classical method is Monte Carlo (MC) sampling. Although it is dimension-independent, it converges very slowly and requires many samples. This method may not be affordable for time-consuming simulations. Nevertheless, even up-to-date techniques, such as surrogate models and stochastic collocation, require a few hundred to a few thousand time-consuming simulations and assume a certain smoothness by the quantity of interest (QoI).

Another class of methods is the class of perturbation methods \cite{CREMER15_Fingers}. The idea is to decompose the QoI with respect to (w.r.t.) random parameters in a Taylor series. The higher-order terms can be neglected for small perturbations, simplifying the analysis and numerics. These methods assume that random perturbations are small (e.g., up to 5\% of the mean, depending on the problem). For larger perturbations, these methods usually do not work.

There are quite a few studies where authors model uncertainties in reservoirs (cf. \cite{OverviewUncert93,SoilOverview16}). Reconnecting stochastic methods with hydrogeological applications was accomplished in \cite{NowakStochMethods18}, where the authors analyzed a collaboration between academics and water suppliers in Germany and made recommendations regarding optimization and risk assessment. The fundamentals of stochastic hydrogeology and an overview of stochastic tools and accounting for uncertainty are described in \cite{rubin2003applHydro}. 

The review \cite{TARTAKOVSKYI_Risk13} deals with hydrogeologic applications of recent advances in uncertainty quantification, probabilistic risk assessment, and decision-making under uncertainty. The author reviewed probabilistic risk assessment methods in hydrogeology under parametric, geologic, and model uncertainties. Density-driven vertical transport of saltwater through the freshwater lens on the island of Baltrum (Germany) is modeled in 
\cite{POST17_Density-driven}.

In \cite{Laattoe2013_SeawaterIntr}, the authors examined the implications of transgression for a range of seawater intrusion scenarios based on simplified coastal freshwater aquifer settings. They stated that vertical intrusion during transgressions could involve density-driven convective processes, causing substantially greater amounts of seawater to enter the aquifer and create more extensive intrusion than horizontal seawater intrusion in the absence of transgression.

The methods to compute the desired statistics of the QoI are direct integration methods, such as the MC, quasi-MC (QMC) and collocation methods and surrogate-based (generalized polynomial chaos approximation and stochastic Galerkin \cite{Philipp12, babuska2004galerkin,GiraldiLitv14,wahnert-stochgalerkin-2014}) methods. Direct methods compute statistics directly by sampling uncertain input coefficients and solving the corresponding PDEs, whereas the surrogate-based method computes a cheap functional (polynomial, exponential, or trigonometrical) approximation of the QoI. Examples of the surrogate-based methods are radial basis functions \cite{liu2014,bompard2010,Loeven2007,giunta2004}, sparse polynomials \cite{
chkifa-adapt-stochfem-2015,Sudret_sparsePCE,DolgLitv15}, and polynomial chaos expansion \cite{Habib09_PCE, ConradMarzouk13, Dongbin}. Sparse grid methods to integrate high-dimensional integrals are considered in \cite{smoljak63, Griebel_Bungartz, Griebel, spiterp, novakRitter97, gerstnerGriebel98-numint, novakRitter99-simple, ConradMarzouk13,petrasSmolpak}. An idea to generate goal-oriented adaptive spatial grids and use them in the multilevel MC (MLMC) framework was presented in \cite{EIGEL14,BECK22}.

The quantification of uncertainties in stochastic PDEs could be a significant challenge due to a) the large possible number of involved random variables and b) the high cost of each deterministic solution of the governed PDE. The MC quadrature and its variance-reduced variants have a dimension-independent error convergence rate $\mathcal{O}(N^{-\frac{1}{2}})$, and the QMC has the worst-case rate $\mathcal{O}(\log^M(N)N^{-1})$, where $N$ is the number of samples, and $M$ indicates the dimension of the stochastic space \cite{matthies2007}. The MC method is not affected by the dimension of the
integration domain, such as collocations on sparse or full grid methods \cite{babuska_collocation, nobile-sg-mc-2015}. A numerical comparison of other QMC sequences is presented in \cite{radovic1996}. 

Construction of a cheap generalized polynomial chaos expansion-based surrogate model \cite{xiuKarniadakis02a,Litvinenko-UQ-2021, LitLog3D_20} is an alternative to the MC method. 
Some well-known functions, such as the multivariate Legendre, Hermite, Chebyshev, or Laguerre functions, have been taken as a basis \cite{OLADYSHKIN_PCE,xiuKarniadakis02a}. Surrogate models have pros and cons. The pros are that the model can be easily sampled once constructed, and all samples are almost free (much cheaper than sampling the original stochastic PDE). For some problem settings, sampling is unnecessary because the solution can be computed analytically (e.g., computing an integral of a polynomial). The nontrivial part of surrogate models is to define how many coefficients are needed and how accurately they should be computed. Another difficulty is that not every function can be approximated well by a polynomial. The MLMC methods do not have such limitations.

This work is structured as follows. Section~\ref{sec:Model} describes the Henry problem and numerical methods to solve it. The well-known MLMC method is reviewed in Section~\ref{sec:MLMC}. Next, Section~\ref{sec:numerics} details the numerical results, which include the numerical analysis of the Henry problem, computing different statistics, the performance of the MLMC method, and the practical performance of the parallel ug4 solver for the Henry problem \cite{henry1964effects,Simpson04_Henry} with uncertain coefficients. Finally, we conclude this work with a discussion in Section~\ref{sec:Conclusion}.\\

\textbf{Our contribution:} We investigate the propagation of uncertainties in the Henry-like problem. Assuming the porosity, permeability, and recharge are uncertain, we estimate the uncertainties in the density-driven flow. To reduce the high computing complexity, we applied the existing MLMC technique. We use the multigrid ug4 software library as a black-box solver, allowing us to solve the Henry problem and others (see more in \cite{SWLRHEGW-SaltwaterInNorthSea2018}). We run all MLMC random simulations in parallel. To the best of our knowledge, we are unaware of any other studies where Henry's problem \cite{henry1964effects,Simpson04_Henry} was solved using MLMC methods with uncertain porosity, permeability, and recharge parameters.

% Note, that these statistics could be further used for more efficient Bayesian inference, data assimilation, the optimal design of experiment and optimal control. The overall goal is modeling of water pollution and monitoring of quality of subsurface water flow.%
%
%
%%%%%%%%%%%%%%%%%%%%%%%%%%%%%%
\section{Henry Problem with Uncertain Porosity and Permeability}
\label{sec:Model}
\subsection{Problem setting}
\label{subsec:Henry}
In coastal aquifers, salty seawater intruding on the formation on one side (the seaside) displaces the pure water due to water recharge from land sources and precipitation from the other side. Due to its higher density, seawater mainly penetrates along the bottom of the aquifer. This process can achieve a steady state but may be time-dependent due to the periodicity of the recharge or controlling the pumping rate from the wells. An accurate simulation of the salinization is vital for the prediction of water resource availability. However, the accuracy of such predictions strongly depends on the hydrogeological parameters of the formation and the geometry of the computational domain, denoted by $\D$.

% {\color{red} For some alternative formulation, namely, the hydraulic-head formulation, see \cite{Langevin20}.}

The aquifer $\D \subset \mathbb{R}^d$, $d \in \{2, 3 \}$, can be modeled as an immobile porous matrix filled with liquid phase---a solution of salt in water. Due to the nonhomogeneous density distribution, gravitation induces the motion of the liquid phase. This motion transports the salt, which is otherwise subject to molecular diffusion.

\begin{figure}[htbp!]
\begin{center}
  \includegraphics[width=0.47\textwidth]{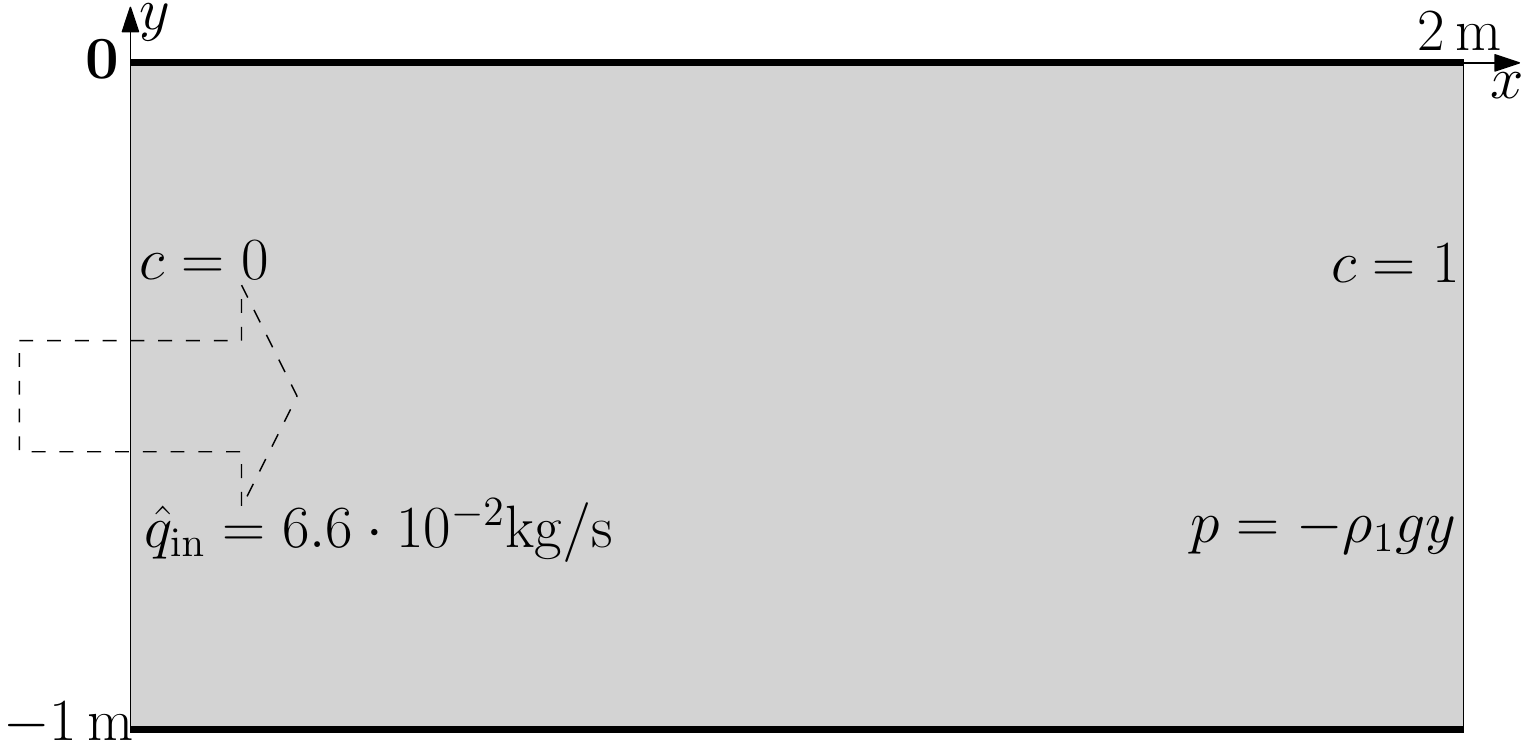}\;
  \includegraphics[width=0.43\textwidth]{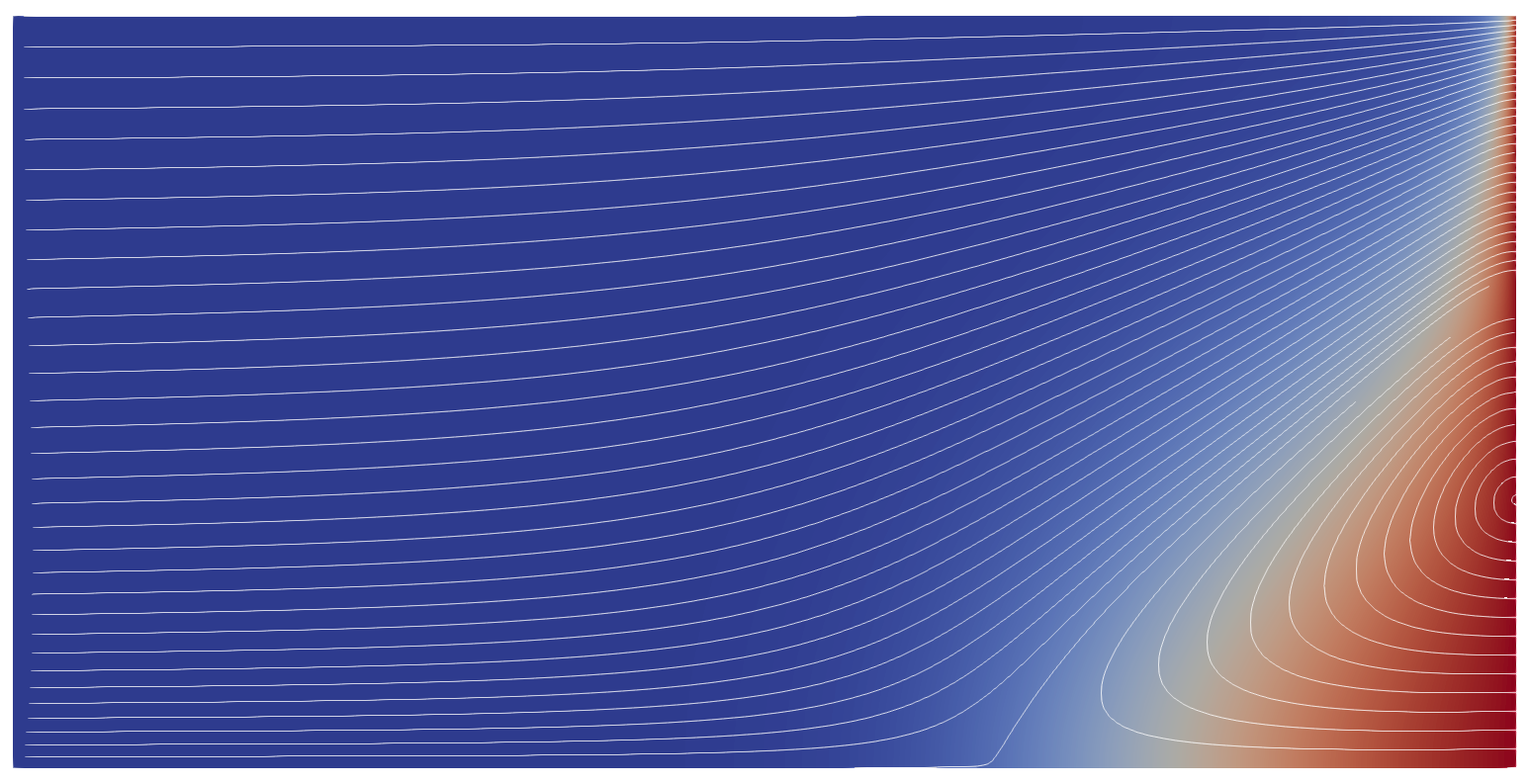}
    \caption{(left) Computational domain $\mathcal{D}:=[0,2]\times [-1,0]$. (Right) One realization of the mass fraction $\conc(t,\bx)$ and the streamlines of the velocity field $\dvel$ for the undisturbed Henry problem at $t = 6016$ $\mathrm{s}$.}
    \label{fig:Henry2d-scheme}
\end{center}    
\end{figure}

A straightforward but very demonstrative model of coastal aquifers is the so-called Henry problem, first considered in \cite{henry1964effects}. In this two-dimensional setting, the aquifer is represented by a rectangular domain $\D = [0, 2] \times [-1, 0]$ $[\mathrm{m}^2]$ entirely saturated with the liquid phase (Fig.~\ref{fig:Henry2d-scheme}). The salty seawater intrudes from the right side, and pure water is recharged from the left. The top and bottom are considered impermeable. Analogous settings with partially saturated domains are considered in \cite{Stoeckl}.

The mass conservation laws for the entire liquid phase and salt yield the following equations
\begin{eqnarray}
 \label {e_cont_eq}
 \partial_t (\poro \dens) & + & \nabla \cdot (\dens \dvel) = 0, \\
 \label {e_tran_eq}
 \partial_t (\poro \dens \conc) & + & \nabla \cdot (\dens \conc \dvel - \dens \disp \nabla \conc) = 0,
\end{eqnarray}
where $\poro: \D \to \mathbb{R}$ denotes the porosity, $\perm: \D \to \mathbb{R}^{d \times d}$ represents the permeability, $\conc (t, \mathbf{x}): [0, +\infty) \times \D \to [0, 1]$ is the mass fraction of the salt (or of the brine) in the solution, $\dens = \dens (\conc)$ indicates the density of the liquid phase, and $\disp (t, \mathbf{x}): [0, +\infty) \times \D \to \mathbb{R}^{d \times d}$ denotes the molecular diffusion and mechanical dispersion tensor. For the velocity $\dvel (t, \mathbf{x}): [0, +\infty) \times \D \to \mathbb{R}^d$, we assume Darcy's law:
\begin{eqnarray} \label {e_Darcy_vel}
 \dvel = - \frac{\perm}{\visc} (\nabla \pres - \dens \grav),
\end{eqnarray}
where $\pres = \pres (t, \mathbf{x}): [0, +\infty) \times \D \to \mathbb{R}$ is the hydrostatic pressure, $\visc = \visc (\conc)$ denotes the viscosity of the liquid phase, and $\grav = (0, \dots, 0, - g)^T \in \mathbb{R}^d$ represents the gravity vector. Inserting (\ref {e_Darcy_vel}) into (\ref {e_cont_eq}--\ref {e_tran_eq}) results in a system of two time-dependent PDEs in the unknowns $\conc$ and $\pres$. This system should be closed with boundary conditions for $\conc$ and $\pres$ and an initial condition for $\conc$.

\begin{table}[b]
\begin{center}

 \begin{tabular}{|l|l|l|} \hline
  Parameter & Values and Units & Description \\ \hline
  $\EXP{\phi}$ & 0.35 [-] & mean value of porosity \\ \hline
  $D$ & $18.8571\cdot 10^{-6}$ [$\mathrm{m}^2 \cdot \mathrm{s}^{-1}$] & diffusion coefficient in the medium \\ \hline
  $\perm$ & $1.020408\cdot 10^{-9}$ [$\mathrm{m}^2$] & permeability of the medium \\ \hline
%  $E(\perm)$ & {\color{red} ???} [$\mathrm{m}^2$] & the mean value of the permeability \\ \hline
  $g$ & $9.8$ [$\mathrm{m} \cdot \mathrm{s}^{-2}$] & gravity \\ \hline
  $\rho_0$ & $1000$ [$\mathrm{kg} \cdot \mathrm{m}^{-3}$] & density of pure water \\ \hline
  $\rho_1$ & $1024.99$ [$\mathrm{kg} \cdot \mathrm{m}^{-3}$] & density of brine \\ \hline
  $\mu$ & $10^{-3}$ [$\mathrm{kg} \cdot \mathrm{m}^{-1} \cdot \mathrm{s}^{-1}$] & viscosity \\ \hline
 \end{tabular}
 \caption{Parameters of the considered density-driven flow problem}
 \label{tab:HenryParam}
\end{center}
\end{table}

Following the classical setting in \cite{henry1964effects}, for this variant of the Henry problem, we set
\begin {eqnarray} \label {e_lin_density}
 \dens (\conc) = \dens_0 + (\dens_1 - \dens_0) \conc, \qquad  \qquad \visc = \text{const}
\end {eqnarray}
and
\begin {eqnarray} \label {e_mol_diff}
 \disp = \poro D \mathbf{I}
\end {eqnarray}
with a constant scalar $D \in \mathbb{R}$, and the identity matrix $\mathbf{I}\in\mathbb{R}^{d\times d}$. Furthermore, we assume the isotropic permeability
\begin {eqnarray*} %\label {e_perm_of_poro}
 \perm = K \mathbf{I}, \qquad K \in \mathbb{R}.
\end {eqnarray*}
This setting is consistent with the problem setting in \cite{Voss_Souza}. However, we do not assume the Boussinesq approximation and keep the density variable for all terms. For the initial conditions, we set
\begin {equation}
 \left. \conc \right |_{t = 0} = 0.
\end {equation}

The boundary conditions are presented in Fig.~\ref {fig:Henry2d-scheme}(a). On the right side of the domain, we impose Dirichlet conditions for the $\conc$ and $\pres$ variables that represent the adjacent seawater aquifer:
\begin {equation}
	\left. \conc \right |_{x=2} = 1, \qquad \left . \pres \right |_{x=2} = - \rho_1 g y.
\end {equation}
On the left side, we prescribe the inflow of fresh water:
\begin {equation}
	\left. \conc \right |_{x=0} = 0, \qquad \left . \dens \dvel \cdot \mathbf{e}_x \right |_{x=0} = \hat{q}_{\mathrm{in}},
\end {equation}
where $\mathbf{e}_x = (1, 0)^\top$, and $\hat{q}_{\mathrm{in}}$ is a constant. For the classical formulation of the Henry problem, this value was set to $\hat{q}_{\mathrm{in}} = 6.6 \cdot 10^{-2}$ $\mathrm{kg}/\mathrm{s}$ in \cite{Voss_Souza} or $\hat{q}_{\mathrm{in}} = 3.3 \cdot 10^{-2}$ $\mathrm{kg}/\mathrm{s}$ in \cite{Simpson04_Henry,Simpson2003}. The Neumann zero boundary conditions are imposed on the upper and lower sides of $\mathcal{D}$.

An example of $\conc(t,\bx)$ and $\dvel(t,\bx)$ for the parameters from Table~\ref {tab:HenryParam} is presented in Fig.~\ref {fig:Henry2d-scheme}(right). The dark red color corresponds to $\conc=1$, and dark blue corresponds to $\conc=0$. Due to its higher density, the saltwater intrudes into the aquifer in the lower right part. It is pushed back by the lighter pure water coming from the left. This process induces a vortex in the flow in the lower right corner of the domain. The saltwater flows in at the lower part of the right boundary and deviates to the top and right, back to the seaside, forming a salt triangle. This flow does not transport the salt to the left part of the domain. The salt propagates further to the left due to diffusion and dispersion and is washed out by the recharge. In the classical formulation, this salt triangle initially increases over time but achieves a steady state (cf. \cite {Voss_Souza,Simpson04_Henry,Simpson2003}). However, the initial nonstationary phase may take significant time. Investigating this phase is especially important to understand the system behavior when changing the recharge. For this, in addition to the mean and variance, we consider the mass fraction at 12 points (listed below) and an integral value---the total amount of pure water (as in Eq.~\ref{eq:integral_fw}). The list of chosen points follows:
\begin{align}
%\{p_1,p_2,\ldots,p_{12}\}:
\{(x,y)_{i=1,\ldots,12}\}&=\{(1.10, -0.95),
	(1.35, -0.95),
	(1.60, -0.95),
	(1.85, -0.95),
	(1.10, -0.75),
	(1.35, -0.75), \label{eq:12points}\\ 
	&(1.60, -0.75),
	(1.85, -0.75),
	(1.10, -0.50),
	(1.35, -0.50),
	(1.60, -0.50),
	(1.85, -0.50).\} \nonumber 
\end{align}
%One could also consider other points.
The motivation is to consider points where the concentration variation is considerable.
In addition, the mass fraction $\conc$ at each point $\bx$ is a function of time. 

These spatial points may help track salinity changes over time in groundwater wells and understand which areas in the aquifer are most vulnerable. Farmers can use this information to take action, such as decreasing salinity or adapting strategies by planting salt-tolerant crops.

%!!!
%
%
\subsection{Modeling porosity, permeability, and mass fraction}
\label{subsec:PorosityVar}
The primary sources of uncertainty are the hydrogeological properties of the porous medium---porosity ($\poro$) and permeability ($\perm$) fields of the solid phase---and the freshwater recharge flux $\hat{q}_x$ through the left boundary. 
%Other uncertainties are not considered in this work.    
The QoIs are related to the mass fraction $\sol$, a function of $\poro$, $\perm$, and the recharge. We model the uncertain $\poro$ using a random field and assume $\perm$ to be isotropic and dependent on $\poro$:
\begin {eqnarray} \label {e_perm_of_poro}
 \perm = K \mathbf{I}, \qquad K = K (\poro) \in \mathbb{R}.
\end {eqnarray}
The distribution of $\poro(\bx,\xib)$, $\bx\in \D$, is determined by a set of stochastic parameters $\xib=(\xi_1,\ldots,\xi_M,...)$. Each component $\xi_i$ is a random variable depending on a random event $\omega$. For concision, we skip $\omega$ and write $\xib:=\xib(\omega)$.

The dependence in Eq.~(\ref{e_perm_of_poro}) is specific for every material.
%and there is no a general law. 
We refer to \cite{Panda_Lake_Perm_vs_Por,Pape_Clauser_Iffland_1999,Costa_2006} for a detailed discussion. In the proposed model, we use a Kozeny--Carman-like dependence
\begin{eqnarray} 
\label{e_perm_Kozeny_Carman}
 K (\poro) = \kappa_{KC} \cdot \dfrac {\poro^3} {1 - \poro^2},
\end {eqnarray}
where the scaling factor $\kappa_{KC}$ is chosen to satisfy the equality $K(\EXP{\poro}) \mathbf{I} = E(\mathbf{K})$, resembling the parameters of the standard Henry problem.
% \begin{remark}
% Typical porosity values, which we have found in the literature are: well sorted sand $0.25-0.5$, poorly sorted sand $0.15-0.3$, clay $0.4-0.6$, crystalline rock $0.001-0.01$, fractured rock $0.01-0.05$.
% \end{remark}
The inflow flux is kept constant across the left boundary but depends on the stochastic variable $q_{\mathrm{in}}$. We also assume that the inflow flux is independent of $\poro$ and $\perm$.

% In the experiments \hat{q}_x = \hat{q}_x (\theta_3) = -6.6e-2 * (1 + 0.5 * \theta_3ls)
%A discretisation of these uncertain parameters see in Sect.~\ref{sec:Discret}.
%
%
%\subsection{Quantities of interest}
%\label{sec:QoI}
%$P(c>c^*)>\varepsilon$, where $c^*$ is some critical value, $0<\varepsilon<1$
%
%

\subsection{Numerical methods for the deterministic problem}
\label{sec:Num}
The system (\ref{e_cont_eq}--\ref{e_tran_eq}) is numerically solved in the domain $\D \times [0, T]$, where the symbol $\times$ denotes the Cartesian product. After the discretization of $\D$ using quadrilaterals of size $h$, we obtain $\D_h$. Equations (\ref{e_cont_eq}--\ref{e_tran_eq}) are discretized using a vertex-centered finite-volume scheme with a “consistent velocity” for the approximation of Darcy's law \eqref{e_Darcy_vel}, as presented in \cite{Frolkovic-DeSchepper-ConvDom,Frolkovic-ConsVel,Frolkovic-Knaber-ConsVel}. The degrees of freedom associated with $\D_h$ are denoted by $n$. There are two degrees of freedom per grid vertex: one for the mass fraction and another for the pressure. We use the implicit Euler method with a fixed time step $\tau$ for time discretization. The number of the computed time steps is $r = T / \tau$.
% with a fixed time step length to simplify evaluations of the stochastic quantities.

We use partial upwind for the convective terms (cf.\ \cite{Frolkovic-DeSchepper-ConvDom}). 
Therefore, the discretization error is of the second order w.r.t. the spatial mesh size $h$. However, the diffusion in \eqref{e_tran_eq} is minimal compared with the velocity. For the grids in the numerical experiments, the observed reduction of the discretization error after grid refinement corresponds to the first order. Thus, we assume the first-order dependence of the discretization error w.r.t. $h$, which is consistent with the numerical experiments. Furthermore, the Euler method provides the first-order dependence of the discretization error w.r.t. $\tau$.

The implicit time-stepping scheme provides unconditional stability but requires the solution to an extensive nonlinear algebraic system of the discretized equations with $n$ unknowns in every time step. The Newton method is used to solve this system. Linear systems inside the Newton iteration are solved using the BiCGStab method (cf.\ \cite {Templates}) preconditioned with the geometric multigrid method (V-cycle, cf.\ \cite{Hackbusch85}). In the multigrid cycle, the ILU${}_\beta$-smoothers \cite{Hackbusch_Iter_Sol} and Gaussian elimination are used as the coarse grid solver.

To construct the spatial grid hierarchy $\D_0, \D_1, \dots, \D_L$, we start with a coarse grid consisting of 512 grid elements (quadrilaterals) and $n_0 = 1122$ degrees of freedom. This grid is regularly refined to obtain all other grid levels. After every spatial grid refinement, the number of grid elements is multiplied by a factor of four. Consequently, the number of degrees of freedom is increased by a factor of four (i.e., $n_\ell \approx n_0\cdot 2^{d\ell}$, $d=2$; see Table~\ref{tab:adaptiveTS_times}). This hierarchy is used in the geometric multigrid preconditioner and MLMC method. We also construct the temporal grid hierarchy $\Tau_0, \Tau_1, \dots, \Tau_L$. The time step on each temporal grid is denoted by $\tau_\ell$ with $\tau_{\ell+1} = \tfrac{1}{2} \tau_\ell$. The number of time steps on the $\ell$th grid (level) is $r_{\ell+1} = 2 r_\ell$ and $r_\ell = r_0 2^\ell$, where $r_0$ is the number of grid points on $\Tau_0$. On the $\ell$th level, the MLMC uses the grid $\D_{\ell}\times \Tau_{\ell}$. Up to six spatial and time grids were used in the numerical experiments.

In the context of this work, it is critical to estimate the numerical complexity of the deterministic solver w.r.t. the grid level $\ell$. The most time-consuming part of the simulation is the solution of the discretized nonlinear system. Typically, it is challenging to predict the number of Newton iterations in every time step, but in the numerical experiments, two iterations were sufficient to achieve the prescribed accuracy. Accordingly, the linear solver was called at most two times per time step. Furthermore, the convergence rate of the geometric multigrid method does not depend on the mesh size (cf.\ \cite{Hackbusch_Iter_Sol}). Hence, the computation complexity of one time step is $\mathcal{O} (n_{\ell})$, where $n_{\ell}$ is the number of the degrees of freedom on the grid level $\ell$. Therefore, the overall numerical cost of the computation of one scenario on grid level $\ell$ for $r_{\ell}$ time steps is
\begin{equation}
\label{eq:CompComplexity}
s_\ell = \mathcal{O} (n_\ell r_\ell), \quad s_\ell \propto s_{\ell-1} 2^{(d+1)}, \quad d=2.
\end{equation}

%\textbf{Software.}  
%ug4 is a flexible software system for simulating PDE based models on high performance parallel clusters \cite{ug4_ref1_2013,ug4_ref2_2013}. Starting in the early 1990's the software toolbox UG (UG states for §Unstructured Grids) has been developed. This software framework has successfully been applied to a variety of problems such as density-driven and thermohaline flow in porous media, Navier-Stokes equation, drug diffusion through human skin, level-set methods.
%Special care is taken to parallel performance issues, scalability, linear algebra algorithms and data structures. In addition, ug4 is specially designed for a user-friendly handling and functionality is made available to non-programming users by the usage of scripts and graphical representations.
%The numerical solution of the systems (\ref {e_cont_eq}--\ref {e_tran_eq}) has been performed using the D3F plugin of the simulation software package ug4, cf.\ \cite{ug4_ref1_2013,ug4_ref2_2013}.

%\textbf{Parallel computations.}
%Computations were conducted in ug4 solver on a parallel machine Shaheen II, provided by the King Abdullah University of Science and Technology.
%

%The data and the code are available on github %repository\footnote{\url{https://github.com/ug4/ughub.wiki.git}}.  

%The number of spatial degrees of freedom in the presented experiments varied from 500 till $500.000$. The maximal number of parallel processing units was $600\times 32$, where $600$ is the number of parallel nodes, and $32$ is the number of computing cores on each node.

%
\section{Multilevel Monte Carlo}
\label{sec:MLMC}
Various numerical methods can quantify uncertainty, and every method has pros and cons. For example, the classical MC method converges slowly and requires numerous samples. To reduce the total computing cost, we apply the MLMC method, which is a natural idea because the deterministic solver uses a multigrid method (see Section~\ref{sec:Num}). The MLMC method efficiently combines samples from various levels. Further, we repeat the main idea of the MLMC method. A more in-depth description of these techniques is found in~\cite{MLMC_PDE_anal11,CMLMC,giles2008,giles2015,ErikOptGeom15,teckentrup2013further,Litv_Scattered19}.

%Before to apply the MLMC, we implement and check some preliminary conditions needed for successful MLMC performance.

We let $\xib(\omega)$ and $g(\xib)=g(\xib(\omega))$ represent a vector of random variables and the QoI, respectively, where $\omega$ is a random event. The MLMC method aims to approximate the expected value $\EXP{g}$ with an optimal computational cost. In this work, $g$ could be $\sol(t,\bx,\xib)$ in the whole domain or at a point $(t,\bx)$ or an integral over a subdomain. The MLMC method constructs a telescoping sum, defined over a sequence of spatial and temporal meshes, $\ell=0, \ldots, L$, as described next, to achieve this goal. Moreover, $g$, numerically evaluated on level $\ell$, is denoted by $g_{h_{\ell},\tau_{\ell},\ell}$ or, for simplicity, by just $g_\ell$, where $h_{\ell}$ and $\tau_{\ell}$ are the discretization steps in space and time on level $\ell$. Further, we assume that $ \EXP{g_{h,\tau}} \rightarrow \EXP{g}$ as $h\rightarrow 0$ and $\tau \rightarrow 0$.  
%\textcolor{red}{From now on, when we write $\EXP{g}$ or $\EXP{g_{\ell}}$ we mean $\EXP{g_{h,\tau}}$ or $\EXP{g_{h_{\ell},\tau_{\ell},\ell}}$ respectively.}

Furthermore, $s_0$ is the computing cost to evaluate one realization of $g_0$ (the most expensive one from all realizations). Similarly, $s_\ell$ denotes the computing cost of evaluating $g_\ell - g_{\ell-1}$. For simplicity, we assume that $s_\ell$ for $g_{\ell} - g_{\ell-1}$ is almost the same as $s_\ell$ for $g_{\ell}$. The number of iterations is variable; thus, the cost of computing a sample of $g_\ell - g_{\ell-1}$ may fluctuate for various realizations.

For a better understanding, we consider a two-level MLMC (cf. \cite{giles2015}) and estimate the optimal number of needed samples on both levels. The two-level MLMC has only two meshes: a coarse one and a fine one. The QoI $\EXP{g}$ can be approximated on the fine mesh by $\EXP{g_1}$ and on the coarse mesh by $\EXP{g_0}$. Furthermore,
\begin{equation} \label {eq:2levelMC}
\EXP{g_1} = \EXP{g_0} + \EXP{g_1 - g_0}\approx m_0^{-1}\sum_{i=1}^{m_0} g_0^{(i)} + m_1^{-1}\sum_{j=1}^{m_1} (g_1^{(j)} - g_0^{(j)}),
\end{equation}
where $g_1^{(j)} - g_0^{(j)}:= g_1(\xib_j)- g_0(\xib_j)$, $\xib_j$ is a random vector, and $m_0$ and $m_1$ represent the numbers of quadrature points (numbers of samples/realizations) on the coarse and fine meshes, respectively. 
% Let us define the cost of computing a single sample of $g_0$ and $g_1−g_0$ by $s_0$ and $s_1$, respectively.
The total computational cost of evaluation \eqref{eq:2levelMC} is $m_0 s_0+m_1 s_1$. The variances of $g_0$ and $g_1−g_0$ are denoted by $V_0$ and $V_1$, and the total variance is obtained by $V_0/m_0 + V_1/m_1$, assuming that $g_0^{(i)}$ and $g_1^{(j)}-g_0^{(j)}$ use independent samples. By solving an auxiliary minimization problem, the variance is minimal if $m_1 = m_0\cdot \frac{\sqrt{V_1 /s_1}}{\sqrt{V_0 /s_0}}$. Thus, with the estimates of the variances and $m_0$, we can estimate $m_1$. 
%Our assumption here is that $V_1$ is small.

The idea presented above can be extended to a case with multiple levels. Thus, we can find (quasi-) optimal numbers of samples $m_0,m_1,\ldots, m_L$.
%Let $\{{P_\ell}\}_{\ell=0}^L$ be a sequence of spatio-temporal meshes that discretize the computational domain.
%The corresponding grid step sizes denote by $\{h_{\ell}\}_{\ell=0}^L$ 
%It is assumed $\{{P_\ell}\}_{\ell=0}^L$ are generated hierarchically with ${h_\ell=h_0 \beta^{-\ell}}$ for $h_0>0$ and a constant  $\beta > 1$.
The MLMC method calculates $\EXP{g_L}\approx \EXP{g}$ using the following telescopic sum: 
\begin{align}
  \EXP{g_L} &= \EXP{g_{0}} + \sum_{\ell=1}^L \EXP{g_{\ell}-g_{\ell-1}} \label{eq:EgL} \\
  &\approx  %= \sum_{\ell=0}^L \EXP{G_{\ell}}
  m_0^{-1}\sum_{i=1}^{m_0} g_{0}^{(0,i)}  + \sum_{\ell=1}^L \left( m_\ell^{-1}\sum_{i=1}^{m_\ell} (g_{\ell}^{(\ell,i)} - g_{\ell-1}^{(\ell,i)} )\right). \ \label{eq:A}
\end{align}
In the above equation, level $\ell$ in the superscript $(\ell, i)$ indicates that independent samples are used at each correction level.
%where $ G_\ell :=g_\ell- g_{\ell-1}$, $ G_0=g_0$.
% \begin{equation}
% \label{eq:EgLapprox}
%   \EXP{g_L} = \EXP{g_{0}} + \sum_{\ell=1}^L \EXP{g_{\ell}-g_{\ell-1}}, %= \sum_{\ell=0}^L \EXP{G_{\ell}}
% \end{equation}
% Note that both $g_\ell$ and $g_{\ell-1}$ are computed using the same input random parameter $\xib$.
% In the telescoping sum~\eqref{eq:EgL}, the expected values in practice are replaced by sample averages, i.e., \begin{equation}
% \EXP{G_\ell} \approx \est G_\ell  =  M_\ell^{-1}\sum_{m=1}^{M_\ell} G_{\ell,m},
% \end{equation}
% where random variable $G_{\ell,m}$ have the same distribution as  %\thicksim 
% $G_{\ell}$ and are independent identically distributed (i.i.d.) samples.
As $\ell$ increases, the variance of $g_{\ell} - g_{\ell-1}$ decreases. Thus, the total computational cost can be reduced by taking fewer samples on finer meshes.

We recall that $h_{\ell}=h_0\cdot 2^{-2\ell}$ and $\tau_{\ell}=\tau_0\cdot 2^{-\ell}$. We assume that the average cost of generating one sample of $g_{\ell}$ (the cost of one deterministic simulation for one random realization) is
\begin{equation}
s_\ell = \mathcal{O}(n_{\ell}r_{\ell})=\mathcal{O}(h_{\ell}^{-1} \tau_{\ell}^{-1})=
\mathcal{O}\left( \frac{1}{h_0 \tau_0} 2^{2\ell}2^{\ell} \right)=
\mathcal{O}\left( \frac{1}{h_0 \tau_0} 2^{3\ell} \right)
=\mathcal{O}\left( \frac{1}{h_0 \tau_0}  2^{(d+1)\ell \gamma} \right),
\label{eq:workpl}
\end{equation}
where $d=2$ is the spatial dimension, and $\gamma=1$ is determined by the computational complexity of the deterministic solver (ug4).

%For the solver used here we have $d=2$, $\gamma\approx 1$, and $\beta=2$.

We let $V_{\ell}$ be the variance of one sample of $g_{\ell}-g_{\ell-1}$. Then,
the total cost and variance of the multilevel estimator in Eq.~(\ref{eq:EgL}) are
$\sum_{\ell=0}^{L} m_{\ell}s_{\ell}$ and $\sum_{\ell=0}^{L}V_{\ell}/m_{\ell}$, respectively. For a fixed variance, the cost is minimized by choosing $m_{\ell}$ to minimize
the following functional for some value of the Lagrange multiplier $\mu^2$:
\begin{equation}
\label{eq:goal_function}
F(m_0,\ldots,m_{L}):=\sum_{\ell=0}^{L} m_{\ell}s_{\ell}+\mu^2 \frac{V_{\ell}}{m_{\ell}}.
\end{equation}
To determine $m_{\ell}$, we take the derivatives w.r.t. $m_{\ell}$ and set them equal to zero:
$$\frac{\partial F(m_0,\ldots,m_{L})}{\partial m_{\ell}}:= s_{\ell}-\mu^2 \frac{V_{\ell}}{m_{\ell}^2}=0.$$
After solving the obtained equations, we obtain
$$m_{\ell}^2=\mu^2 \frac{V_{\ell}}{s_{\ell}} \quad \text{and} \quad
m_{\ell}=\mu \sqrt{\frac{V_{\ell}}{s_{\ell}}}.$$
To achieve an overall variance of $\varepsilon^2$, that is,
$$\sum_{\ell=0}^{L}V_{\ell}/m_{\ell} = \varepsilon^2,$$
we substitute $m_{\ell}$ with the computed $m_{\ell}=\mu \sqrt{\frac{V_{\ell}}{s_{\ell}}}$, and obtain
$$\sum_{\ell=0}^{L}\frac{V_{\ell}}{\mu \sqrt{\frac{V_{\ell}}{s_{\ell}}}} = \varepsilon^2.$$
From the last equation, we obtain
$$\mu=\varepsilon^{-2}  \sum_{\ell=0}^{L} \sqrt{V_{\ell}s_{\ell}},\quad \text{and}$$
\begin{equation}
\label{eq:M_ell}
m_{\ell}=\varepsilon^{-2} \sqrt{\frac{V_{\ell}}{s_{\ell}}} \sum_{i=0}^{L} \sqrt{V_{i}s_{i}}.
\end{equation}
The total computational cost is $S:=\varepsilon^{-2}\left( \sum_{\ell=0}^L \sqrt{V_{\ell} s_{\ell}}\right)^2$ (for further analysis of this sum, see \cite{giles2015}, p.4). %The Theorem~1 in \cite{giles2015} lists the conditions when MLMC is faster than the standard MC. In the worst case, when the dominant computational cost is on the coarsest level, the MLMC has the same computational cost as the standard MC.

\begin{defn}
We let 
\begin{align}
\label{eq:Yell} 
\EXP{Y_{\ell}}:=
      \begin{cases}
        \EXP{g_0}, & \ell=0 \\
        \EXP{g_{\ell} - g_{\ell-1}}, & \ell>0 \\
      \end{cases}.
\end{align}
In addition, $Y:=\sum_{\ell=0}^L Y_{\ell}$ denotes a multilevel estimator of $\EXP{g}$ based on $L+1$ levels and $m_{\ell}$ independent samples on level $\ell$, where $\ell=0,\ldots,L$. Moreover, $Y_{\ell}=m_{\ell}^{-1}\sum_{i=1}^{m_{\ell}} (g_{\ell}^{(\ell,i)} - g_{\ell-1}^{(\ell,i)})$, where $g_{-1}\equiv 0$. 

The standard theory indicates that $\EXP{Y}=\EXP{g_L}$, $\Var{Y_{\ell}}=\sum_{\ell=0}^L m_{\ell}^{-1} V_{\ell}$, and $V_{\ell}\equiv \Var{g_{\ell} - g_{\ell-1}}$.
\end{defn}

%!!!
The mean squared error (MSE) is used to measure the quality of the multilevel estimator:
\begin{equation}
\label{eq:MSE}
\MSE:=\EXP{(Y-\EXP{g})^2}=\Var{Y} + \left( \EXP{Y} - \EXP{g} \right)^2.    
\end{equation}
To obtain an MSE smaller than $\varepsilon^2$, we ensure that both $\Var{Y}$ and $\left( \EXP{Y} - \EXP{g} \right)^2=(\EXP{g_L-g})^2$ are smaller than $\varepsilon^2/2$. Combining this idea with a geometric sequence of levels in which the cost increases exponentially with the level while the weak error $\EXP{g_L−g}$ and multilevel correction variance $V_{\ell}$ decrease exponentially leads to the following theorem (cf. Theorem 1, p.~6 in \cite{giles2015}):
\begin{theorem}
\label{thm:costMLMC}
We let $d$ denote the problem dimension. Suppose positive constants $\alpha,\beta,\gamma > 0$ exist such that $\alpha \geq \frac{1}{2} \min(\beta, \gamma d)$, and
\begin{subequations}
\label{eq:q1q2}
\begin{align}
    \vert \EXP{g_\ell-g} \vert & \leq c_1 2^{-\alpha\ell} \label{eq:weak_error_model} \\
    V_{\ell} & \leq c_2 2^{-\beta\ell} \label{eq:strong_error_model} \\
     S_{\ell} &\leq c_3 2^{d\gamma \ell}.
\end{align}
\end{subequations}
Then, for any accuracy $\varepsilon < e^{-1}$, a constant $c_4>0$ and a sequence of realizations $\{m_{\ell}\}_{\ell=0}^L$ exist, such that
%the following statements are satisfied [cf.~\eqref{eq:goal}]
%$e(\tilde{G}_{\ell})< \tol$
\begin{equation*}
\MSE:=\EXP{(Y-\EXP{g})^2}< \varepsilon^2,    
\end{equation*}
and the computational cost is
\begin{align}
\label{eq:mlmc_iso_work} 
%S_\varepsilon\left(\hat{Q}^{ML}_{h,\{m_{\ell}\}}\right)\leq
S=
      \begin{cases}
        \mathcal{O}({\varepsilon^{-2}}), & \beta > d\gamma \\
        \mathcal{O}({\varepsilon^{-2}) \left(\log(\varepsilon)\right)^2}, & \beta= d\gamma \\
        \mathcal{O}({\varepsilon^{-\left(2 +\frac{d\gamma-\beta}{\alpha}\right)}}),  & \beta < d\gamma. \\
      \end{cases}
\end{align}
\end{theorem}

This theorem (see also \cite{hoel2014implementation,hoel2012adaptive,charrier2013,MLMC_PDE_anal11,giles2008}) indicates that, even in the worst-case scenario, the MLMC algorithm has a lower computational cost than that of the traditional (single-level) MC method, which scales as $\mathcal{O}(\varepsilon ^{-2-d\gamma/\alpha})$. %\textcolor{red}{DELETE? whereas each case in \eqref{eq:mlmc_iso_work} shows a smaller total work}.
Furthermore, in the best-case scenario presented above, the computational cost of the MLMC algorithm scales as $\Order{\varepsilon ^{-2}}$.
%, i.e. identical to that of the MC method assuming that the cost per sample is $\Order{1}.$ 
%In other words, for this case, the \CMLMC algorithm can in effect remove the computational cost required by the discretization, namely  $\mathcal{O}(\tol^{-d\gamma/q_1})$.
%

Using preliminary tests, we can estimate the convergence rates $\alpha$ for the mean (the so-called weak convergence) and $\beta$ for the variance (the so-called strong convergence).
% \begin{subequations}
% \label{eq:q1q2}
% \begin{align}
%     \EXP{g_\ell-g_{\ell-1}} & = \mathcal{O}(h_\ell^{q_1}) \label{eq:weak_error_model} \\
%     \var{g_\ell-g_{\ell-1}} &= \mathcal{O}( h_{\ell-1}^{q_2}) \label{eq:strong_error_model} 
% \end{align}
% \end{subequations}
In addition, $\alpha$ is strongly connected to the order of the discretization error (see Section~\ref{sec:Num}), which equals $1$, and precise estimates of parameters $\alpha$ and $\beta$ are crucial to distribute the computational effort optimally.
% The total computational cost of the adaptive algorithm is close to that of the \MLMC method with correct values of parameters given a priori.

%!!
%
%
% \subsection{Probability density functions and exceedance probabilities}
% After expansion in Eq.~\ref{eq:PCEdecom} is calculated, it could be used for computing probability density functions (Fig.~\ref{fig:points}), exceedance probabilities or quantilies (Tab.~\ref{tab:quantiles}). Below we compute \pdf in a point $(x^{*},t^{*})$ and monitor how such pdfs are changing in time.

% The random variable $\hat{\sol}(\thetab)$ can be sampled, e.g., $N_s=10^6$ times (no additional extensive simulations are required). The obtained sample can be used to evaluate the required statistics, and pdf. The exceedance probability for some critical value $c^{*}$ can be estimated as follow:
% \begin{equation}
%     P(\sol>\sol^{*})\approx \frac{\#\{\sol(\thetab_i):\; \sol(\thetab_i)>\sol^{*}, \; i=1,\ldots,N_s\}}{N_s}.
% \end{equation}
% %

%Example ($m=5$, 3 layers, the number of GL quadrature points is $241$.):\\

%
%
\section{Numerical Experiments}
\label{sec:numerics}
The goal is to reduce the total computational cost of stochastic simulations. We use the MLMC method to compute the mean value of various QoIs, such as $\sol$ in the whole domain, $\sol$ at a point, or an integral value (we call it the freshwater integral):
\begin{equation}
\label{eq:integral_fw}
Q_{FW}(t,\omega):=\int_{\bx\in \D} I(\sol(t,\bx,\omega) \le 0.012178) d\bx,
\end{equation}
where $I(\cdot)$ is the indicator function identifying a subdomain $\{\bx:\; \sol(t,\bx,\omega) \le 0.012178\}$, meaning the mass of the fresh water at a time $t$. %The computed QoIs are compared with those computed using the QMC approach. 
Each simulation may contain up to $n=0.5\cdot 10^6$ spatial mesh points and a few thousand time steps ($r=6016$ on the finest mesh).

%
% In the following, we describe a numerical example, where we assume that the porosity coefficient and the recharge are uncertain.
% The mean and the variance of $\sol(t,\bx,\thetab)$ is computed by the MLMC method. 
% The reference solution is estimated by the QMC method (Halton sequence).

\textbf{Software and parallelization:}
The computations presented in this work were performed using the ug4 simulation software toolbox (\url{https://github.com/ug4/ughub.wiki.git}) \cite{ug4_ref1_2013,ug4_ref2_2013}. This software has been applied for subsurface flow simulations of real-world aquifers (cf.~\cite{SWLRHEGW-SaltwaterInNorthSea2018}). The toolbox was parallelized using MPI, and the presented results were obtained on the Shaheen II cluster provided by the King Abdullah University of Science and Technology. Every sample was computed on 32 cores of a separate cluster node. Each simulation (scenario) was localized to one node to reduce the communication time between nodes. All scenarios were concurrently computed on different nodes. A similar approach was used in \cite{LitLog3D_20,Litvinenko-UQ-2021}. Simulations were performed on different meshes; thus, the computation time of each simulation varied over a wide range (see Table~ \ref{tab:adaptiveTS_times}). 
%The nodes where the processes have been completed were immediately released and made available for other users.

\begin{figure}[htbp!]
\begin {center}
\includegraphics[width=0.49\textwidth]{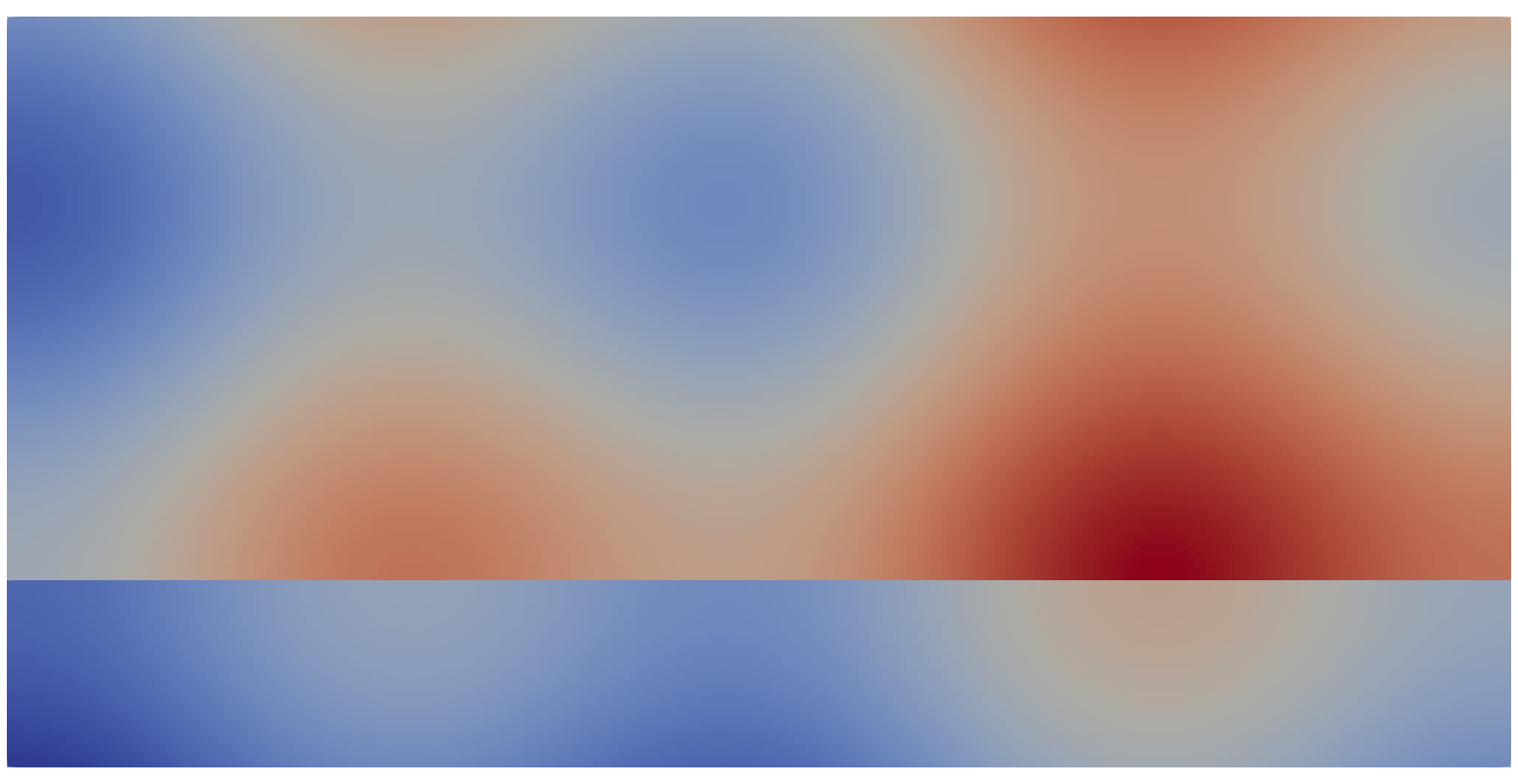}%
\;
\includegraphics[width=0.49\textwidth]{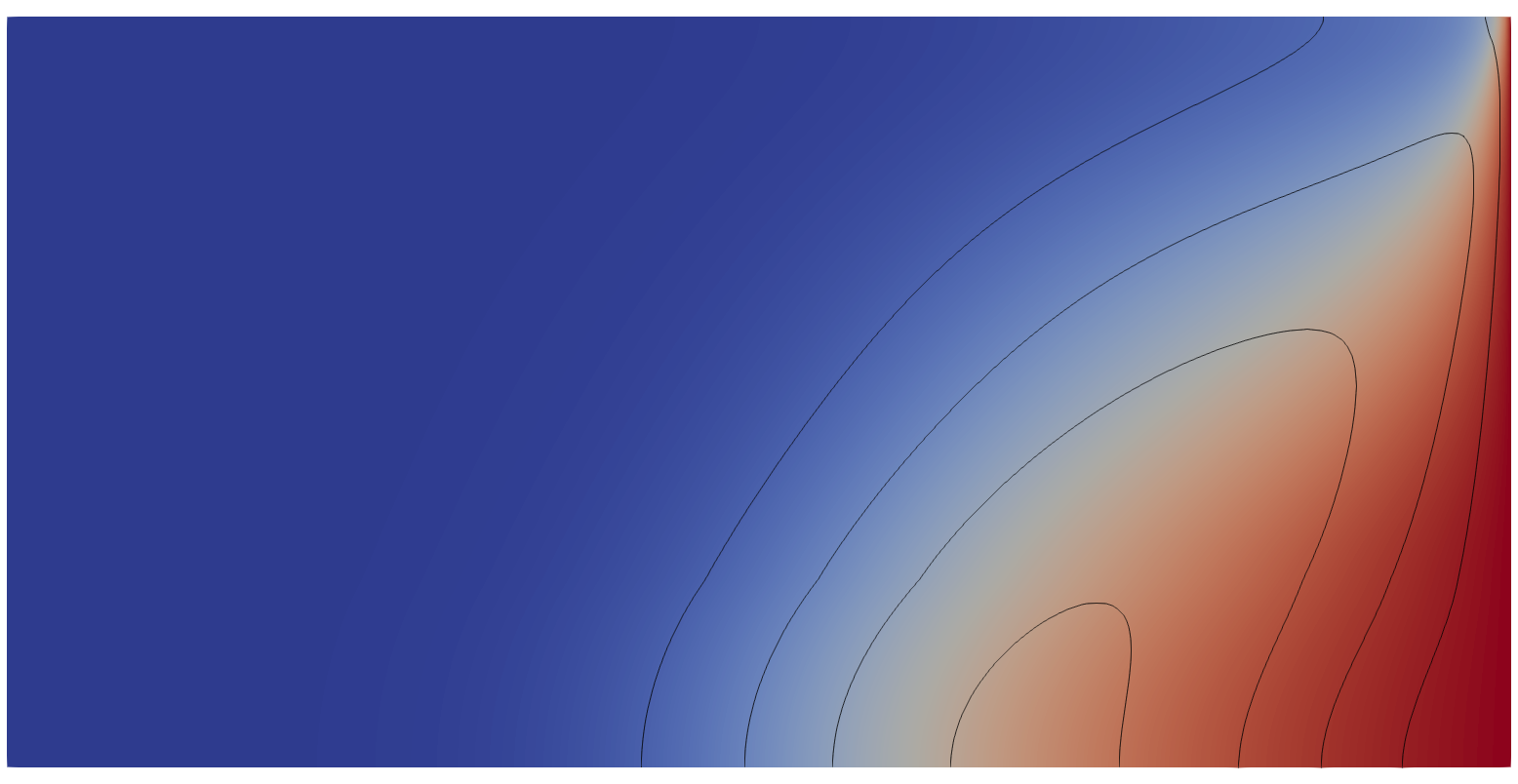}%
\caption {(left) Realisation of porosity $\phi(\xib^*) \in [0.248, 0.499]$. (right) Corresponding mass fraction $\sol(T,\bx,\phi(\xib^*)) \in [0,1]$ with isolines $\{\bx:\; |\sol(T,\poro(\xib^*)) - \overline{\sol}(T)|=0.1\cdot i\}$, $i=1,2,3,4$.}
\label{fig:sol_poro1}
\end{center}
\end{figure}

\textbf{Porosity and recharge:} We assume two horizontal layers: $y\in (-0.75,0]$ (the upper layer) and $y\in [-1, -0.75] $ (the lower layer). The porosity inside each layer is uncertain and is modeled as in Eq.~(\ref{eq:poro_2levels}):
\begin{align}
  \poro(\bx,\xib) &= 0.35\cdot(1+0.15(\xi_2\cos(\pi x/2) + \xi_2 \sin(2\pi y)+\xi_1 \cos(2\pi  x)))\cdot C_0(\xi_1), \label{eq:poro_2levels} \\
  \text{where}\;\;C_0(\xi_1) &= \left\{ 
  \begin{array}{ll}
    1 + 0.2  \xi_1 & \text{if}\; y<-0.75\\
    1 - 0.2  \xi_1 & \text{if}\; y\geq -0.75,
    \end{array} 
  \right. \label{eq:poro_2levels_C}
\end{align}
Additionally, the recharge flux is also uncertain and is equal to
\begin {eqnarray}
 \hat{q}_{\mathrm in} = -6.6\cdot 10^{-2}(1+0.5\cdot \xi_3),
\end {eqnarray}
where $\xi_1$, $\xi_2$, and $\xi_3$ are sampled independently and uniformly in $[-1,1]$.
Figure~\ref{fig:sol_poro1} depicts a random realization of the porosity random field $\poro(\xib)$ (left) and the corresponding solution $\conc(t,\bx,\xib)=\sol(t,\poro(\xib))$ at $t=T$ (right). Additionally, four isolines $\{\bx:\; |\sol(t,\poro(\xib)) - \overline{\sol}(t)|=0.1\cdot i\}$, $i=1,2,3,4$, are presented on the right. The isolines demonstrate the absolute value of the difference between the computed realization $\sol(t,\poro(\xib))$ and the expected value $\overline{\sol}(t)$. These computations were performed for  $\xib=\xib^*= (-0.5898, -0.7257, -0.9616)$ and $t = T=6016$ s. 

\begin{figure}[htbp!]
\begin{center}
  \includegraphics[width=0.49\textwidth]{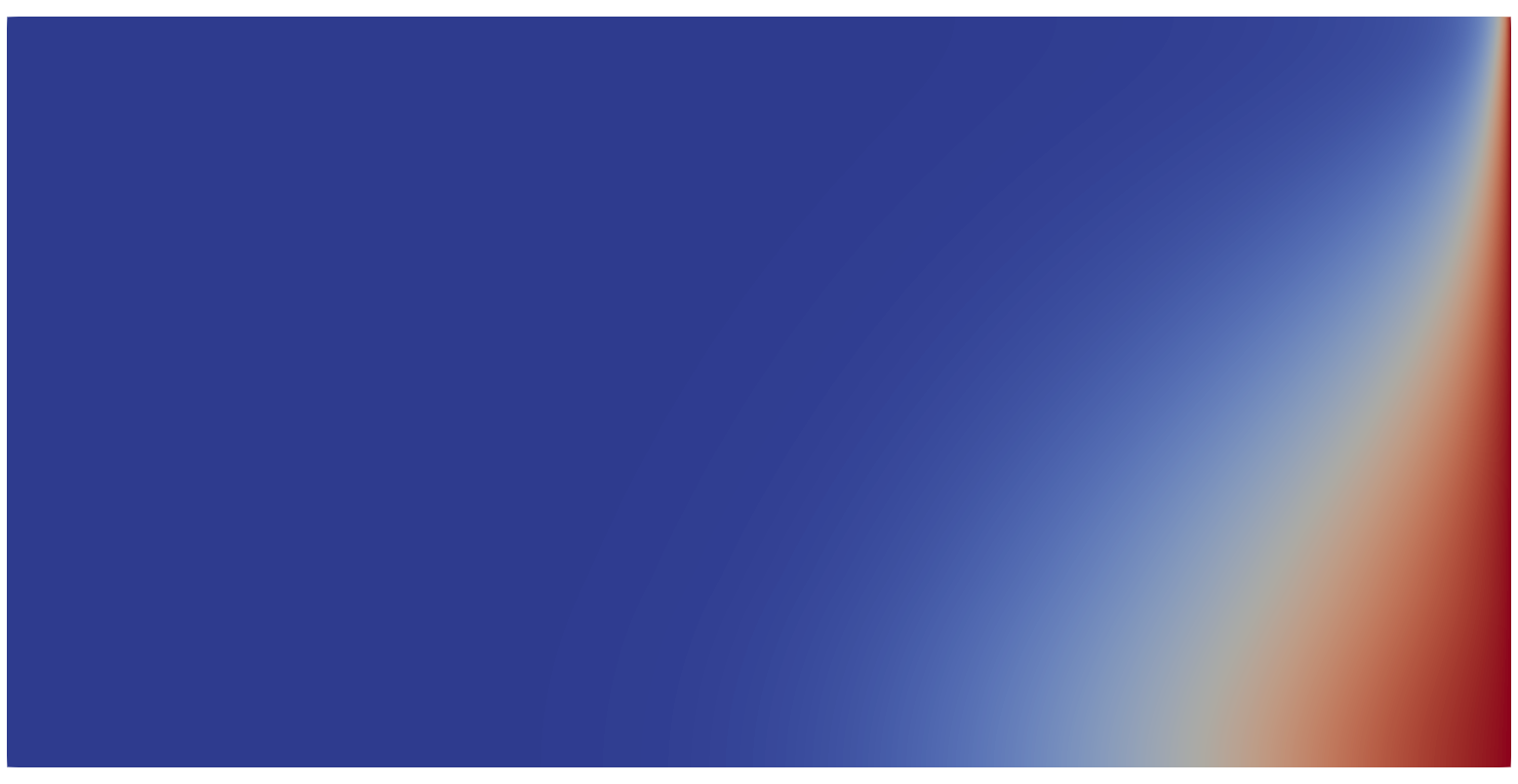}\;
  \includegraphics[width=0.49\textwidth]{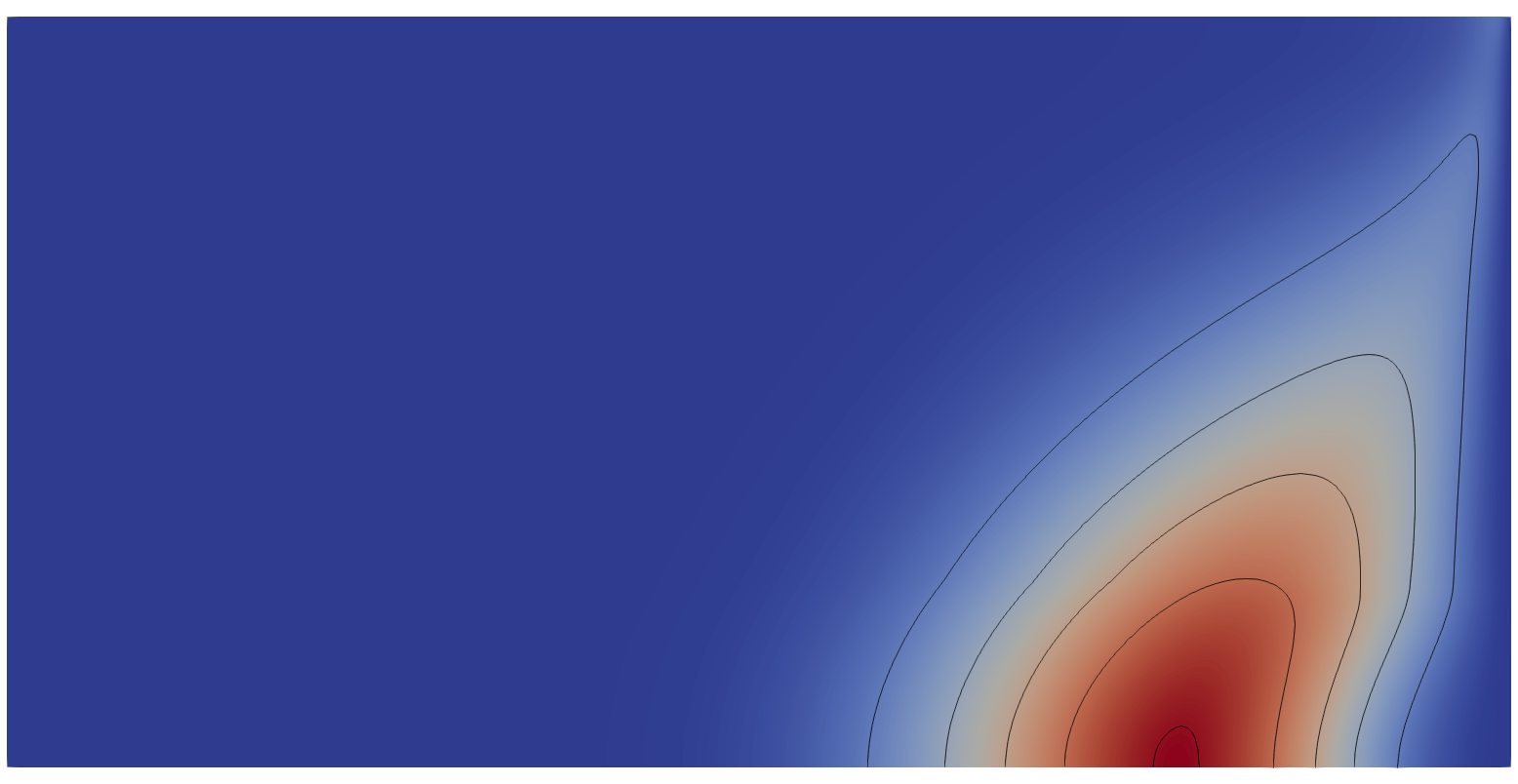}\;
  \caption{(left) Mean value $\overline{\sol} \in [0,1]$ and (right) variance $\var{\sol} \in [0.01,0.05]$ of the mass fraction, with contour lines $\{\bx:\; \var{ \sol }=0.01\cdot i\}$, $i=1..5$, $t=T=6016$.}
    \label{fig:countours_mean_var}
\end{center}    
\end{figure}

The mean and variance of the mass fraction are provided in Fig.~\ref{fig:countours_mean_var} on the left and right, respectively. The expectation takes values from $[0,1]$, and the variance range is $[0,0.05]$. The areas with high variance (dark red) indicate regions with high variability/uncertainty. Such regions may need additional attention from specialists (e.g., placement of additional sensors). Additionally, the right image displays five contour lines
  $\{\bx:\; \var{ \sol }(t,\bx)=0.01\cdot i\}$, $i=1..5$, $t=T=6016$.

We observed that the variability (uncertainty) of the mass fraction might vary from one grid point to another. At some points (dark blue regions), the solution does not change. At other points (white-yellow regions), the variability is very low or high (dark red regions). In regions with high uncertainty, refining the mesh and applying the MLMC method make sense.

\begin{figure}[htbp!]
\begin{center}
  \includegraphics[width=0.245\textwidth]{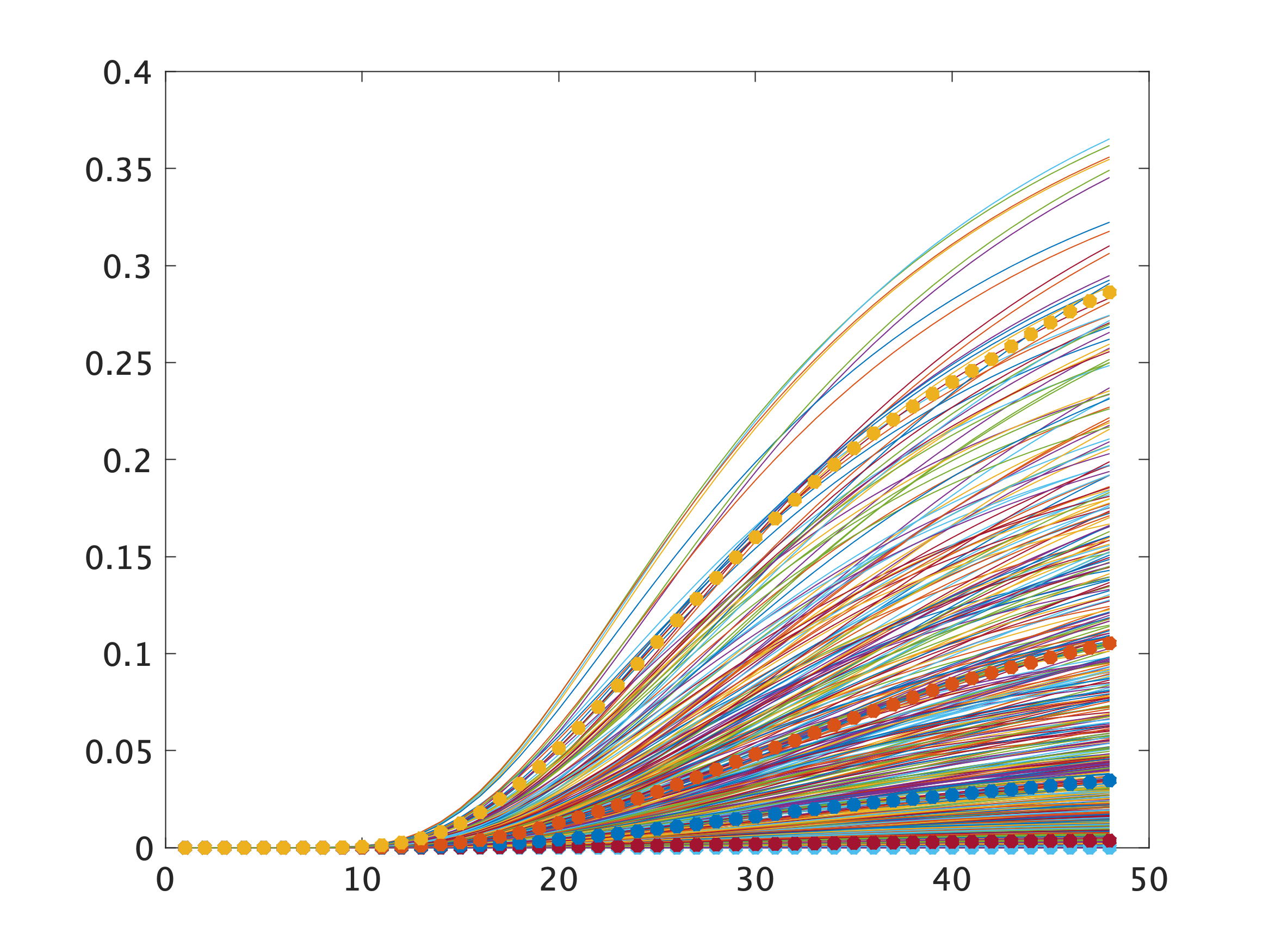}
  \includegraphics[width=0.245\textwidth]{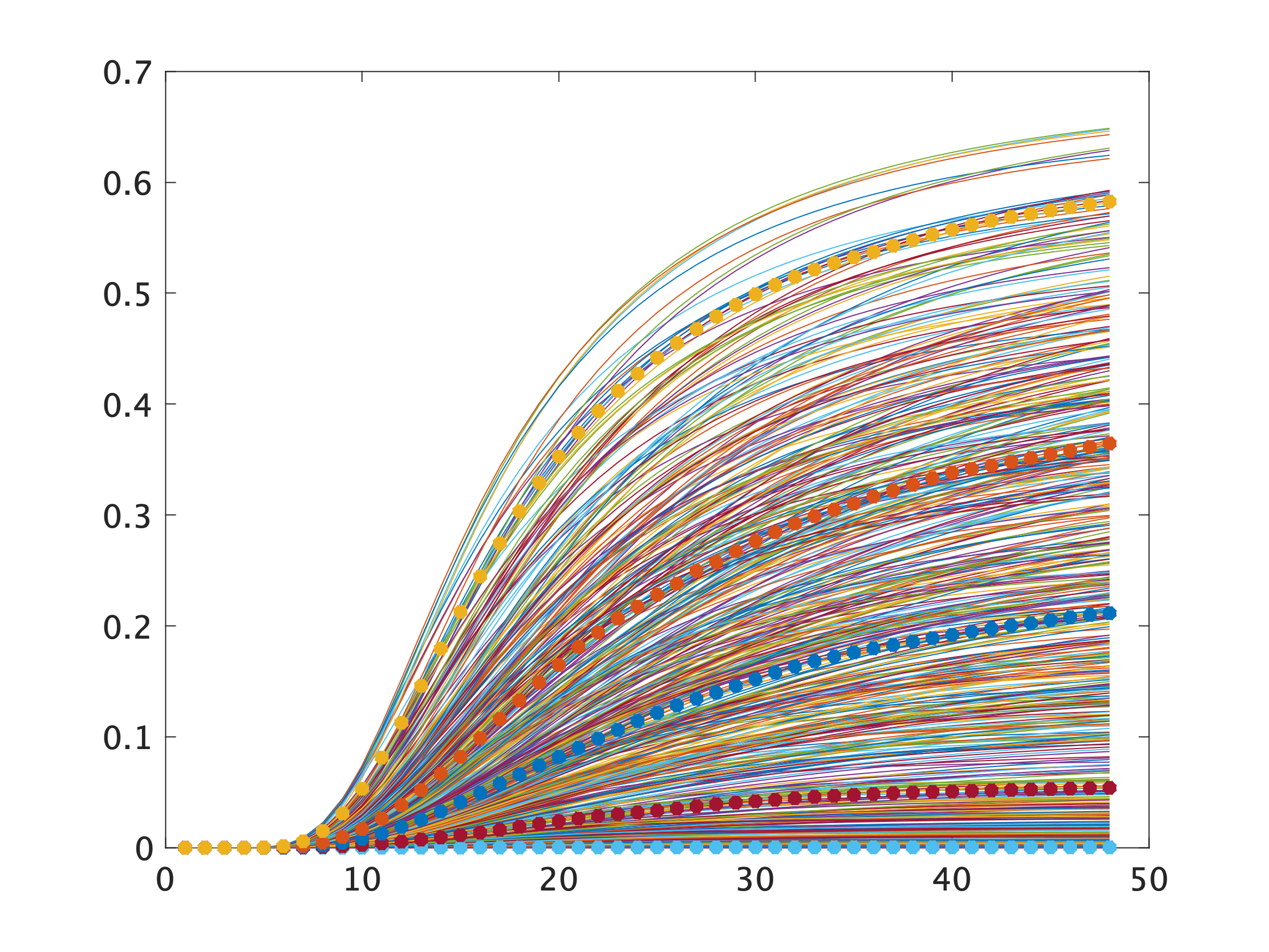}
  \includegraphics[width=0.245\textwidth]{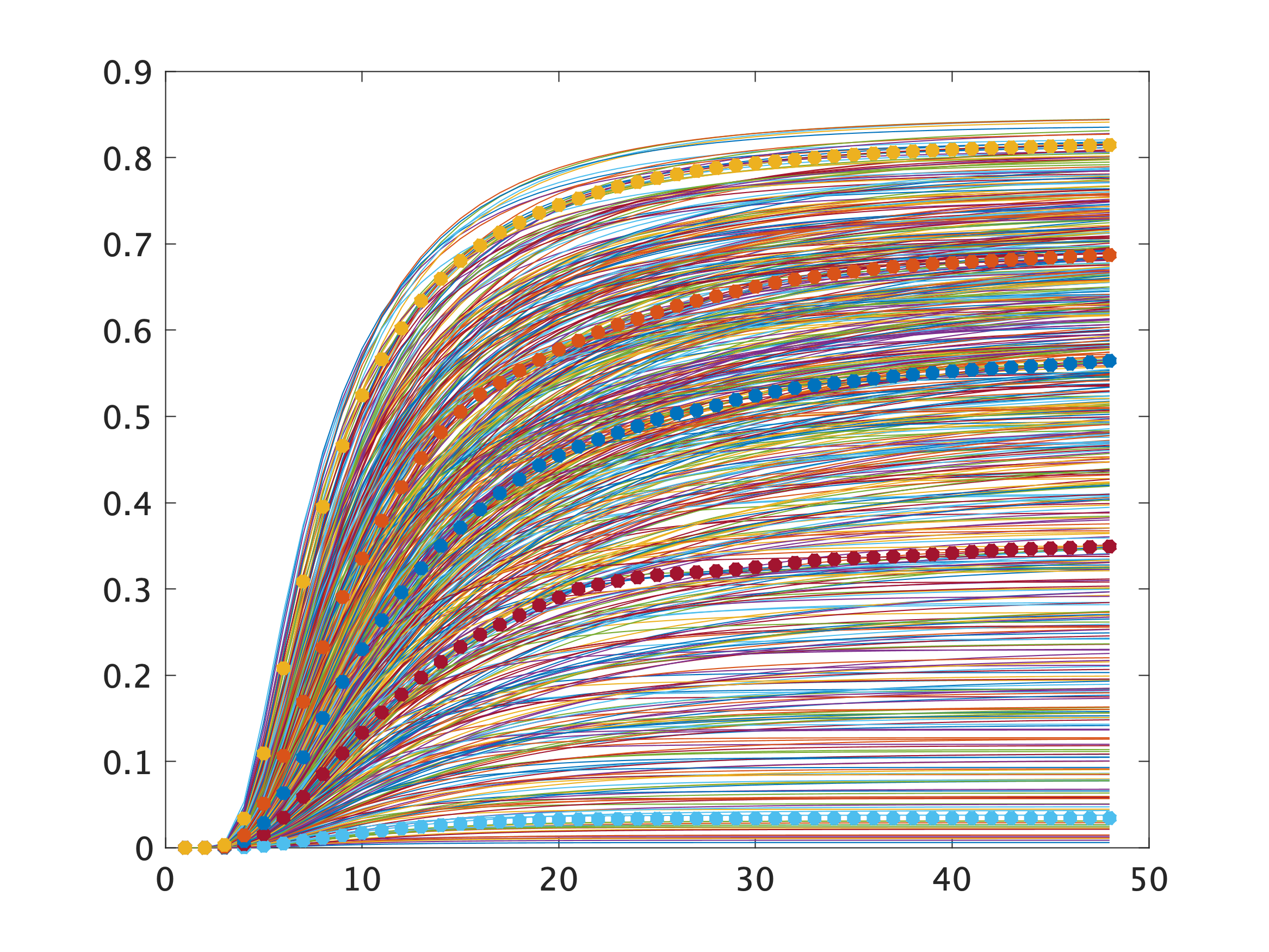}
  \includegraphics[width=0.245\textwidth]{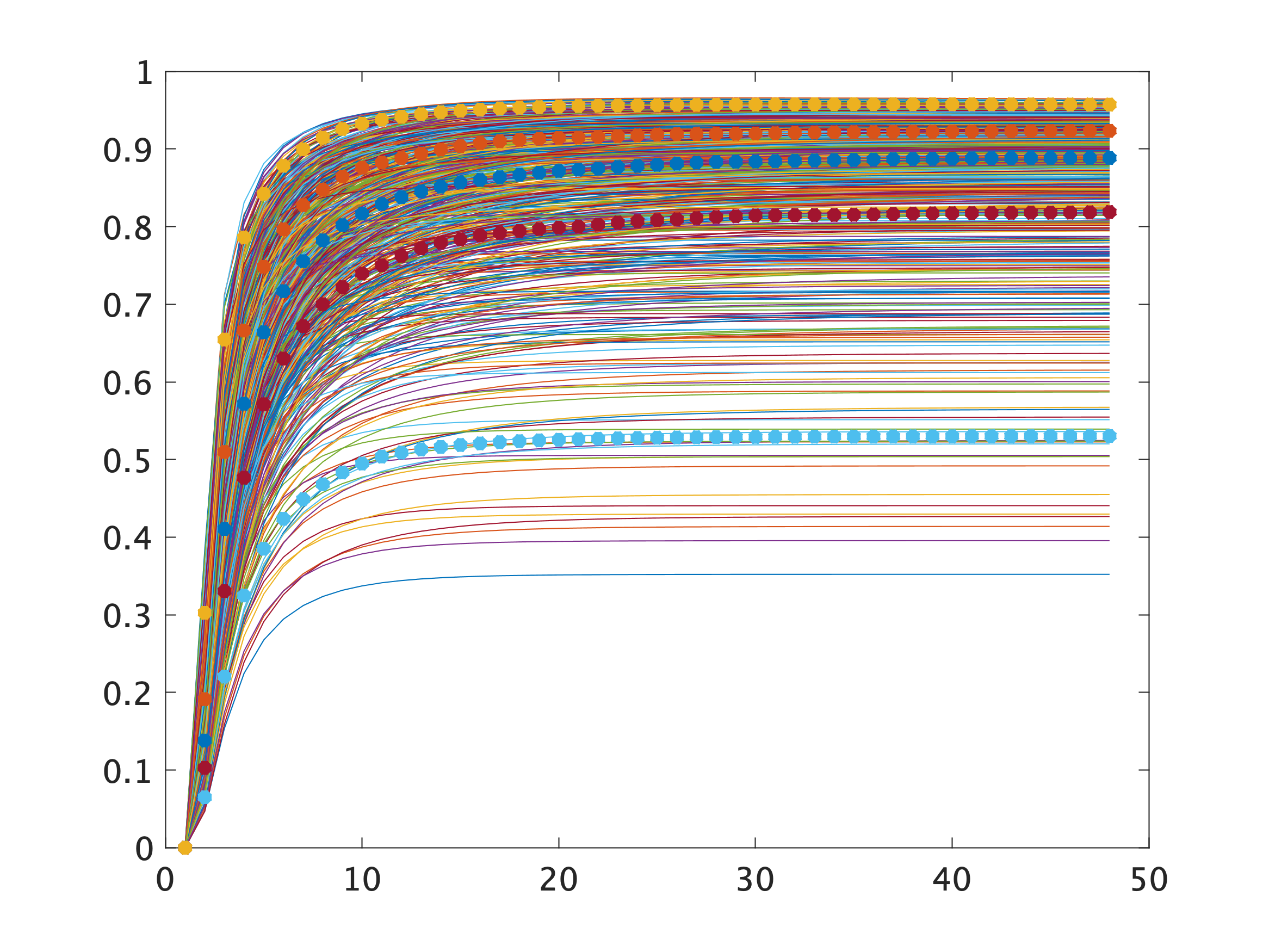}\\
  \includegraphics[width=0.245\textwidth]{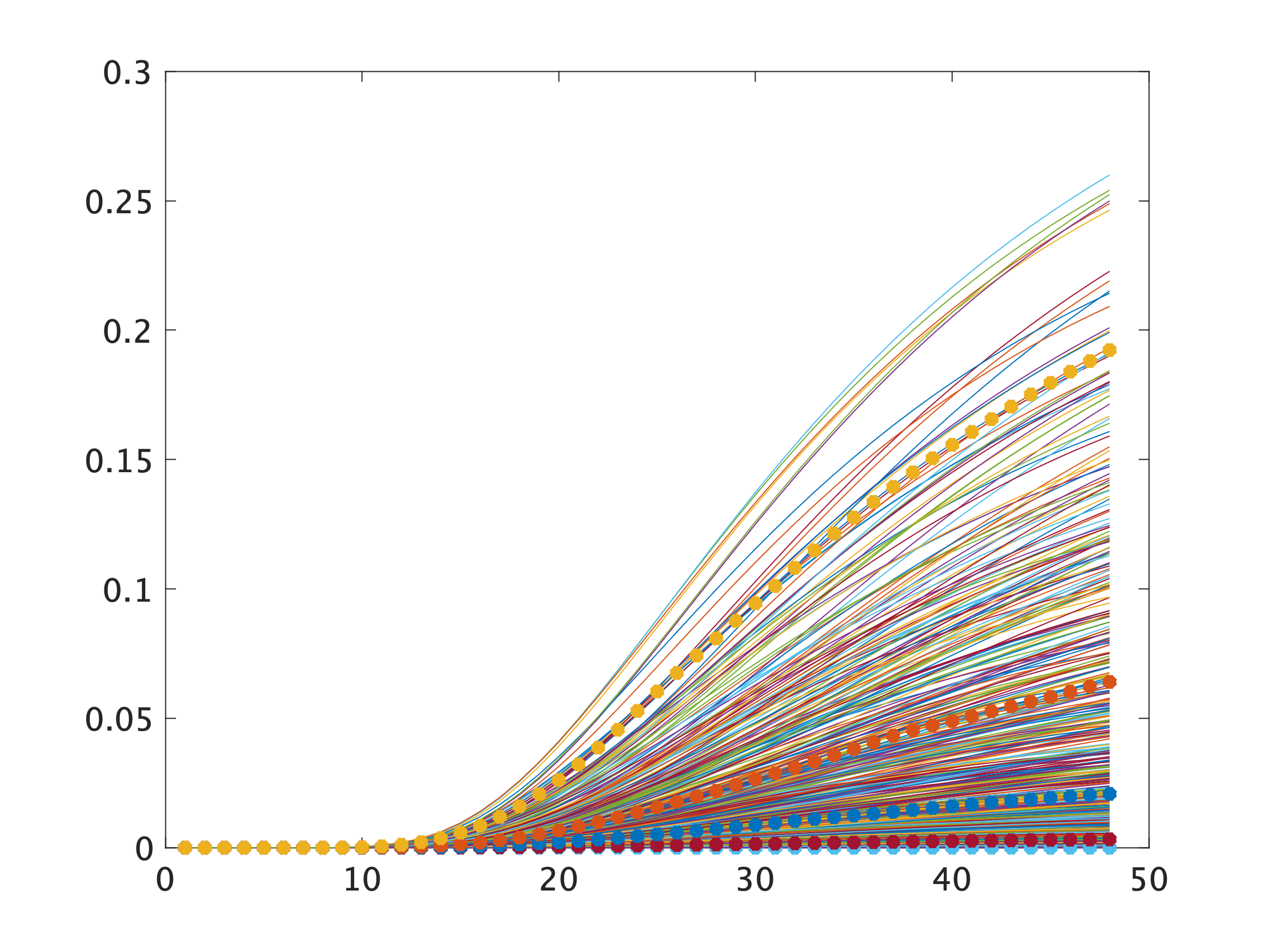}
  \includegraphics[width=0.245\textwidth]{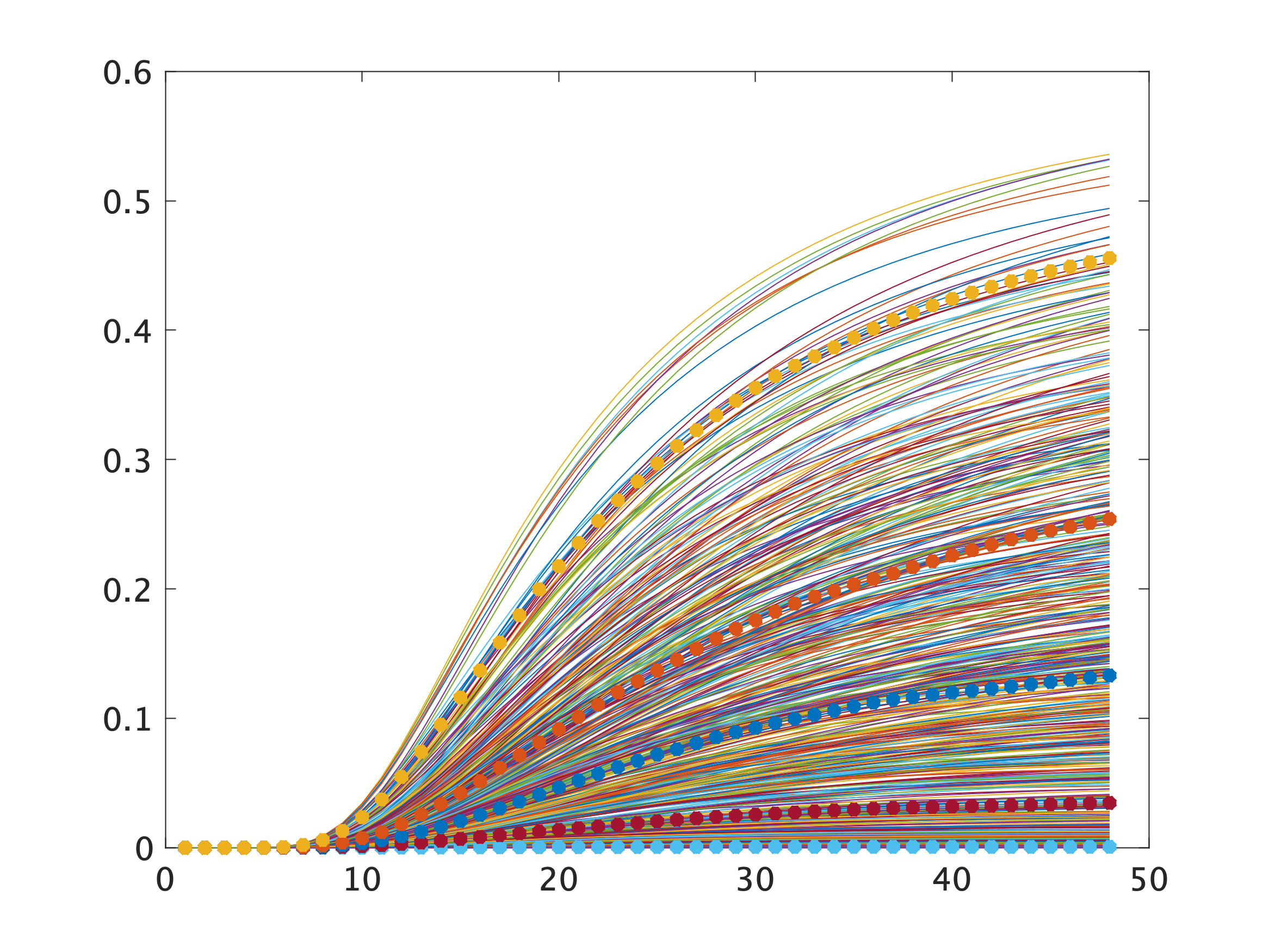}
  \includegraphics[width=0.245\textwidth]{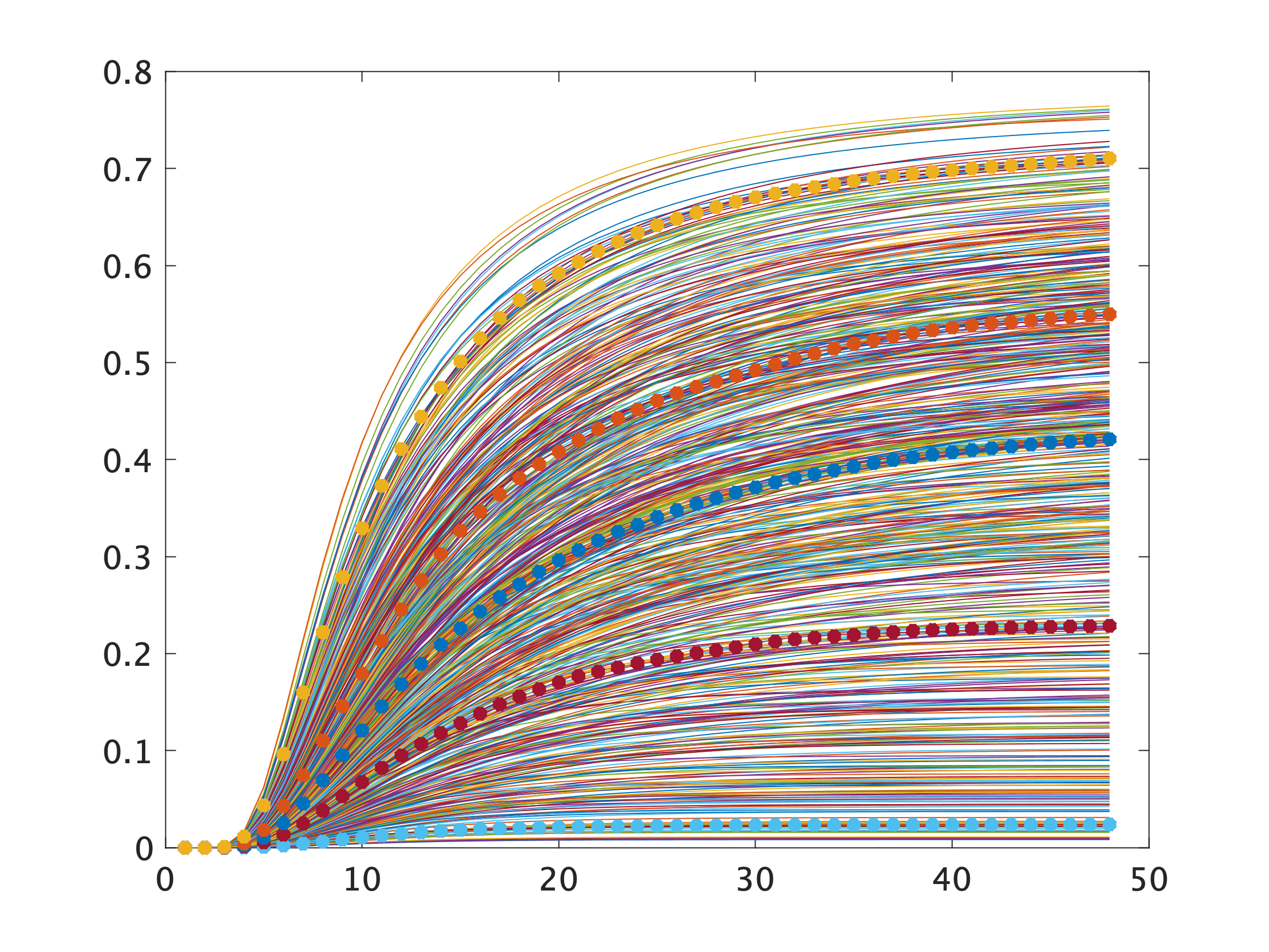}
  \includegraphics[width=0.245\textwidth]{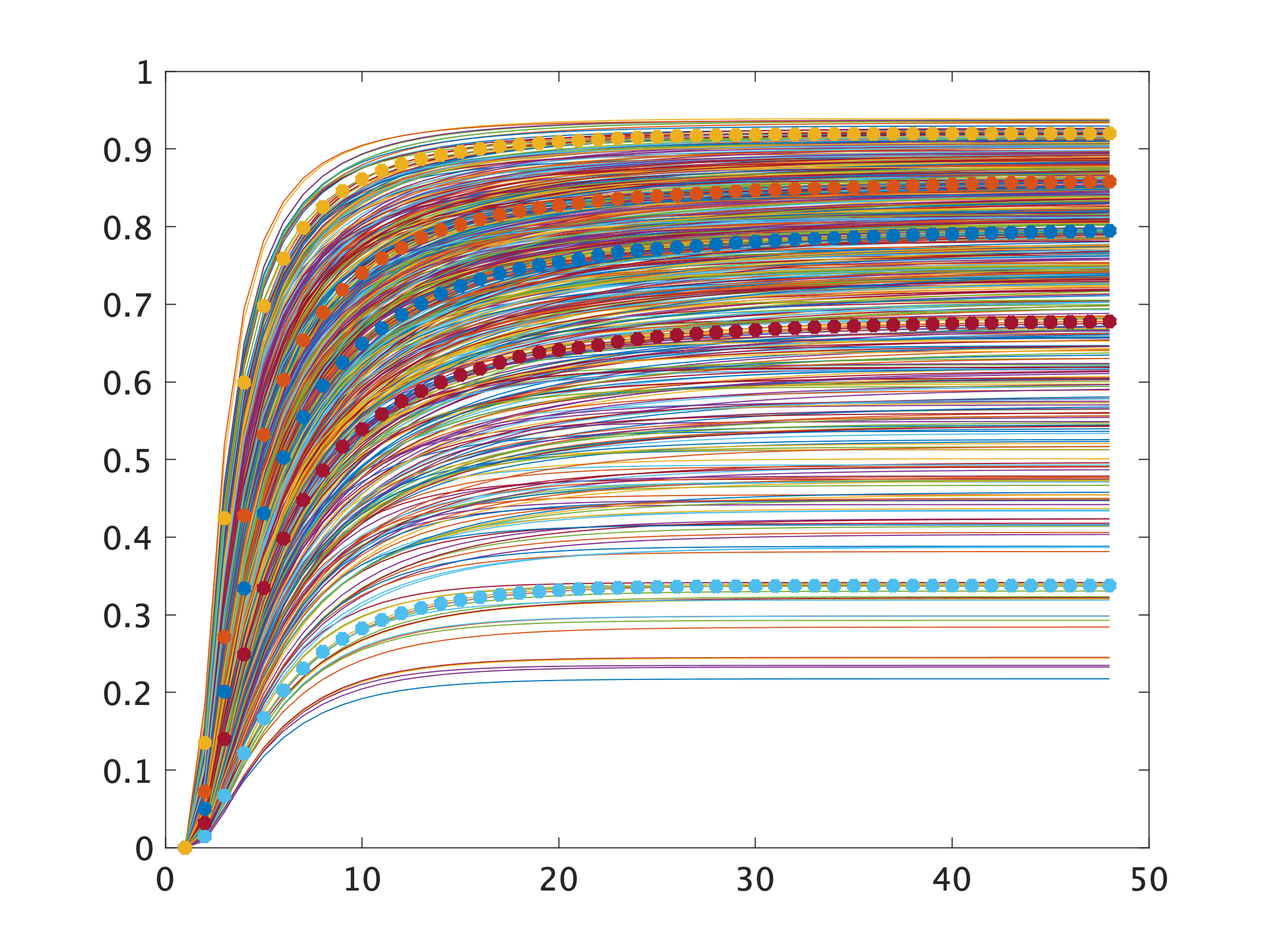}\\
  \includegraphics[width=0.245\textwidth]{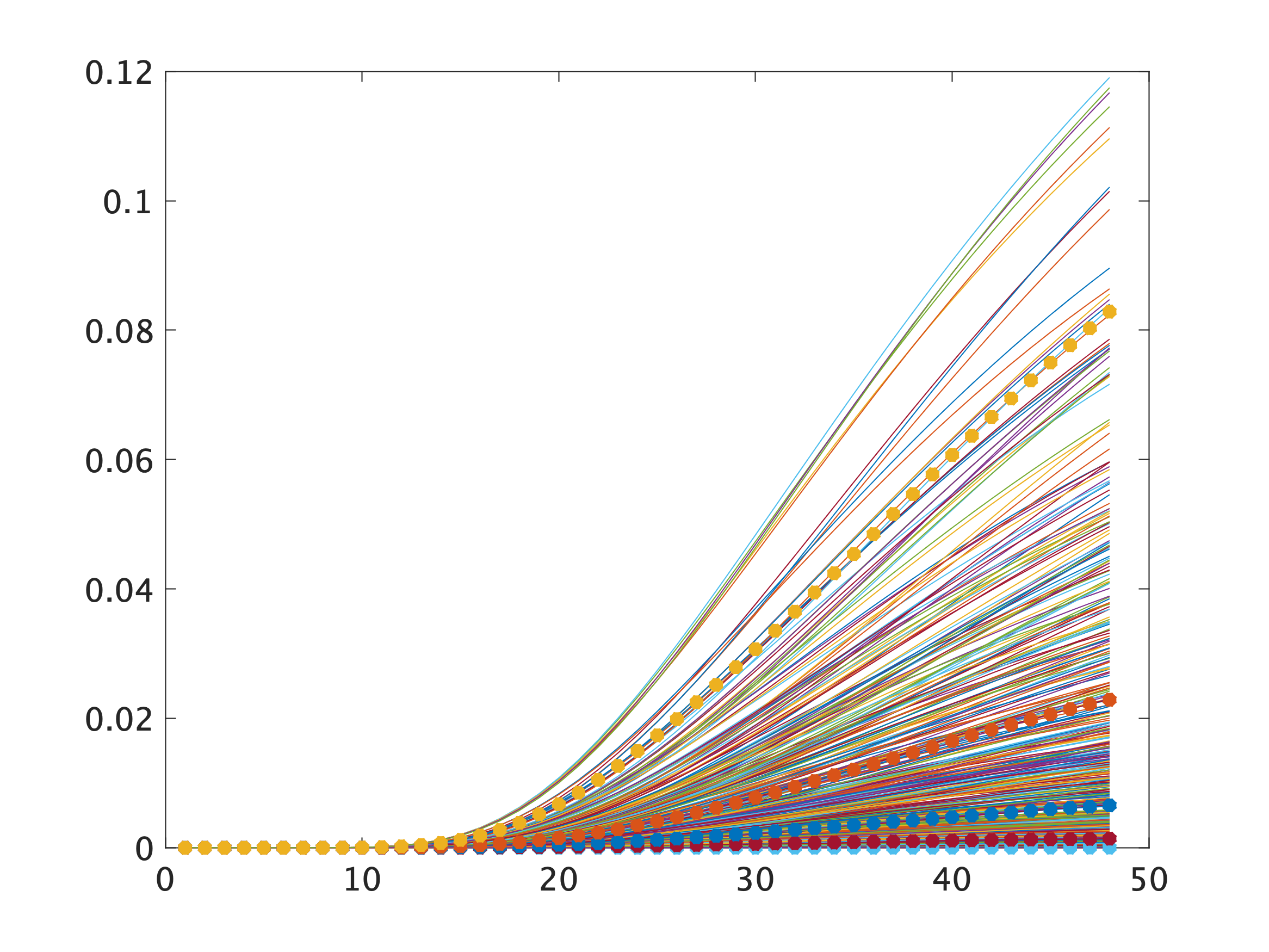}
  \includegraphics[width=0.245\textwidth]{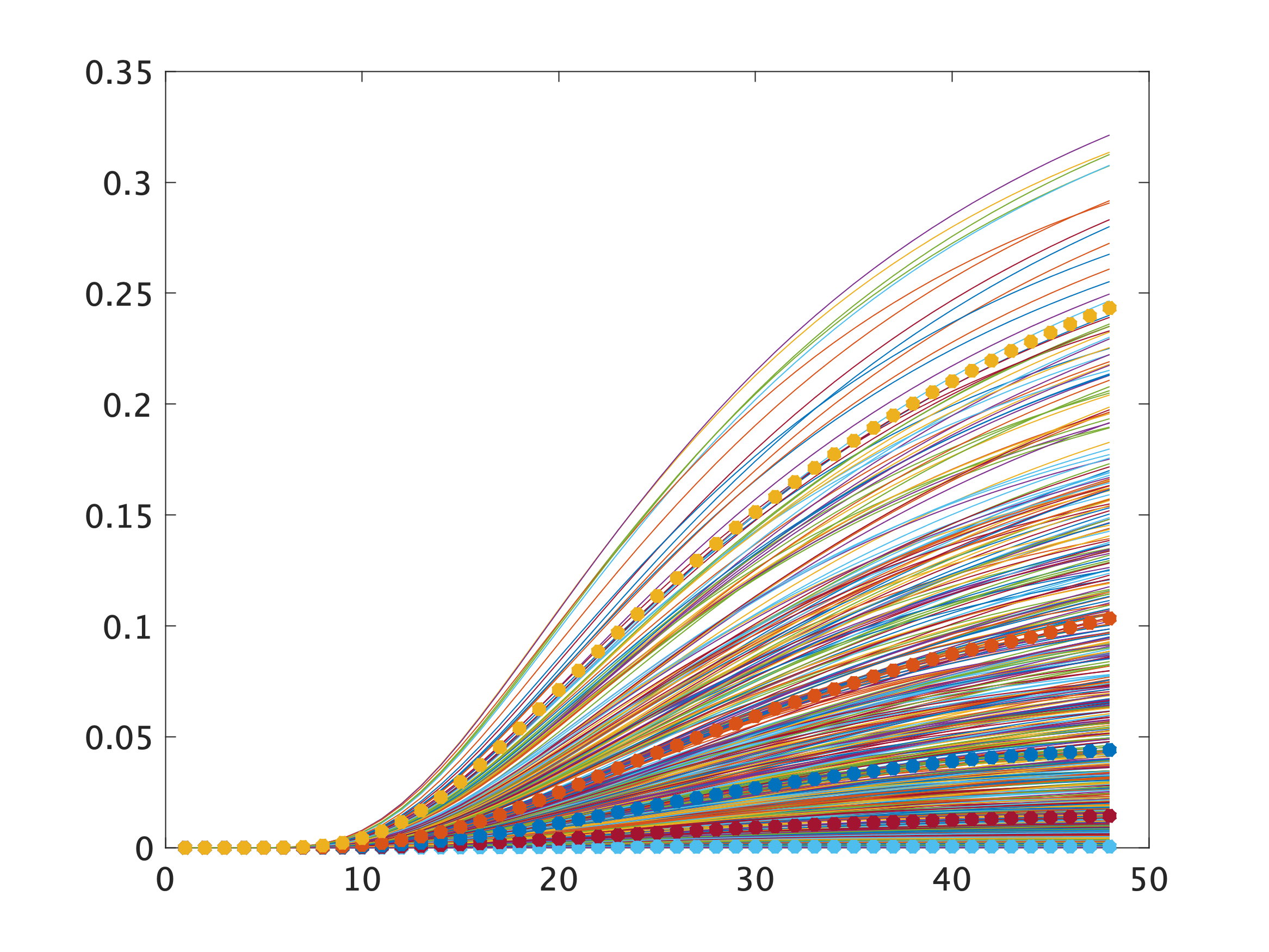}
  \includegraphics[width=0.245\textwidth]{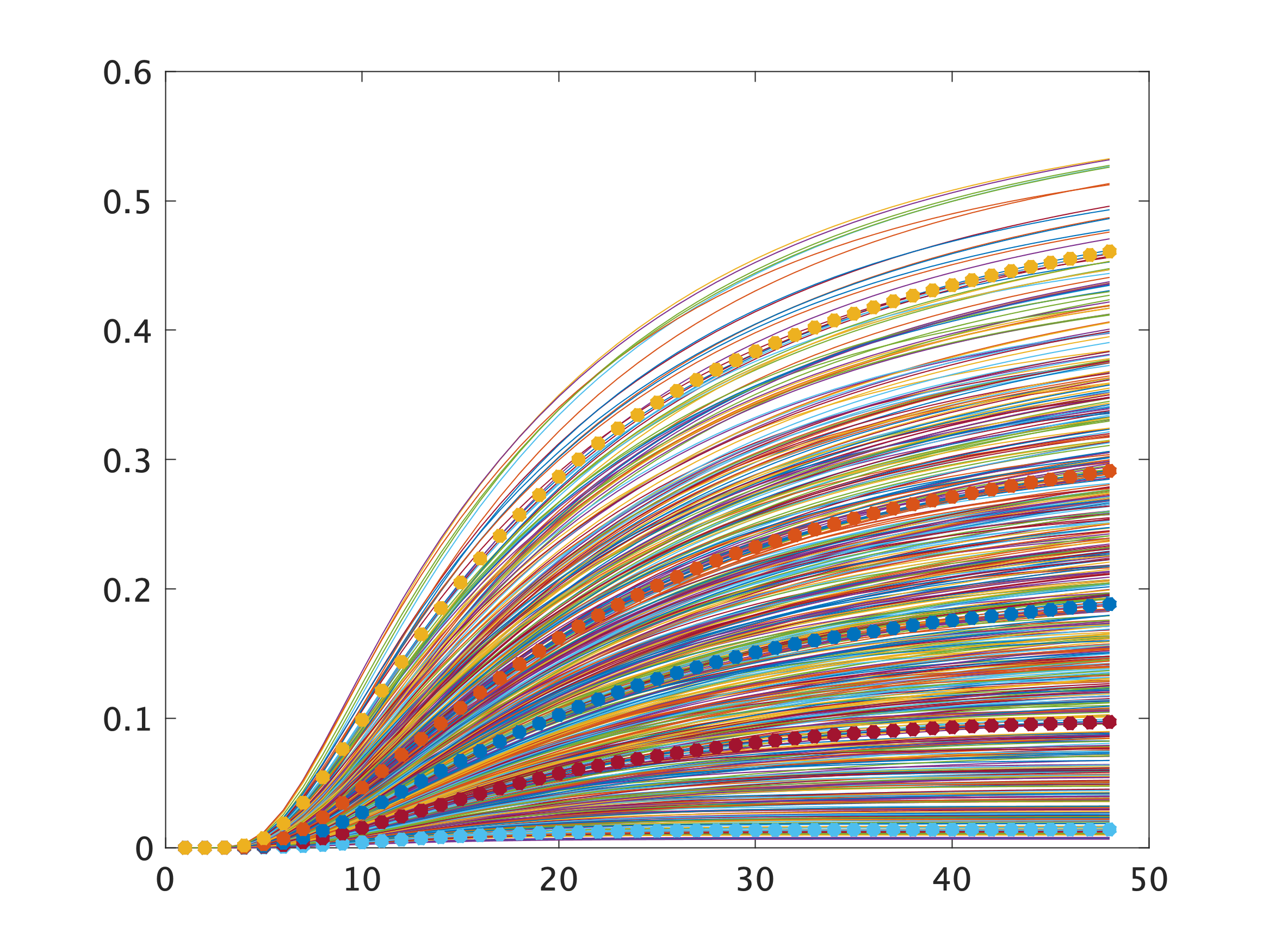}
  \includegraphics[width=0.245\textwidth]{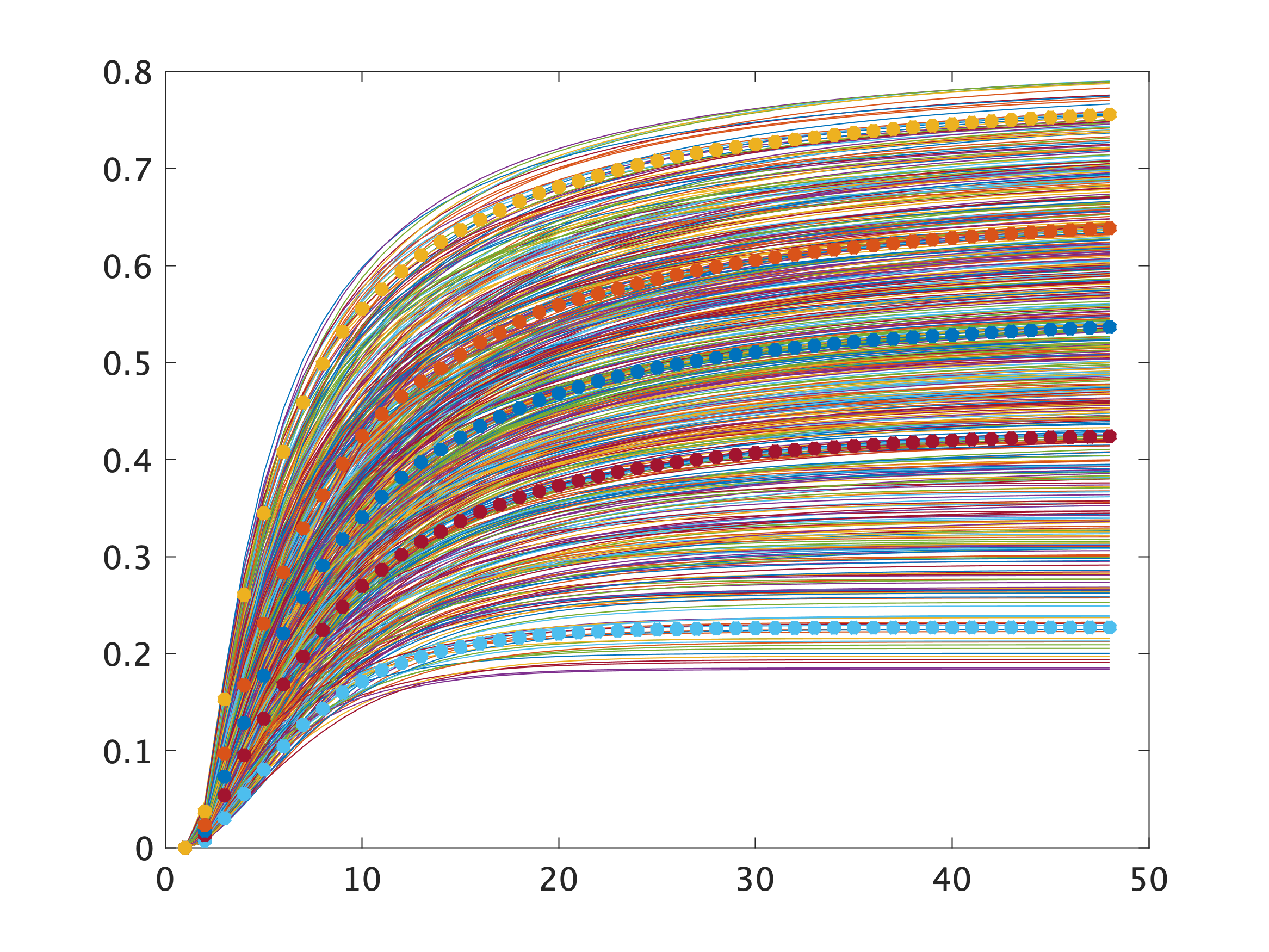}
\caption{Six hundred QMC realisations of $\sol(t,\bx)$ at 12 $\bx$-points listed in Eq.~(\ref{eq:12points}). First row points: $\{(1.10, -0.95)$,$(1.35, -0.95)$,$(1.60, -0.95)$,$(1.85, -0.95)\}$, second row points: $\{(1.10, -0.75)$,$(1.35, -0.75)$,$(1.60, -0.75)$,$(1.85, -0.75)\}$, and third row points: $\{(1.10, -0.50)$,$(1.35, -0.50)$,$(1.60, -0.50)$,$(1.85, -0.50)\}$. Dotted lines from the bottom to the top indicate the quantiles 0.025, 0.25, 0.50, 0.75, and 0.975, respectively.}
 \label{fig:600reslisations_quantiles}
\end{center}
\end{figure}

\begin{figure}[htbp!]
\begin{center}
  \includegraphics[width=0.47\textwidth]{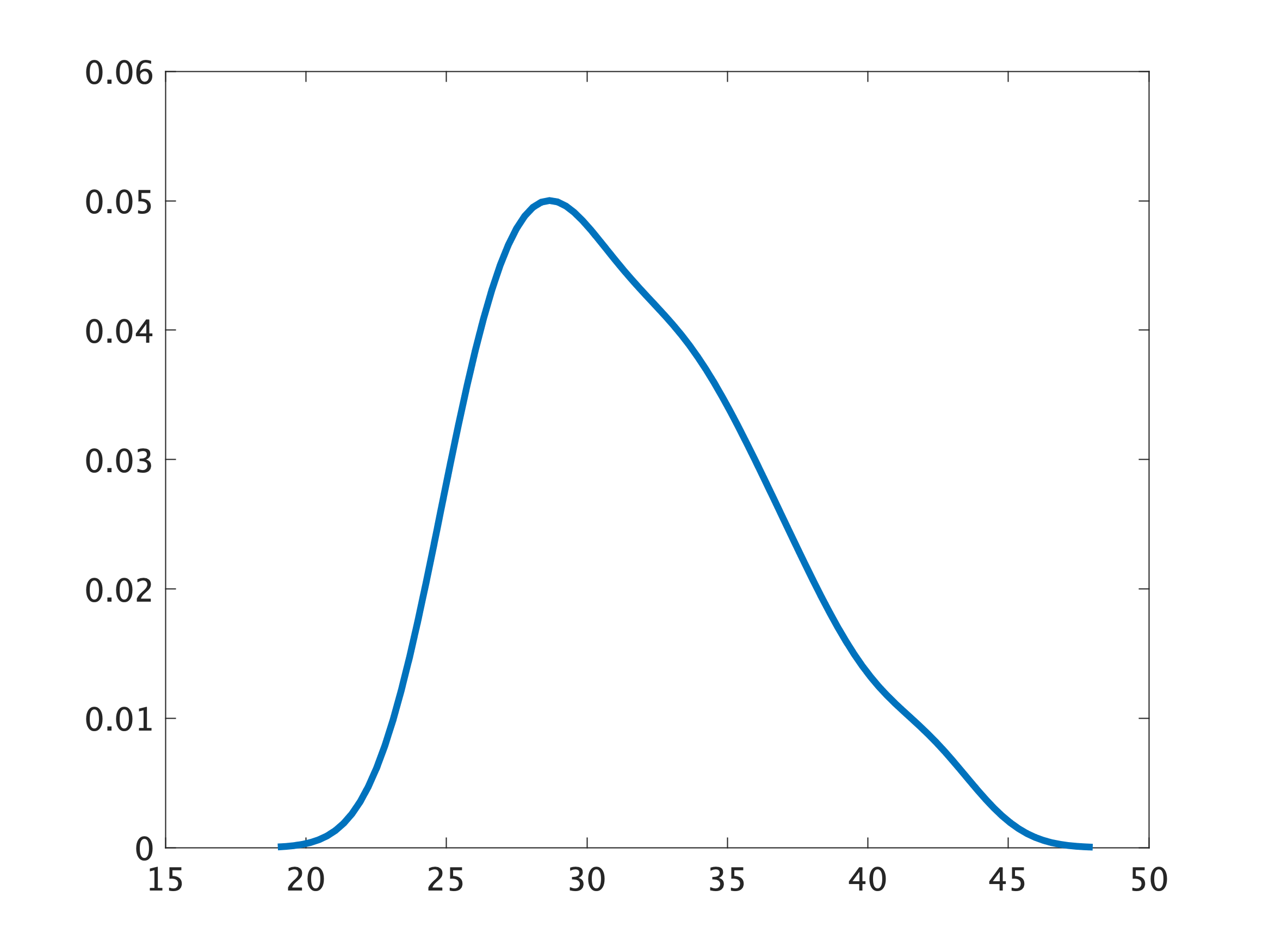}\;
  \includegraphics[width=0.47\textwidth]{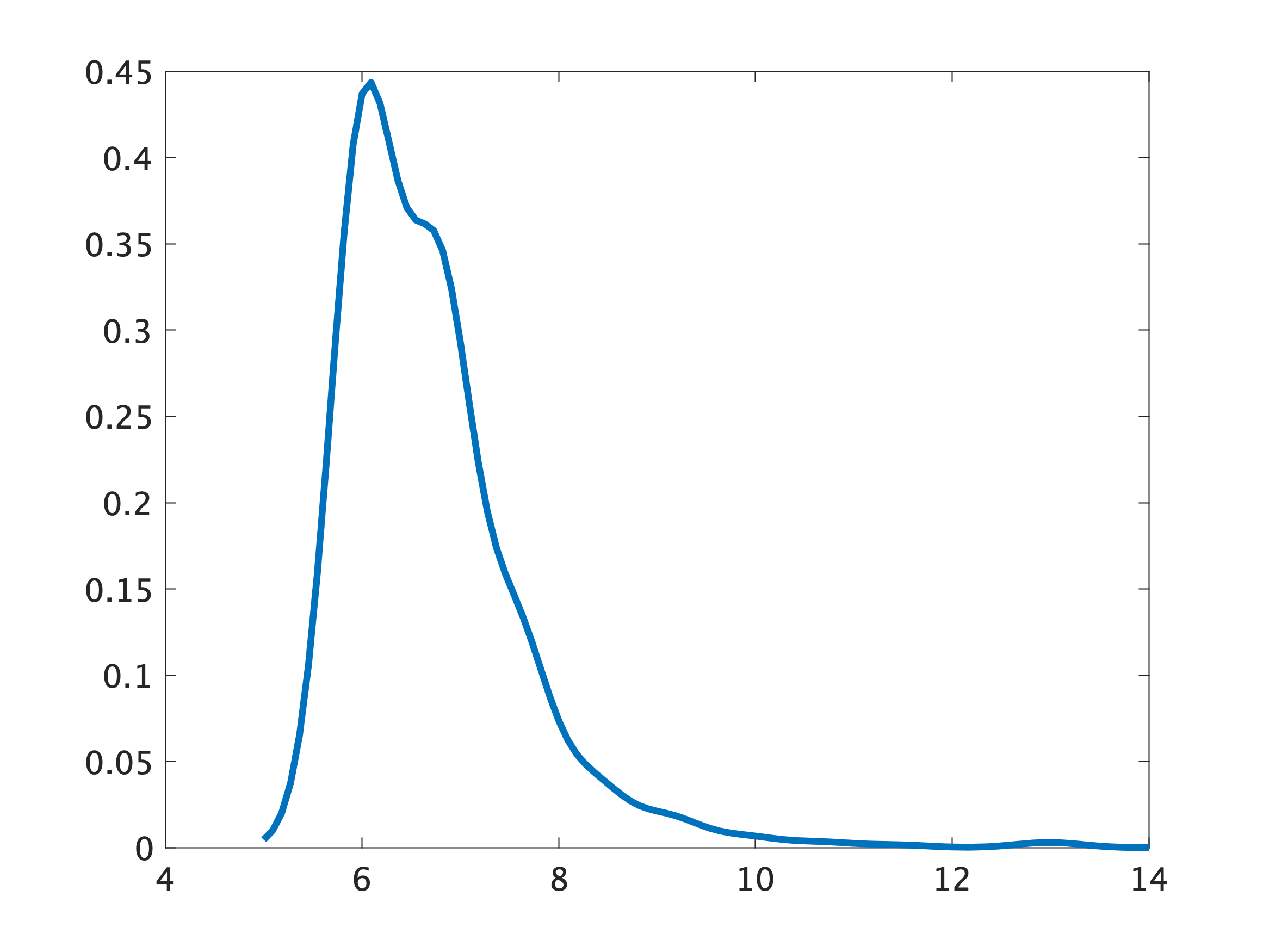}\;
\caption{The pdf of the earliest time point when the freshwater integral $Q_{FW}$ becomes smaller than 1.2 (left) and 1.7 (right). The $x$-axis represents time points.}
\label{fig:600pdf_days}
\end{center}
\end{figure}

Before we run the MLMC method, we first examine the solution $\sol(t,\bx)$ at 12 preselected points (see Eq.~(\ref{eq:12points})). Figure~\ref{fig:600reslisations_quantiles} includes 12 subfigures. Each subfigure presents 600 QMC realizations of $\sol(t,\bx)$ and five quantiles depicted by dotted lines. The dotted line at the bottom indicates the quantile $0.025$. The following dotted line is the quantile $0.25$, and the dotted line on the top indicates the quantile $0.975$. All five quantiles from the bottom to the top are 0.025, 0.25, 0.50, 0.75, and 0.975, respectively. We observe that $\sol$ at the final point $t=T$ varies considerably.
\newpage
\textbf{Example.} In Fig.~\ref{fig:600pdf_days}, we demonstrate the probability density function (pdf) of $t^*(\omega) = \min_t\{t:\; Q_{FW}(t,\omega)<1.2\}$ (left), and the pdf of $t^*(\omega) = \min_t\{t:\; Q_{FW}(t,\omega)<1.7\}$ (right). On average, after approximately 29 time steps (on the left) and six time steps (on the right), the volume of the fresh water becomes less than 1.2 and 1.7, respectively. The initial volume of the fresh water was 2.0.

\begin{figure}[htbp!]
\begin{center}
   \includegraphics[width=0.68\textwidth]{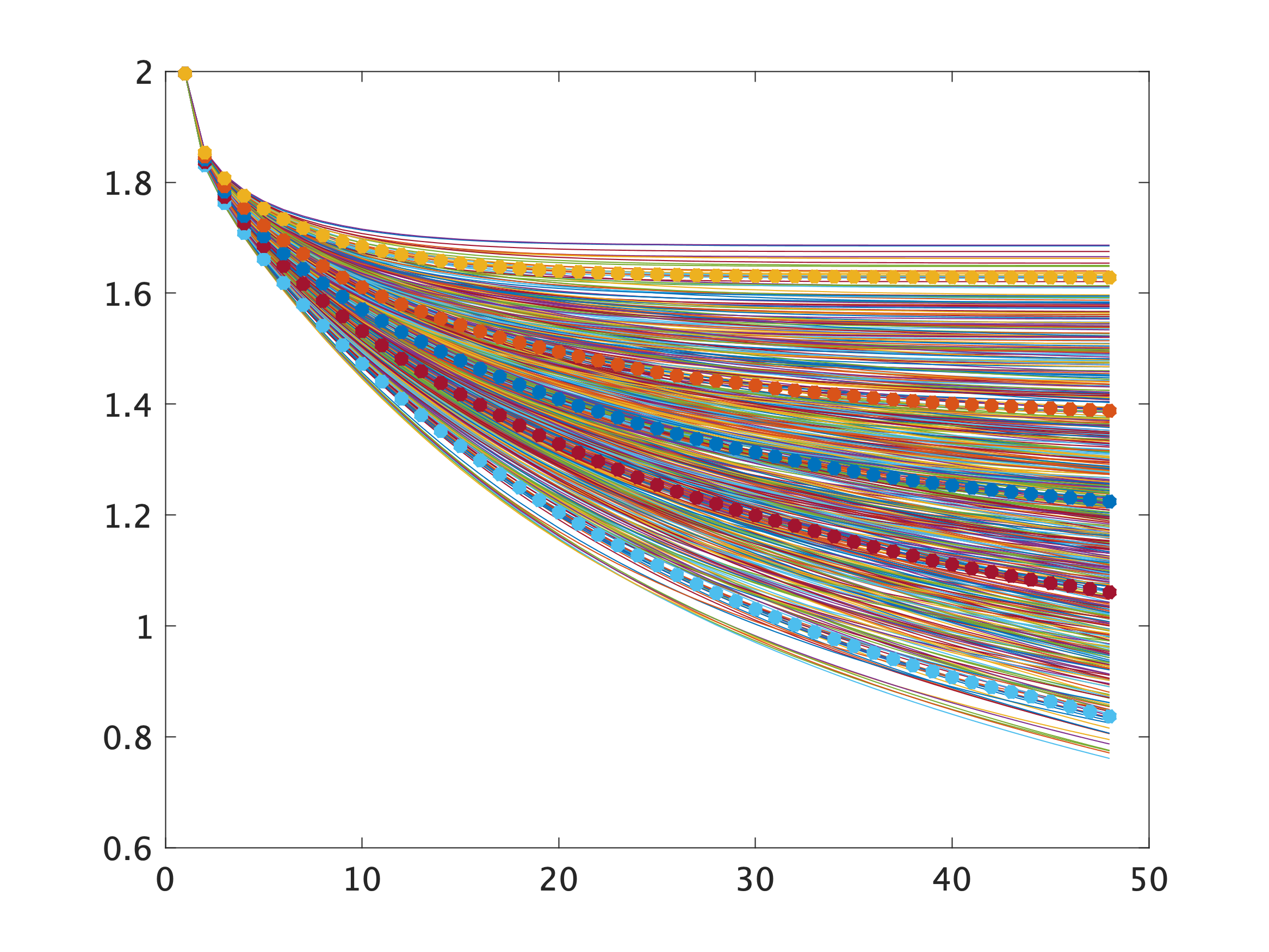}\;  
\caption{Six hundred realizations of $Q_{FW}(t)$. The $x$-axis represents time $t=1\tau,\ldots,48\tau$; dotted curves denote five quantiles: 0.025, 0.25, 0.50, 0.75, and 0.975 from the bottom to the top.}
 \label{fig:600pure_water_realisations}
\end{center}
%Matlabcode/exceed_prob2_pure.m
\end{figure}

All 600 realiations of $Q_{FW}(t)$ are depicted in Fig.~\ref{fig:600pure_water_realisations} 
%(and all 600 scenarios of $\sol$ are shown in Fig.~\ref{fig:600reslisations_quantiles}, last sub-figure in the first row). 
The time is along the $x$-axis, $t\in[\tau,48\tau]$. Additionally, five quantiles are represented by dotted curves from the bottom to the top and are 0.025, 0.25, 0.50, 0.75, and 0.975, respectively. 

\begin{figure}[htbp!]
\begin{center}
   \includegraphics[width=0.48\textwidth]{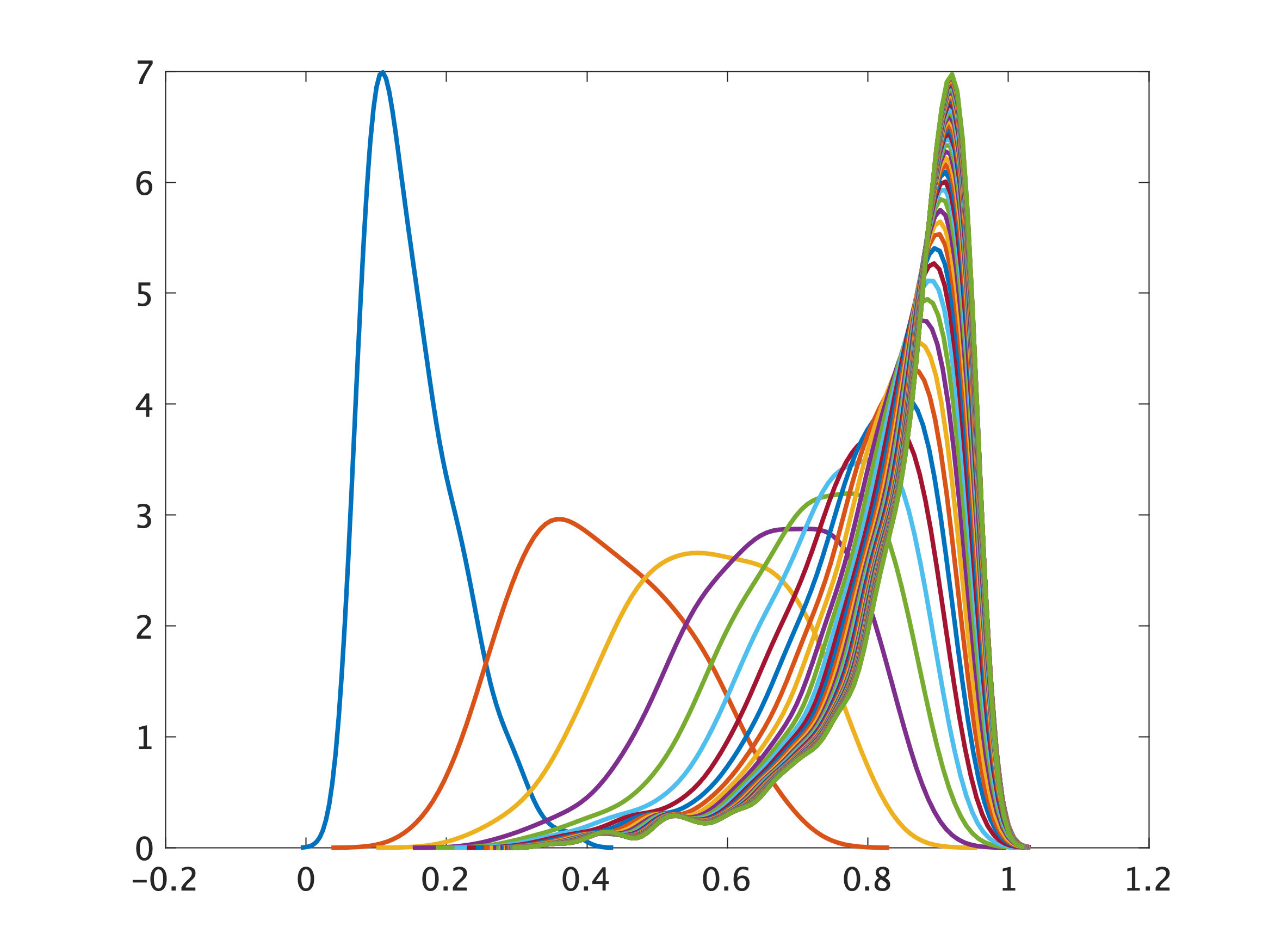}\; 
   \includegraphics[width=0.48\textwidth]{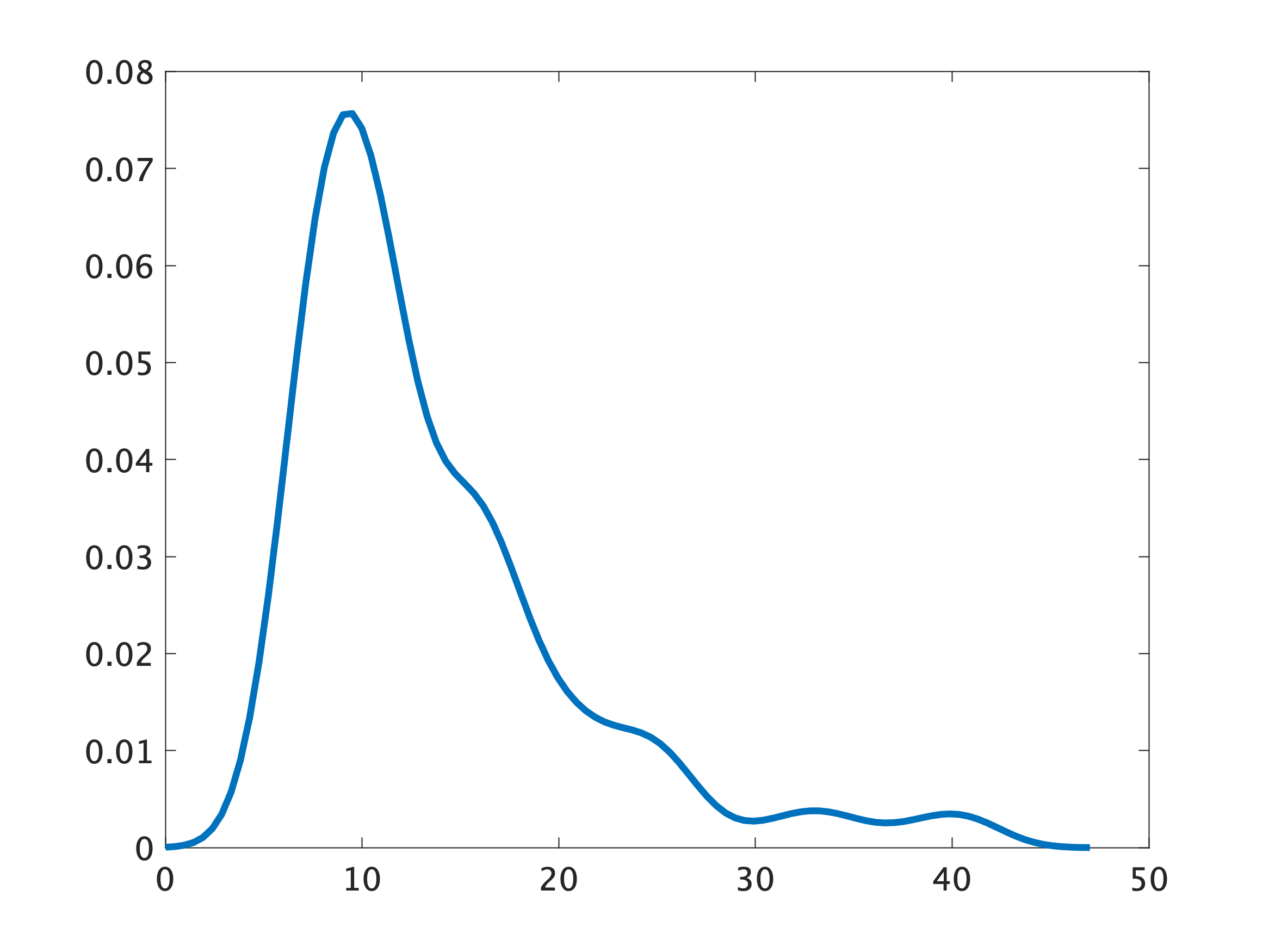}
\caption{(left) Evolution of the pdf of $\sol(t,\bx)$ for $t=\{3\tau,\ldots,48\tau\}$. (right)
The pdf of the earliest time point when $\sol(t,\bx)< 0.9$ ($\bx=(1.85, -0.95)$ is fixed).}
 \label{fig:pdfs_p4_L7}
\end{center}
%Matlabcode/exceed_prob1.m
\end{figure}

\textbf{Example.} Figure~\ref{fig:pdfs_p4_L7} (left) displays the evolution of the pdf of $\sol(t,\bx,\omega)$ at a fixed point $\bx=(1.85, -0.95)$ in time $t=\{3\tau,\ldots,48\tau\}$.  From left to right, the farthest left (blue) pdf corresponds to $t=3\tau$, the second curve from the left (red) corresponds to $t=4\tau$, and so on. In the beginning, $t=3\tau$, and the mass fraction $\sol$ is low, about 0.15 on average. Then, with time, $\sol$ increases and, at $t=T=48\tau$, is approximately equal to 1. 
\textbf{Example.} The next QoI is the earliest time moment when $\sol(t,\bx)$, at fixed $\bx=(1.85, -0.95)$, becomes smaller than the threshold value $0.9$ (maximum is 1.0). Figure~\ref{fig:pdfs_p4_L7} (right) presents its pdf. On average, after $t\approx 10$ time steps, the mass fraction becomes smaller than $0.9$, but 40 time steps are needed in some scenarios.

%
% \begin{figure}[htbp!]
% \begin{center}
%   \includegraphics[width=0.4\textwidth]{111Ekaterina/Poro-lev8-100of100H-0261175-0459453.png}\;
% %  \includegraphics[width=0.3\textwidth]{111Ekaterina/Poro-lev8-1of100H-0315-0395.png}\;
%   \includegraphics[width=0.4\textwidth]{111Ekaterina/Poro-lev8-50of100H-0231704-0473158.png}
%   \caption{Two realisations of the porosity field. (left) $\poro \in [0.2612, 0.4594]$; (right) $\poro\in[0.2317,0.4732]$.}
%     \label{fig:Henry2d-poro}
% \end{center}    
% \end{figure}

% \begin {figure} [ht!]
% \begin {center}
% \includegraphics[width=0.49\textwidth]{fixedTS600_fig/fihemoca600-lev7-t47-p1-c-adiff_0_1_0_2.png}%
% \;%
% \includegraphics[width=0.49\textwidth]{fixedTS600_fig/fihemoca600-lev7-t47-p1-poro_0_315_0_385.png}%
% \caption {Left: Mass fraction $c \in [0,1]$ (color) and isolines for values $0.1$ and $0.2$ of $|c - \overline{c}| \in [0, 0.3)$ for the scenario with $(\theta_0, \theta_1, \theta_2) = (0, -0.3333, -0.6)$ at time $t = 6016$ $[s]$; right: Porosity $\phi \in [0.315, 0.385]$ for this scenario.}
% \end {center}
% \end {figure}

%

\newpage

\begin{figure}[htbp!]
\begin{center}
  \includegraphics[width=0.49\textwidth]{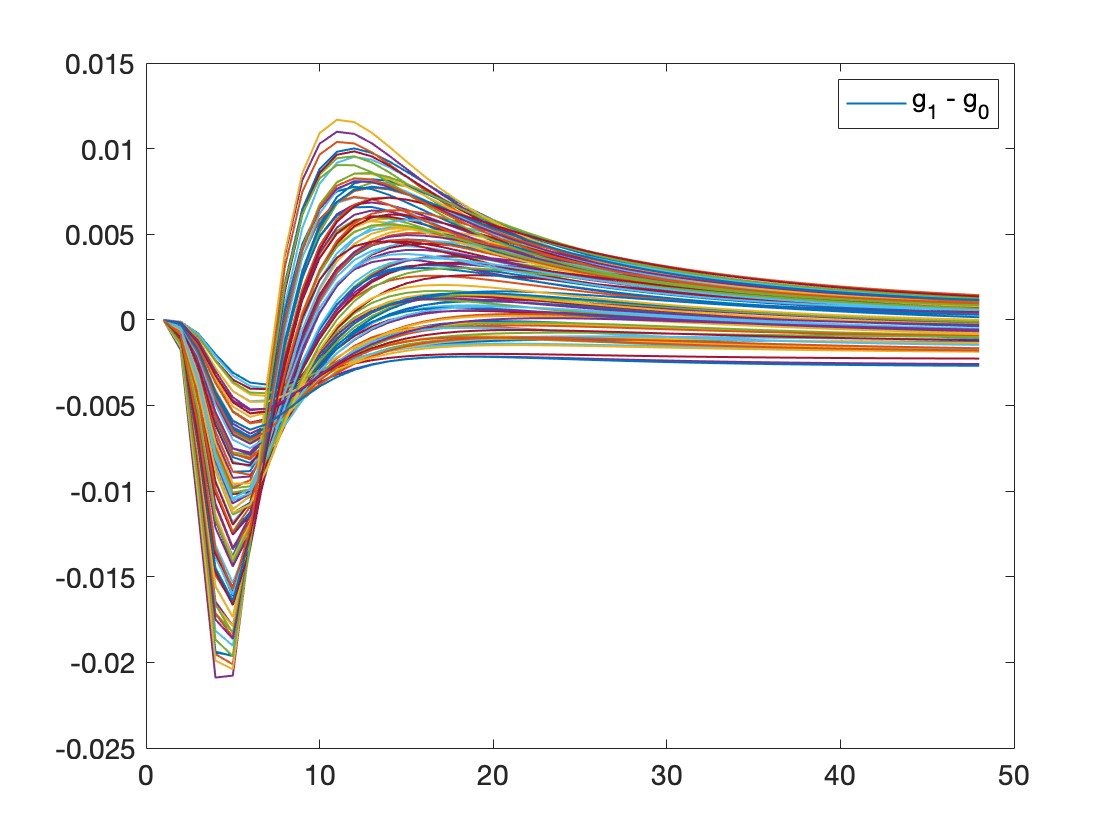}\;
  \includegraphics[width=0.49\textwidth]{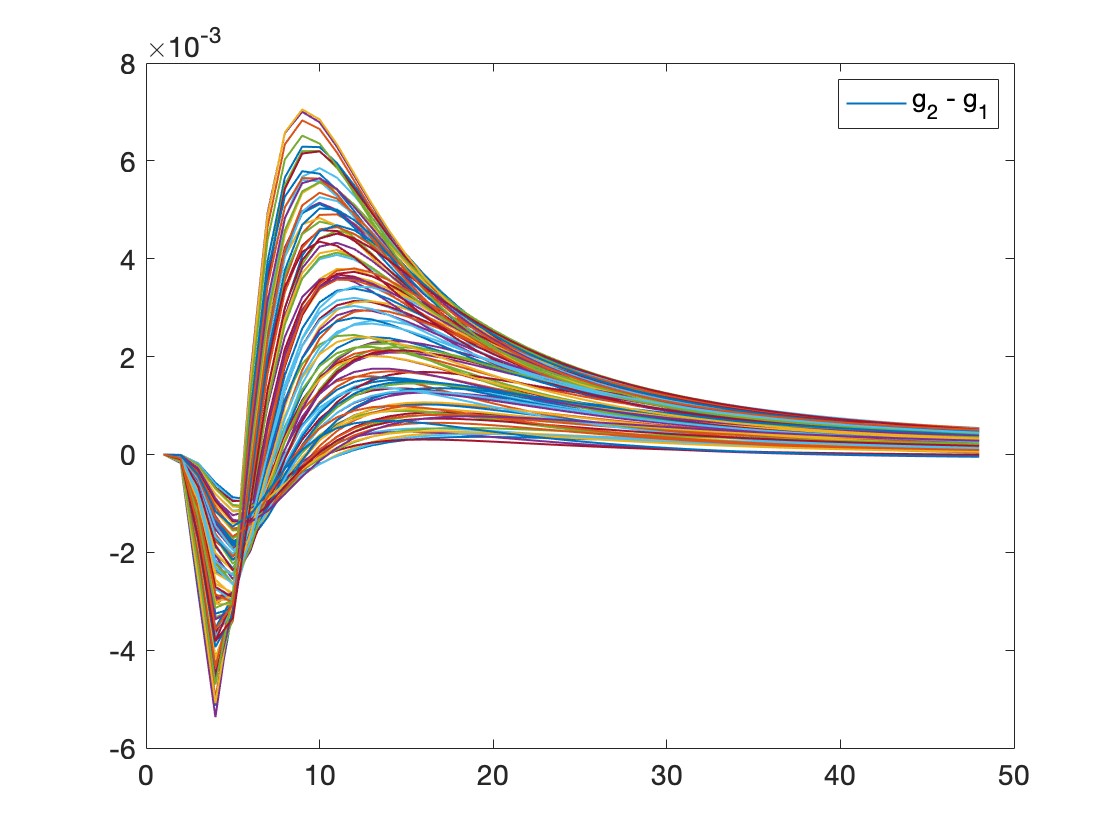}\\
  \includegraphics[width=0.49\textwidth]{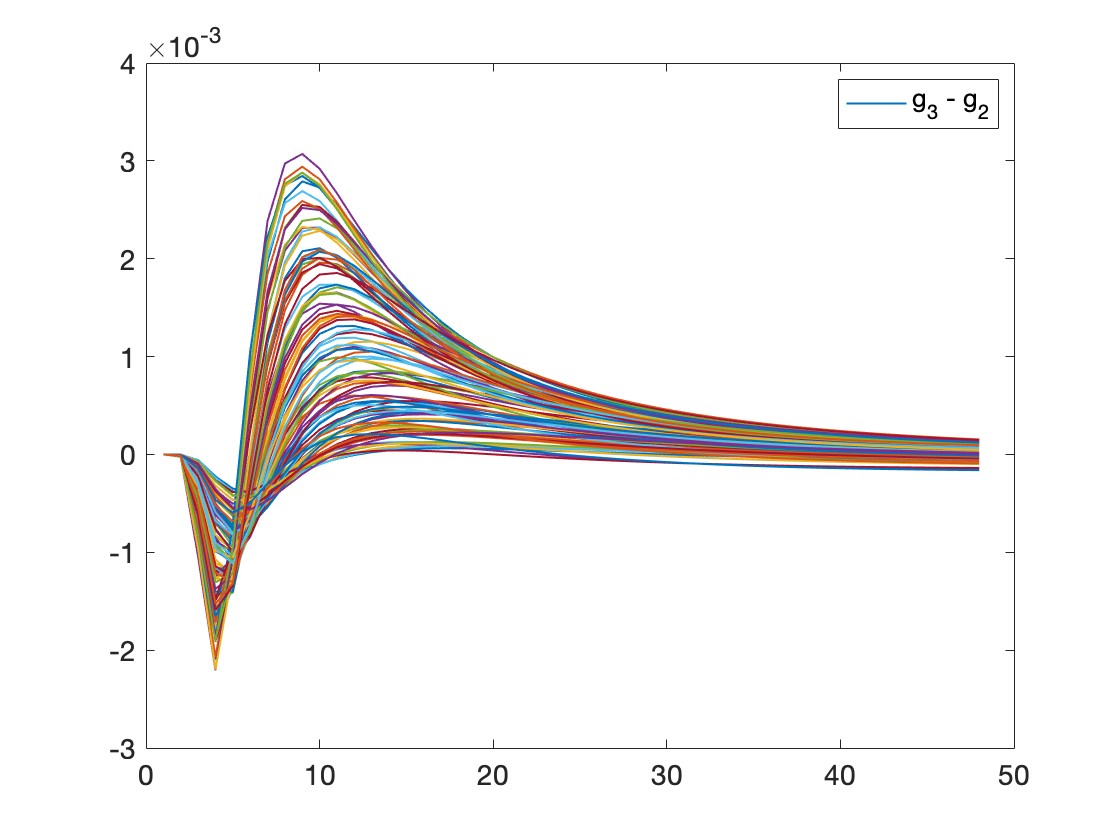}\;
  \includegraphics[width=0.49\textwidth]{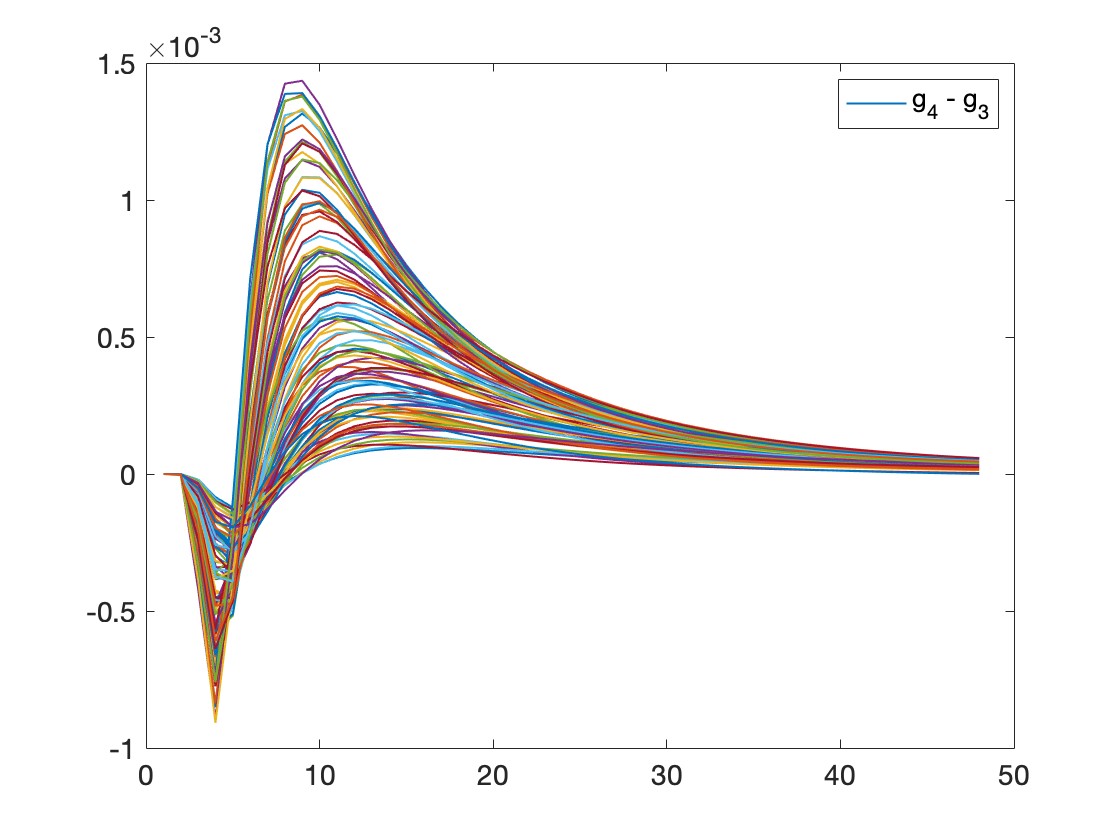}\\
  \includegraphics[width=0.49\textwidth]{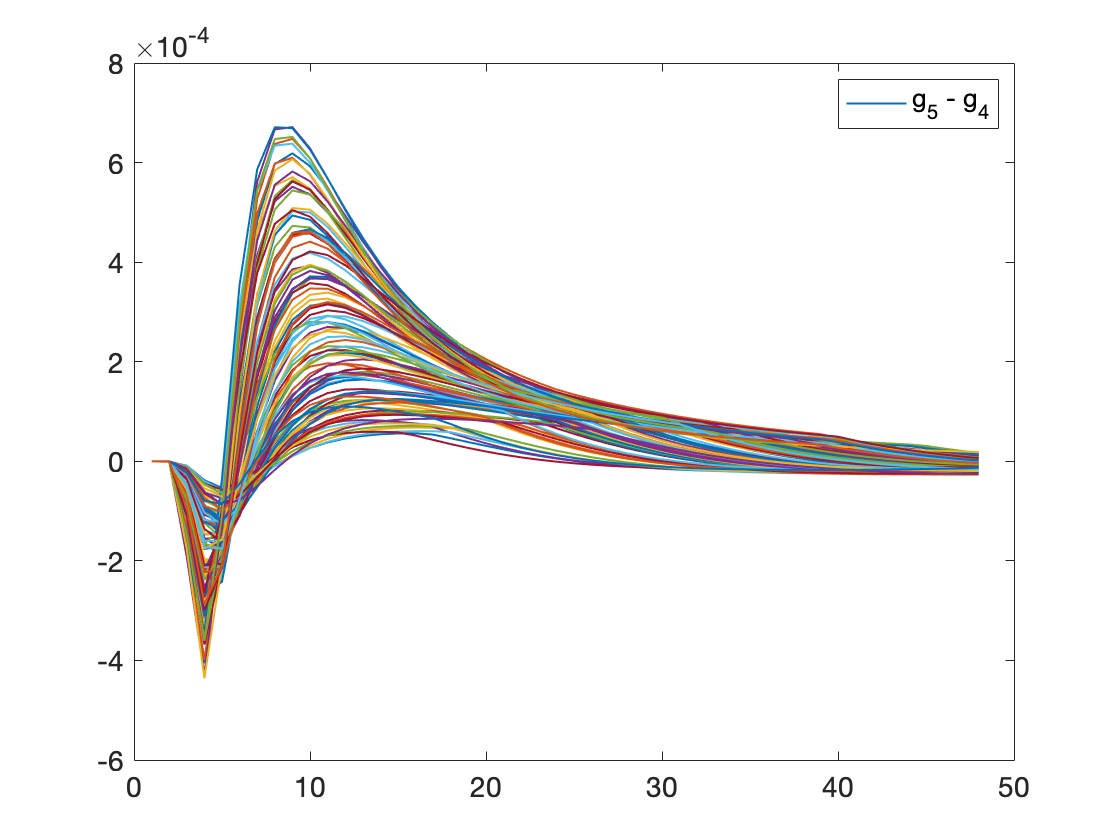}
  \caption{Differences between mass fractions $\sol$ computed at the point $(1.60, −0.95)$ on levels a) 1 and 0, b) 2 and 1 (first row), c) 3 and 2, d) 4 and 3 (second row), and e) 5 and 4 (third row) for 100 realizations ($x$-axis represents time).}
    \label{fig:diffs_levels}
\end{center}
%dima_compare_MC_with_MLMC.m
\end{figure}

Next, we research how $g_{\ell}-g_{\ell-1}$ depends on the time and level. All graphics in Fig.~\ref{fig:diffs_levels} display 100 realizations of the differences between solutions computed on two neighbor meshes for every time point $t_i$, $i=1\ldots 48$ (along the $x$-axis). The top left graphic indicates the differences between the mass fractions computed on Levels 1 and 0. The other graphics reveal the same, but for Levels 2 and 1, 3 and 2, 4 and 3, and 5 and 4, respectively. The largest value decreases from $2.5\cdot 10^{-2}$ (top left) to $5\cdot 10^{-4}$. Considerable variability is observed for $t\in[3,7]$ and $t\in [8,25]$. Starting with $t\approx 30$, the variability between solutions decreases and stabilizes. From these five graphics, we can estimate that the maximal amplitude decreases by a factor $\approx 2$, at 0.015, 0.008, 0.004, 0.0015, and 0.0008. However, it is challenging to make a similar statement about each time point $t$. This observation makes it difficult to estimate the weak and strong convergence rates and the optimal number of samples correctly on each mesh level. They are different for each time $t$ (and for each $\bx$). For some time points, the solution is smooth and requires only a few levels and a few samples on each level. For other points with substantial dynamics, the numbers of levels and samples are higher.

 %We can see that the highest difference, achieved at $t\approx 8$, is decreasing with an approximate rate 2: $(0.013,0.009, 0.0045, 0.0021, 0.001)$. 

\begin{figure}[htbp!]
\begin{center}
  
  \includegraphics[width=0.47\textwidth]{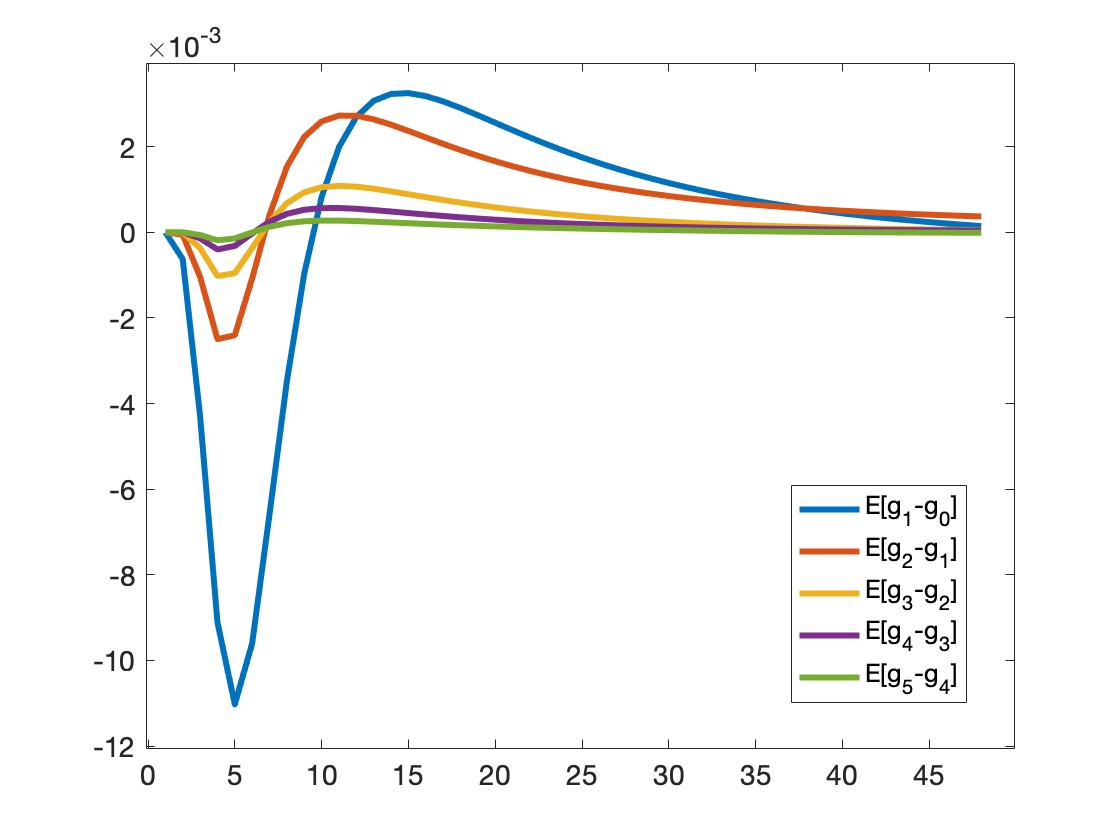}\;
  \includegraphics[width=0.47\textwidth]{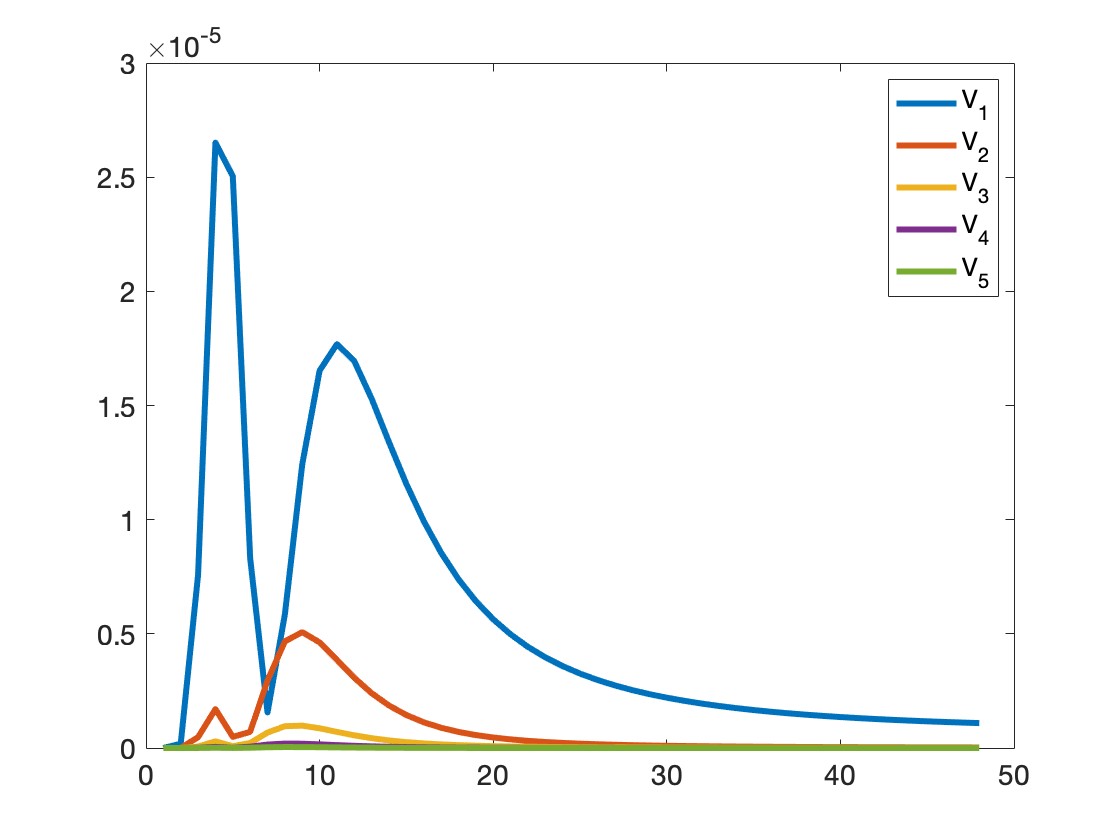}
\caption{(left) Mean and (right) variance of the differences $g_{\ell}-g_{\ell-1}$ vs. time, computed on various levels at the point $(1.60, −0.95)$.}
    \label{fig:mean_var_diffs}
\end{center}
%dima_Henry_mean_var_diffs.m
%dima_compare_MC_with_MLMC.m

\end{figure}

\newpage
Because $g_{\ell}-g_{\ell-1}$ is random, we visualize its mean and variance.
Figure~\ref{fig:mean_var_diffs} demonstrates the mean (left) and variance (right) of the differences in concentrations $g_{\ell}-g_{\ell-1}$, $\ell=1,\ldots,5$. On the left, the amplitude decreases when $\ell$ increases. A slight exception is the blue line for $t\approx 9,10,11$ (right). A possible explanation is that the solutions $g_0$ or $g_1$ are insufficiently accurate. The right image presents how the amplitude of the variances decays. This decay is necessary for the successful work of the MLMC method. We also observe a potential issue; the weak and strong convergence rates vary for various time points $t$. Thus, determining the optimal number of samples $m_{\ell}$ for each level is not possible (only suboptimal). 

At the beginning $t=0$, the variability is zero and starts to increase. We observe changes during a specific time interval, and then the process starts to stabilize after $\approx 45$ time steps. The variability is either unchanging from level to level or decreases.
% and is equal to $(1-0.2\cdot \theta_1)\poro_0$ (in the upper level) and $(1+0.2\cdot \theta_1)\poro_0$ (in the bottom level), where
% \begin{equation}
%     \poro_0=1+0.15\cdot(\theta_2\cdot \cos(\pi x/2) + \theta_2\cdot sin(2\pi y) + \theta_1\cdot \cos(2\pi x)).
% \end{equation}

Table~\ref{tab:adaptiveTS_times} contains average computing times, which are necessary to estimate the number of samples $m_{\ell}$ at each level $\ell$. The fourth column contains the average computing time, and the fifth and sixth columns contain the shortest and longest computing times. The computing time for each simulation varies depending on the number of iterations, which depends on the porosity and permeability. We observed that, after $\approx 6016$~s, the solution is almost unchanging; thus, we restrict this to only $t\in [0, T]$, where $T=6016$. For example, if the number of time steps is $r_{\ell}=188$ (Level 0 in Table~\ref{tab:adaptiveTS_times}), then the time step $\tau = \frac{T}{r_{\ell}}=\frac{6016}{188}=32$~s. 
%Table~\ref{tab:fixedTS_times} is for the settings when the time step $\tau$ is fixed and $\tau=1$ sec. (very fine mesh). We observe that the averaged computing time is growing linear with $n_{\ell}$.

% \begin{table}[htbp!]
% \begin{center}
% \begin{tabular}{|l|l|l|l|l|l|}
% \hline 
% \multirow{2}{*}{Level $\ell$}& 
% \multirow{2}{*}{$n_{\ell}$} &
% \multirow{2}{*}{$r_{\ell}$} &
% \multicolumn{3}{|c|}{Computing times}\\
% \cline{4-6}
%     &&                                  & averaged & min. & max.\\
% \hline
% 0 &    1122 & 6016 &    23.7 & 22   & 25   \\ \hline % 512 grid elements
% 1 &    4290 & 6016 &    48.0 & 46   & 50   \\ \hline % 2048
% 2 &   16770 & 6016 &   133.5 & 126  & 142  \\ \hline % 8192
% 3 &   66306 & 6016 &   501.4 & 455  & 514  \\ \hline % 32768
% 4 &  263682 & 6016 &  1927.4 & 1644 & 1960 \\ \hline % 131072
% 5 & 1051650 & 6016 &  8138.0 & 6453 & 8777 \\ \hline % 524288
% \end{tabular}
% \caption{Number of the degrees of freedom $n_{\ell}$, number of time steps $r_{\ell}$, averaged, minimal, maximal computing times on each level $\ell$.} %RunRechargeRanPoro}
% \label{tab:fixedTS_times}
% \end{center}
% \end{table}
%

%Table~\ref{tab:adaptiveTS_times} lists computing times for the settings when 
The time step $\tau$ is adaptive and changing from $\tau=\frac{6016}{128}=32$~s (very coarse mesh) to $\tau=\frac{6016}{6016}=1$~s (finest mesh). Starting with level $\ell=2$, the average time increases by a factor of eight. These numerical tests confirm the theory in Eq.~(\ref{eq:CompComplexity}), stating that the numerical solver is linear w.r.t. $n_{\ell}$ and $ r_{\ell}$.

\begin{table}[htbp!]
\begin{center}
\begin{tabular}{|l|l|l|l|l|l|l|}
\hline 
\multirow{2}{*}{Level $\ell$}& 
\multirow{2}{*}{$n_{\ell}$} &
\multirow{2}{*}{$r_{\ell}$} &
\multirow{2}{*}{$\tau_{\ell}=6016/r_{\ell}$} &
\multicolumn{3}{|c|}{Computing times ($s_{\ell}$)}\\
\cline{5-7}
    &&   &                               & average & min. & max.\\
\hline
0 &    1122 &  188 & 32 &     1.15 &    0.88 &    1.33 \\ \hline % 512 grid elements
1 &    4290 &  376 & 16 &     4.1  &    3.4  &    4.87 \\ \hline % 2048
2 &   16770 &  752 &  8 &    19.6  &   17.6  &   22    \\ \hline % 8192
3 &   66306 & 1504 &  4 &   136.0  &  128    &  144    \\ \hline % 32768
4 &  263682 & 3008 &  2 &  1004.0  &  891    & 1032    \\ \hline % 131072
5 & 1051650 & 6016 &  1 &  8138.0  & 6430    & 8480    \\ \hline % 524288
\end{tabular}
\caption{Number of degrees of freedom $n_{\ell}$, number of time steps $r_{\ell}$, step size in time $\tau_{\ell}$, average, minimal, and maximal computing times on each level $\ell$.} %RunRechargeRanPoro}. and computing times on each level $\ell=0..5$.} 
\label{tab:adaptiveTS_times} % RunRechargeRanPoro.
\end{center}
\end{table}

%
% Figure~\ref{fig:mean_var_diffs_levels} shows the mean (left) and the variance (right) of the solution for $t=1\ldots 48$ computed on meshes $L=\{0,1,2,3,4,5\}$. Particularly, one can see a very small difference between the mean values (on the left) computed on all meshes. The variances are also very similar.

% \begin{figure}[htbp!]
% \begin{center}
%   \includegraphics[width=0.4\textwidth]{figs/concentr_means_L4_9_point5_RanRechargeRanPoro.jpg}\;
%   \includegraphics[width=0.4\textwidth]{figs/concentr_var_L4_9_RanRechargeRanPoro.jpg}
%   \caption{The mean (left) and the variance (right) for $p_3=(1.60, −0.95),$ $\ell=0,\dots,5$, $t=1\ldots 48$.}
%     \label{fig:mean_var_diffs_levels}
% \end{center}  
% %dima_compare_MC_with_MLMC.m
% \end{figure}
%
\newpage 

With estimates for each level, we can estimate the rates of $\alpha$ and $\beta$ (Eqs.~(\ref{eq:weak_error_model})-(\ref{eq:strong_error_model})) in weak and strong convergences. 

\begin{figure}[htbp!]
\begin{center}
  \includegraphics[width=0.4\textwidth]{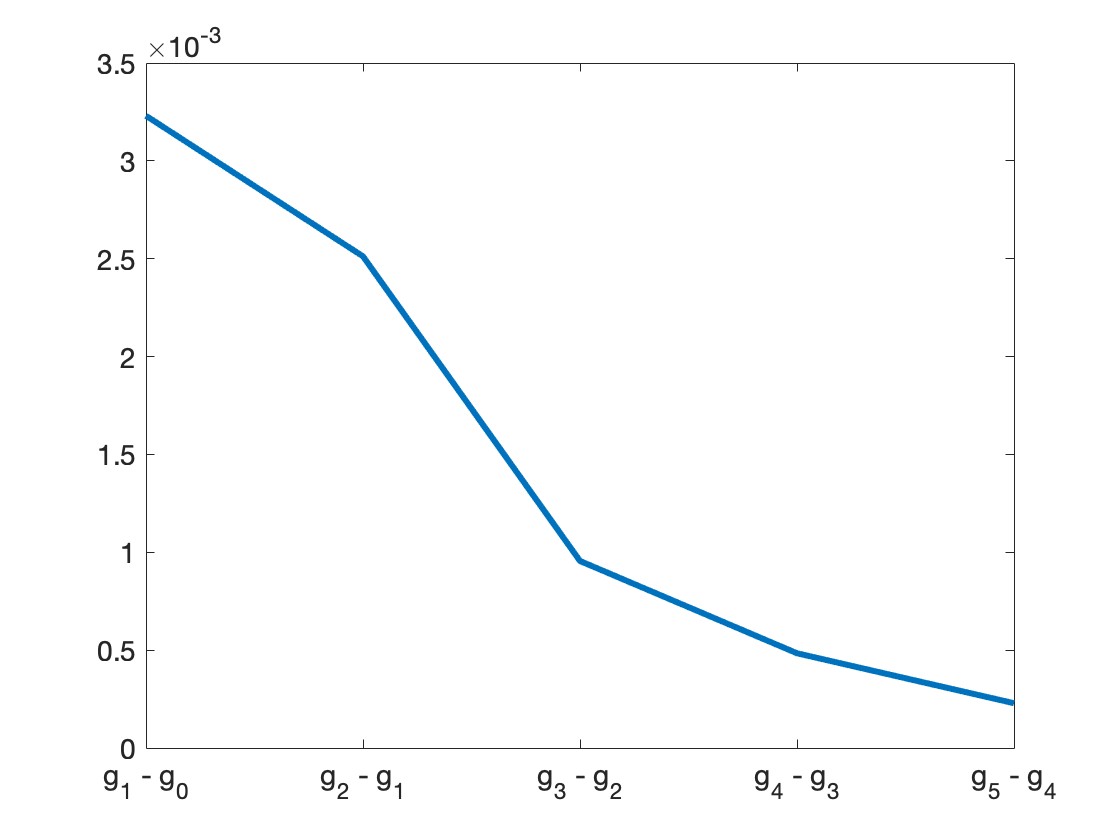}\;
  \includegraphics[width=0.4\textwidth]{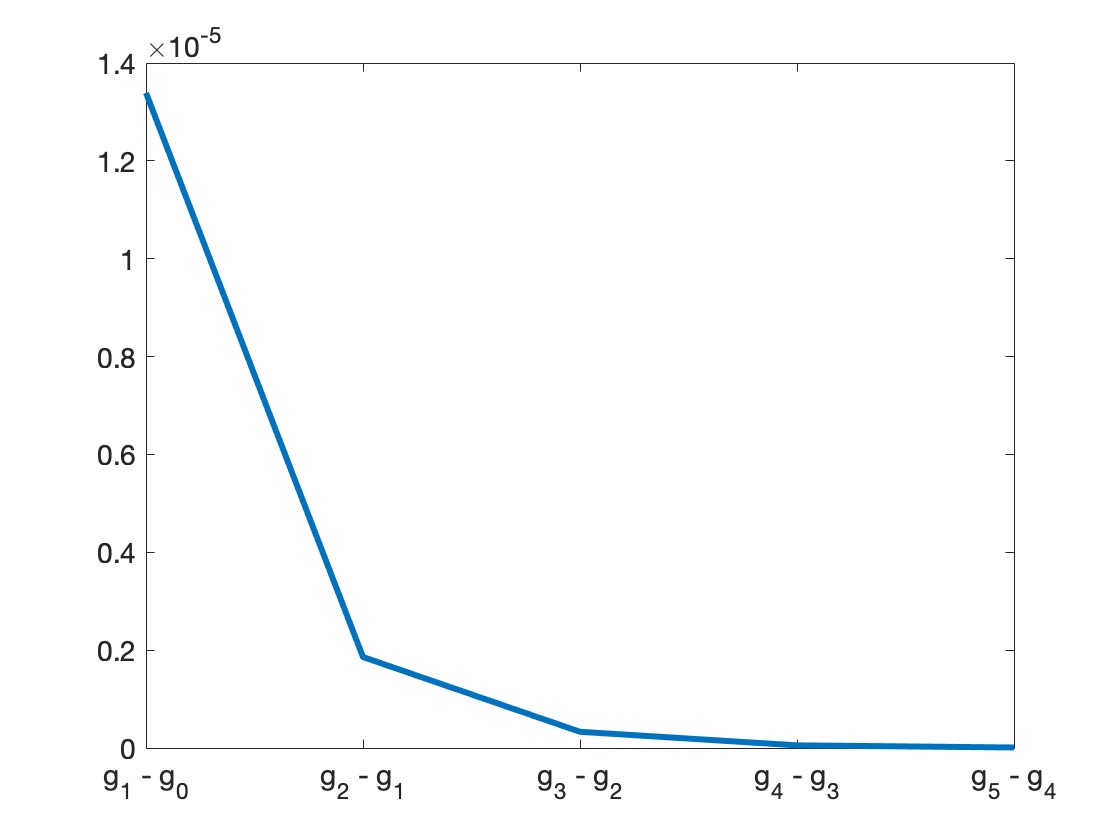}
  \caption{Weak (left) and strong (right) convergences computed for Levels 1 and 0, 2 and 1, 3 and 2, 4 and 3, and 5 and 4 (horizontal axis) at the fixed point $(t,x,y)=(14, 1.60, −0.95)$.}
    \label{fig:weak_strong}
\end{center}  
%dima_compare_MC_with_MLMC.m
\end{figure}

The slope in Fig.~\ref{fig:weak_strong} can be used to estimate the rates of the weak (left) and strong (right) convergences. The differences are indicated on the horizontal axis.

We use computed variances $V_{\ell}$ and computing times (work) $s_{\ell}$ 
from Table~\ref{tab:adaptiveTS_times} to estimate the optimal number of samples $m_{\ell}$ and compute the telescopic sum from Eq.~(\ref{eq:A}) to approximate the expectation.
%\textcolor{red}{TODO: Alex, add how much we won in time}.

Table~\ref{tab:M_ell} lists $m_{\ell}$ for a given total variance $\varepsilon^2$:
\begin{table}[htbp!]
\begin{center}
\begin{tabular}{|l|l|l|l|l|l|l|} \hline
level, $\ell$ & 0 & 1 & 2& 3 & 4 &5 \\ \hline 
$s_{\ell}$ & 1.156 & 4.113  &20.382 &139.0   &993.0  & 8053.0 \\ 
$V_{\ell}$ & 1.4e-5& 0.2e-5 & 0.5e-6& 0.1e-6 &0.5e-7 & 1e-7 \\ \hline
$m_{\ell}(\epsilon^2$ =5e-6$)$ &35  &   7  &   2  &   1  &   1  &   1 \\ \hline
$m_{\ell}(\epsilon^2$ =1e-6$)$ &172  &  35  &   8   &  2  &   1   &  1 \\ \hline
$m_{\ell}(\epsilon^2$ =5e-7$)$ & 343  &  69 &   16  &   3  &   1  &   1 \\ \hline
$m_{\ell}(\epsilon^2$ =1e-7$)$ &
1714    &     344    &      78      &    14     &      4       &    2 \\ \hline
\end{tabular}
\caption{Number of samples $m_{\ell}$ computed using Eq.~(\ref{eq:M_ell}) as a function of the total variance $\epsilon^2$.}
  \label{tab:M_ell}
\end{center}
%Generated in dima_M_ell.m
\end{table}

%\textcolor{red}{2Alex: put here profit compare to MC. By the way, can %MC600 achieve such total variance?}

\begin{figure}[htbp!]
\begin{center}
  \includegraphics[width=0.4\textwidth]{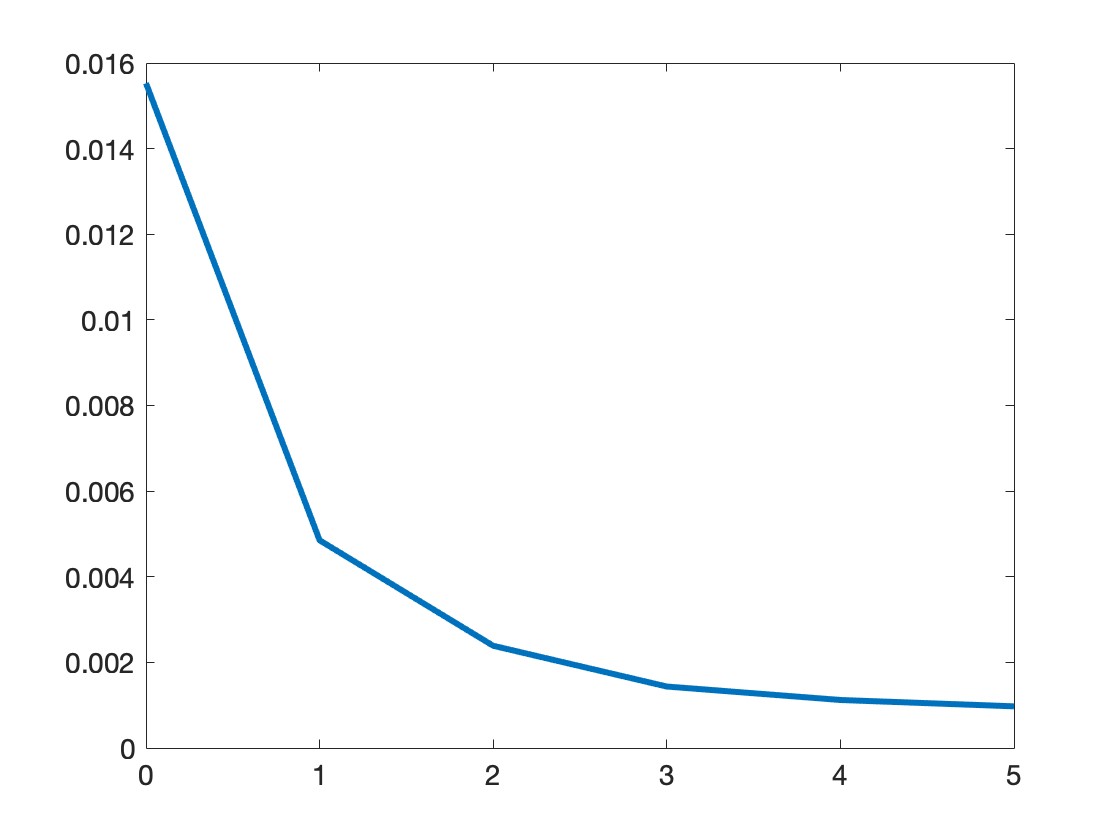}\;
  \includegraphics[width=0.4\textwidth]{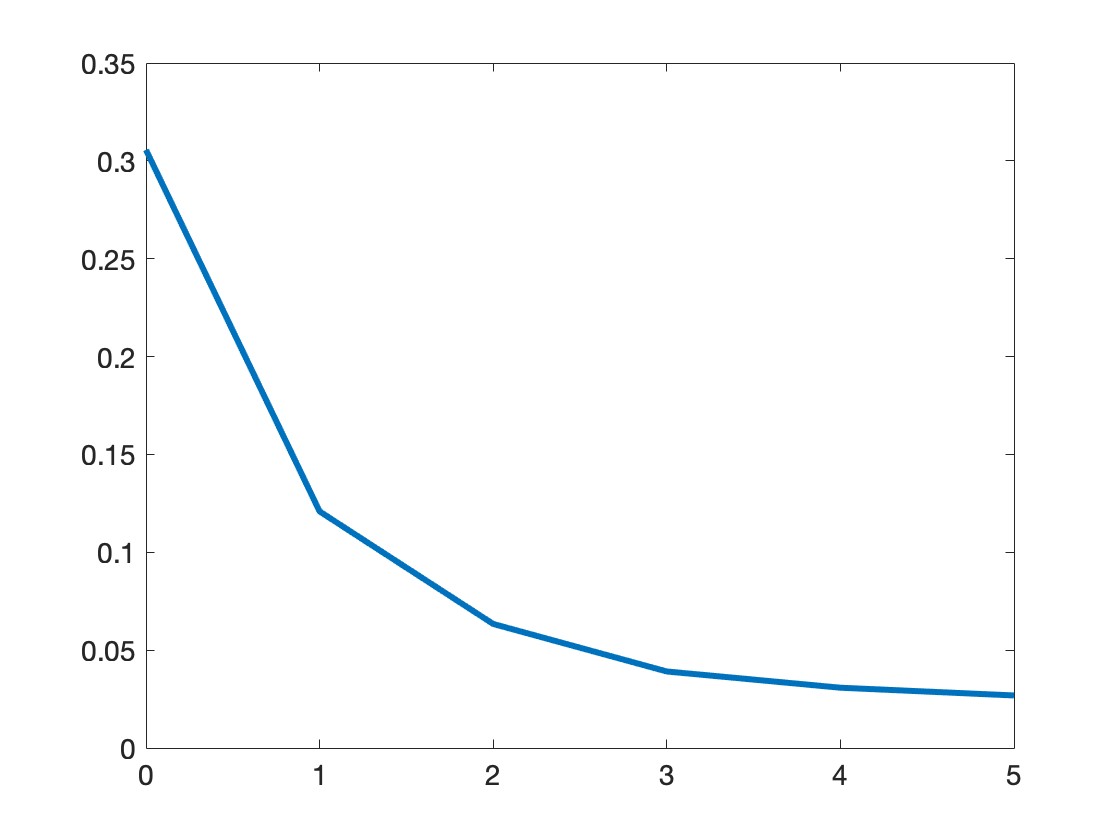}
  \caption{Decay of the absolute (left) and relative (right) errors between the mean values computed on a fine mesh via QMC and on a hierarchy of meshes via MLMC at the fixed point $(t,x,y)=(12, 1.60, −0.95)$. $x$-axis contains the mesh levels.}
    \label{fig:abs_rel_errors_mean}
\end{center}    
%dima_compare_MC_with_MLMC.m
\end{figure}

After the telescopic sum is computed, we can compare the results with the QMC results. Figure~\ref{fig:abs_rel_errors_mean} depicts the decay of the absolute (left) and relative (right) errors vs. levels along the $x$-axes. The 'true' solution was computed using the QMC method on the finest mesh level $L=5$.

% \begin{figure}[htbp!]
% \centering
% \includegraphics[width=0.49\textwidth]{figs/meanMC200.png}\,
% \caption{(First row)}
% \label{fig:1}
% \end{figure}%
%
%

%
%
\section{Conclusion}
\label{sec:Conclusion}
We investigated the applicability and efficiency of the MLMC approach for the Henry-like problem with uncertain porosity, permeability, and recharge. These uncertain parameters were modeled by random fields with three independent random variables. The numerical solution for each random realization was obtained using the well-known ug4 parallel multigrid solver. The number of required random samples on each level was estimated by computing the decay of the variances and computational costs for each level. These estimates depend on the minimization function in Eq.~(\ref{eq:goal_function}).

We also computed the expected value and variance of the mass fraction in the whole domain, the evolution of the pdfs, the solutions at a few preselected points $(t,\bx)$, and the time evolution of the freshwater integral value. We have found that some QoIs require only 2-3 of the coarsest mesh levels, and samples from finer meshes would not significantly improve the result. Note that a different type of porosity in Eq.~(\ref{eq:poro_2levels}) may lead to a different conclusion.

The results show that the MLMC method is faster than the QMC method at the finest mesh. Thus, sampling at different mesh levels makes sense and helps to reduce the overall computational cost.\\
\textbf{Limitations.} 1. It may happen that the QoIs computed on different grid levels are the same (for the given random input parameters). In this case the standard (Q)MC on a coarse mesh will be sufficient. 2. The time dependence is challenging. The optimal number of samples depends on the point $(t,\bx)$ and may be small for some points and large for others. 3. Twenty-four hours may not be sufficient to compute the solution at the sixth mesh level.\\ 
\textbf{Future work.} Our model of porosity in Eq.~(\ref{eq:poro_2levels}) is quite simple. It would be beneficial to consider a more complicated/multiscale/realistic model with more random variables. A more advanced version of MLMC may give better estimates of the number of levels $L$ and the number of samples on each level $m_{\ell}$. Another hot topic is data assimilation and the identification of unknown parameters \cite{Litv_HLIBCov2020,LitvGenton19,Litv_params2016,Rosic2013}.
Known experimental data and measurements of porosity, permeability, velocity or mass fraction could be used to minimise uncertainties. 

\section*{Acknowledgments}
We thank the KAUST HPC support team for assistance with Shaheen II. This work was supported by the Alexander von Humboldt foundation. 

\bibliographystyle{siam}
%\begin{thebibliography}{}

%\end{thebibliography}
% WHERE IS THIS FILE ? \bibliographystyle{spmpsci}
%\bibliography{new_article_Sydney_about_sampling, matthies_BU_paper-1, mybib, references}
%\bibliography{litvinenko_dolgov_khoromskij}
\bibliography{MoCa}

\begin{thebibliography}{10}

\bibitem{Abarca07}
{\sc E.~Abarca, J.~Carrera, X.~S{\'a}nchez-Vila, and M.~Dentz}, {\em
  Anisotropic dispersive henry problem}, Advances in Water Resources, 30
  (2007), pp.~913--926.

\bibitem{babuska_collocation}
{\sc I.~Babu\v{s}ka, F.~Nobile, and R.~Tempone}, {\em A stochastic collocation
  method for elliptic partial differential equations with random input data},
  SIAM J. Numer. Anal., 45 (2007), pp.~1005--1034 (electronic).

\bibitem{babuska2004galerkin}
{\sc I.~Babu\v{s}ka, R.~Tempone, and G.~Zouraris}, {\em Galerkin finite element
  approximations of stochastic elliptic partial differential equations}, SIAM
  Journal on Numerical Analysis, 42 (2004), pp.~800--825.

\bibitem{Templates}
{\sc R.~Barrett, M.~Berry, T.~F. Chan, J.~Demmel, J.~Donato, J.~Dongarra,
  V.~Eijkhout, R.~Pozo, C.~Romine, and H.~van~der Vorst}, {\em Templates for
  the Solution of Linear Systems: Building Blocks for Iterative Methods},
  Society for Industrial and Applied Mathematics, 1994.

\bibitem{BECK22}
{\sc J.~Beck, Y.~Liu, E.~{von Schwerin}, and R.~Tempone}, {\em Goal-oriented
  adaptive finite element multilevel monte carlo with convergence rates},
  Computer Methods in Applied Mechanics and Engineering, 402 (2022), p.~115582.
\newblock A Special Issue in Honor of the Lifetime Achievements of J. Tinsley
  Oden.

\bibitem{Sudret_sparsePCE}
{\sc G.~Blatman and B.~Sudret}, {\em An adaptive algorithm to build up sparse
  polynomial chaos expansions for stochastic finite element analysis},
  Probabilistic Engineering Mechanics, 25 (2010), pp.~183 -- 197.

\bibitem{NowakStochMethods18}
{\sc F.~Bode, T.~Ferr{\'e}, N.~Zigelli, M.~Emmert, and W.~Nowak}, {\em
  Reconnecting stochastic methods with hydrogeological applications: A
  utilitarian uncertainty analysis and risk assessment approach for the design
  of optimal monitoring networks}, Water Resources Research, 54 (2018),
  pp.~2270--2287.

\bibitem{bompard2010}
{\sc M.~Bompard, J.~Peter, and J.-A. D\'{e}sid\'{e}ri}, {\em Surrogate models
  based on function and derivative values for aerodynamic global optimization},
  in Fifth European Conference on Computational Fluid Dynamics, ECCOMAS CFD
  2010, Lisbon, Portugal, 2010.

\bibitem{Griebel_Bungartz}
{\sc H.-J. Bungartz and M.~Griebel}, {\em Sparse grids}, Acta Numer., 13
  (2004), pp.~147--269.

\bibitem{OverviewUncert93}
{\sc J.~Carrera}, {\em An overview of uncertainties in modelling groundwater
  solute transport}, Journal of Contaminant Hydrology, 13 (1993), pp.~23 -- 48.
\newblock Chemistry and Migration of Actinides and Fission Products.

\bibitem{charrier2013}
{\sc J.~Charrier, R.~Scheichl, and A.~L. Teckentrup}, {\em Finite element error
  analysis of elliptic pdes with random coefficients and its application to
  multilevel {M}onte {C}arlo methods}, 51 (2013), pp.~322--352.

\bibitem{chkifa-adapt-stochfem-2015}
{\sc A.~Chkifa, A.~Cohen, and C.~Schwab}, {\em Breaking the curse of
  dimensionality in sparse polynomial approximation of parametric {PDEs}},
  Journal de Mathematiques Pures et Appliques, 103 (2015), pp.~400 -- 428.

\bibitem{MLMC_PDE_anal11}
{\sc K.~Cliffe, M.~Giles, R.~Scheichl, and A.~Teckentrup}, {\em Multilevel
  monte carlo methods and applications to elliptic pdes with random
  coefficients}, Computing and Visualization in Science, 14 (2011), pp.~3--15.

\bibitem{CMLMC}
{\sc N.~Collier, A.-L. Haji-Ali, F.~Nobile, E.~von Schwerin, and R.~Tempone},
  {\em A continuation multilevel monte carlo algorithm}, BIT Numerical
  Mathematics, 55 (2015), pp.~399--432.

\bibitem{ConradMarzouk13}
{\sc P.~Conrad and Y.~Marzouk}, {\em Adaptive smolyak pseudospectral
  approximations}, SIAM Journal on Scientific Computing, 35 (2013),
  pp.~A2643--A2670.

\bibitem{Costa_2006}
{\sc A.~Costa}, {\em Permeability-porosity relationship: A reexamination of the
  kozeny-carman equation based on a fractal pore-space geometry assumption},
  Geophysical Research Letters, 33 (2006).

\bibitem{CREMER15_Fingers}
{\sc C.~Cremer, , and T.~Graf}, {\em Generation of dense plume fingers in
  saturated--unsaturated homogeneous porous media}, Journal of Contaminant
  Hydrology, 173 (2015), pp.~69 -- 82.

\bibitem{Dhal_review22}
{\sc L.~Dhal and S.~Swain}, {\em Understanding and modeling the process of
  seawater intrusion: a review}, 01 2022, pp.~269--290.

\bibitem{DolgLitv15}
{\sc S.~Dolgov, B.~Khoromskij, A.~Litvinenko, and H.~Matthies}, {\em Polynomial
  chaos expansion of random coefficients and the solution of stochastic partial
  differential equations in the tensor train format}, SIAM/ASA Journal on
  Uncertainty Quantification, 3 (2015), pp.~1109--1135.

\bibitem{EIGEL14}
{\sc M.~Eigel, C.~J. Gittelson, C.~Schwab, and E.~Zander}, {\em Adaptive
  stochastic galerkin fem}, Computer Methods in Applied Mechanics and
  Engineering, 270 (2014), pp.~247--269.

\bibitem{Philipp12}
{\sc M.~Espig, W.~Hackbusch, A.~Litvinenko, H.~Matthies, and P.~Waehnert}, {\em
  Efficient low-rank approximation of the stochastic {G}alerkin matrix in
  tensor formats}, Computers and Mathematics with Applications, 67 (2014),
  pp.~818 -- 829.
\newblock High-order Finite Element Approximation for Partial Differential
  Equations.

\bibitem{wahnert-stochgalerkin-2014}
{\sc M.~Espig, W.~Hackbusch, A.~Litvinenko, H.~Matthies, and P.~W{\"a}hnert},
  {\em Efficient low-rank approximation of the stochastic {G}alerkin matrix in
  tensor formats}, Computers and Mathematics with Applications, 67 (2014),
  pp.~818--829.

\bibitem{Frolkovic-DeSchepper-ConvDom}
{\sc P.~Frolkovi{\v c} and H.~{De Schepper}}, {\em Numerical modelling of
  convection dominated transport coupled with density driven flow in porous
  media}, Advances in Water Resources, 24 (2001), pp.~63--72.

\bibitem{Frolkovic-ConsVel}
{\sc P.~Frolkovi\v{c}}, {\em Consistent velocity approximation for density
  driven flow and transport}, in Advanced Computational Methods in Engineering,
  Part 2: Contributed papers, R.~{Van Keer} and at~al., eds., Maastricht, 1998,
  Shaker Publishing, pp.~603--611.

\bibitem{Frolkovic-Knaber-ConsVel}
{\sc P.~Frolkovi\v{c} and P.~Knabner}, {\em Consistent velocity approximations
  in finite element or volume discretizations of density driven flow}, in
  Computational Methods in Water Resources XI, A.~A. Aldama and et~al., eds.,
  Southhampten, 1996, Computational Mechanics Publication, pp.~93--100.

\bibitem{gerstnerGriebel98-numint}
{\sc T.~Gerstner and M.~Griebel}, {\em Numerical integration using sparse
  grids}, Numer. Algorithms, 18 (1998), pp.~209--232.

\bibitem{giles2008}
{\sc M.~B. Giles}, {\em Multilevel {M}onte {C}arlo path simulation}, Operations
  Research, 56 (2008), pp.~607--617.

\bibitem{giles2015}
\leavevmode\vrule height 2pt depth -1.6pt width 23pt, {\em Multilevel {M}onte
  {C}arlo methods}, Acta Numerica, 24 (2015), pp.~259--328.

\bibitem{GiraldiLitv14}
{\sc L.~Giraldi, A.~Litvinenko, D.~Liu, H.~G. Matthies, and A.~Nouy}, {\em To
  be or not to be intrusive? the solution of parametric and stochastic
  equations---the “plain vanilla” galerkin case}, SIAM Journal on
  Scientific Computing, 36 (2014), pp.~A2720--A2744.

\bibitem{giunta2004}
{\sc A.~A. Giunta, M.~S. Eldred, and J.~P. Castro}, {\em Uncertainty
  quantification using response surface approximation}, in 9th ASCE Specialty
  Conference on Probabolistic Mechanics and Structural Reliability,
  Albuquerque, New Mexico, USA, 2004.

\bibitem{Griebel}
{\sc M.~Griebel}, {\em Sparse grids and related approximation schemes for
  higher dimensional problems}, in Foundations of computational mathematics,
  {S}antander 2005, vol.~331 of London Math. Soc. Lecture Note Ser., Cambridge
  Univ. Press, Cambridge, 2006, pp.~106--161.

\bibitem{Hackbusch85}
{\sc W.~Hackbusch}, {\em {Multi-Grid Methods and Applications}}, Springer,
  Berlin, 1985.

\bibitem{Hackbusch_Iter_Sol}
{\sc W.~Hackbusch}, {\em Iterative Solution of Large Sparse Systems of
  Equations}, Springer, New-York, 1994.

\bibitem{ErikOptGeom15}
{\sc A.-L. Haji-Ali, F.~Nobile, E.~von Schwerin, and R.~Tempone}, {\em
  Optimization of mesh hierarchies in multilevel {M}onte {C}arlo samplers},
  Stoch. Partial Differ. Equ. Anal. Comput., 4 (2016), pp.~76--112.

\bibitem{henry1964effects}
{\sc H.~R. Henry}, {\em Effects of dispersion on salt encroachment in coastal
  aquifers, in 'seawater in coastal aquifers'}, US Geological Survey, Water
  Supply Paper, 1613 (1964), pp.~C70--C80.

\bibitem{hoel2012adaptive}
{\sc H.~Hoel, E.~Von~Schwerin, A.~Szepessy, and R.~Tempone}, {\em Adaptive
  multilevel {M}onte {C}arlo simulation}, in Numerical Analysis of Multiscale
  Computations, Springer, 2012, pp.~217--234.

\bibitem{hoel2014implementation}
{\sc H.~Hoel, E.~von Schwerin, A.~Szepessy, and R.~Tempone}, {\em
  Implementation and analysis of an adaptive multilevel monte carlo algorithm},
  Monte Carlo Methods and Applications, 20 (2014), pp.~1--41.

\bibitem{spiterp}
{\sc A.~Klimke}, {\em Sparse grid interpolation
  toolbox,\text{www.ians.uni-stuttgart.de/spinterp/}},  (2008).

\bibitem{Laattoe2013_SeawaterIntr}
{\sc T.~Laattoe, A.~Werner, and C.~Simmons}, {\em Seawater Intrusion Under
  Current Sea-Level Rise: Processes Accompanying Coastline Transgression},
  Springer Netherlands, Dordrecht, 2013, pp.~295--313.

\bibitem{Litv_HLIBCov2020}
{\sc A.~Litvinenko, R.~Kriemann, M.~G. Genton, Y.~Sun, and D.~E. Keyes}, {\em
  Hlibcov: Parallel hierarchical matrix approximation of large covariance
  matrices and likelihoods with applications in parameter identification},
  MethodsX, 7 (2020), p.~100600.

\bibitem{LitLog3D_20}
{\sc A.~Litvinenko, D.~Logashenko, R.~Tempone, G.~Wittum, and D.~Keyes}, {\em
  Solution of the 3d density-driven groundwater flow problem with uncertain
  porosity and permeability}, GEM - International Journal on Geomathematics, 11
  (2020), p.~10.

\bibitem{Litvinenko-UQ-2021}
\leavevmode\vrule height 2pt depth -1.6pt width 23pt, {\em Propagation of
  uncertainties in density-driven flow}, in Sparse Grids and Applications ---
  Munich 2018, H.-J. Bungartz, J.~Garcke, and D.~Pfl{\"u}ger, eds., Cham, 2021,
  Springer International Publishing, pp.~101--126.

\bibitem{LitvGenton19}
{\sc A.~Litvinenko, Y.~Sun, M.~G. Genton, and D.~E. Keyes}, {\em Likelihood
  approximation with hierarchical matrices for large spatial datasets},
  Computational Statistics \& Data Analysis, 137 (2019), pp.~115--132.

\bibitem{Litv_Scattered19}
{\sc A.~Litvinenko, A.~C. Yucel, H.~Bagci, J.~Oppelstrup, E.~Michielssen, and
  R.~Tempone}, {\em Computation of electromagnetic fields scattered from
  objects with uncertain shapes using multilevel monte carlo method}, IEEE
  Journal on Multiscale and Multiphysics Computational Techniques, 4 (2019),
  pp.~37--50.

\bibitem{liu2014}
{\sc D.~Liu and S.~G\"{o}rtz}, {\em Efficient quantification of aerodynamic
  uncertainty due to random geometry perturbations}, in New Results in
  Numerical and Experimental Fluid Mechanics IX, A.~Dillmann et~al., eds.,
  Springer International Publishing, 2014, pp.~65--73.

\bibitem{Loeven2007}
{\sc G.~J.~A. Loeven, J.~A.~S. Witteveen, and H.~Bijl}, {\em A probabilistic
  radial basis function approach for uncertainty quantification}, in
  Proceedings of the NATO RTO-MP-AVT-147 Computational Uncertainty in Military
  Vehicle design symposium, 2007.

\bibitem{matthies2007}
{\sc H.~Matthies}, {\em Uncertainty quantification with stochastic finite
  elements}, in Encyclopedia of Computational Mechanics, E.~Stein, R.~de~Borst,
  and T.~R.~J. Hughes, eds., John Wiley \& Sons, Chichester, 2007.

\bibitem{Litv_params2016}
{\sc H.~G. Matthies, E.~Zander, B.~V. Rosi{\'c}, and A.~Litvinenko}, {\em
  Parameter estimation via conditional expectation: a bayesian inversion},
  Advanced Modeling and Simulation in Engineering Sciences, 3 (2016), p.~24.

\bibitem{Habib09_PCE}
{\sc H.~Najm}, {\em Uncertainty quantification and polynomial chaos techniques
  in computational fluid dynamics}, Annual Review of Fluid Mechanics, 41
  (2009), pp.~35--52.

\bibitem{nobile-sg-mc-2015}
{\sc F.~Nobile, L.~Tamellini, F.~Tesei, and R.~Tempone}, {\em An adaptive
  sparse grid algorithm for elliptic {PDE}s with log-normal diffusion
  coefficient}, MATHICSE Technical Report~04, 2015.

\bibitem{novakRitter97}
{\sc E.~Novak and K.~Ritter}, {\em The curse of dimension and a universal
  method for numerical integration}, in Multivariate approximation and splines
  ({M}annheim, 1996), vol.~125 of Internat. Ser. Numer. Math., Birkh\"auser,
  Basel, 1997, pp.~177--187.

\bibitem{novakRitter99-simple}
\leavevmode\vrule height 2pt depth -1.6pt width 23pt, {\em Simple cubature
  formulas with high polynomial exactness}, Constr. Approx., 15 (1999),
  pp.~499--522.

\bibitem{OLADYSHKIN_PCE}
{\sc S.~Oladyshkin and W.~Nowak}, {\em Data-driven uncertainty quantification
  using the arbitrary polynomial chaos expansion}, Reliability Engineering \&
  System Safety, 106 (2012), pp.~179--190.

\bibitem{Panda_Lake_Perm_vs_Por}
{\sc M.~Panda and W.~Lake}, {\em Estimation of single-phase permeability from
  parameters of particle-size distribution}, AAPG Bull., 78 (1994),
  pp.~1028--1039.

\bibitem{Pape_Clauser_Iffland_1999}
{\sc H.~Pape, C.~Clauser, and J.~Iffland}, {\em Permeability prediction based
  on fractal pore-space geometry}, Geophysics, 64 (1999), pp.~1447--1460.

\bibitem{petrasSmolpak}
{\sc K.~Petras}, {\em {S}molpack---a software for {S}molyak quadrature with
  delayed {C}lenshaw-{C}urtis basis-sequence.
  \text{http://www-public.tu-bs.de:8080/~petras/software.html}}.

\bibitem{POST17_Density-driven}
{\sc V.~Post and G.~Houben}, {\em Density-driven vertical transport of
  saltwater through the freshwater lens on the island of baltrum (germany)
  following the 1962 storm flood}, Journal of Hydrology, 551 (2017), pp.~689 --
  702.
\newblock Investigation of Coastal Aquifers.

\bibitem{radovic1996}
{\sc I.~Radovi\'{c}, I.~Sobol, and R.~Tichy}, {\em Quasi-monte carlo methods
  for numerical integration: Comparison of different low discrepancy
  sequences}, Monte Carlo Methods and Applications, 2 (1996), pp.~1--14.

\bibitem{ReiterLogashenkoVogelWittum2017}
{\sc S.~Reiter, D.~Logashenko, A.~Vogel, and G.~Wittum}, {\em Mesh generation
  for thin layered domains and its application to parallel multigrid simulation
  of groundwater flow}.
\newblock submitted to Comput. Visual Sci., 2017.

\bibitem{ug4_ref1_2013}
{\sc S.~Reiter, A.~Vogel, I.~Heppner, M.~Rupp, and G.~Wittum}, {\em A massively
  parallel geometric multigrid solver on hierarchically distributed grids},
  Computing and Visualization in Science, 16 (2013), pp.~151--164.

\bibitem{Riva2015}
{\sc M.~Riva, A.~Guadagnini, and A.~Dell'Oca}, {\em Probabilistic assessment of
  seawater intrusion under multiple sources of uncertainty}, Advances in Water
  Resources, 75 (2015), pp.~93--104.

\bibitem{Rosic2013}
{\sc B.~Rosi\'c, A.~Ku\v{c}erov\'a, J.~S\'ykora, O.~Pajonk, A.~Litvinenko, and
  H.~Matthies}, {\em Parameter identification in a probabilistic setting},
  Engineering Structures, 50 (2013), pp.~179 -- 196.
\newblock Engineering Structures: Modelling and Computations (special issue
  IASS-IACM 2012).

\bibitem{rubin2003applHydro}
{\sc Y.~Rubin}, {\em Applied stochastic hydrogeology}, Oxford University Press,
  2003.

\bibitem{ScheiderKroehnPueschel2012}
{\sc A.~Schneider, K.-P. Kr{\"o}hn, and A.~P{\"u}schel}, {\em Developing a
  modelling tool for density-driven flow in complex hydrogeological
  structures}, Comput. Visual Sci., 15 (2012), pp.~163--168.

\bibitem{SWLRHEGW-SaltwaterInNorthSea2018}
{\sc A.~Schneider, H.~Zhao, J.~Wolf, D.~Logashenko, S.~Reiter, M.~Howahr,
  M.~Eley, M.~Gelleszun, and H.~Wiederhold}, {\em Modeling saltwater intrusion
  scenarios for a coastal aquifer at the german north sea}, E3S Web Conf., 54
  (2018), p.~00031.

\bibitem{Simpson2003}
{\sc M.~J. Simpson and T.~Clement}, {\em {Theoretical Analysis of the
  worthiness of Henry and {E}lder problems as benchmarks of density-dependent
  groundwater flow models}}, Adv. Water. Resour., 26 (2003), pp.~17--31.

\bibitem{Simpson04_Henry}
{\sc M.~J. {Simpson} and T.~P. {Clement}}, {\em {Improving the worthiness of
  the Henry problem as a benchmark for density-dependent groundwater flow
  models}}, Water Resources Research, 40 (2004), p.~W01504.

\bibitem{smoljak63}
{\sc S.~A. Smolyak}, {\em Quadrature and interpolation formulas for tensor
  products of certain classes of functions.}, Sov. Math. Dokl., 4 (1963),
  pp.~240--243.

\bibitem{Stoeckl}
{\sc L.~Stoeckl, M.~Walther, and L.~K. Morgan}, {\em Physical and numerical
  modelling of post-pumping seawater intrusion}, Geofluids, 2019 (2019).

\bibitem{TARTAKOVSKYI_Risk13}
{\sc D.~Tartakovsky}, {\em Assessment and management of risk in subsurface
  hydrology: A review and perspective}, Advances in Water Resources, 51 (2013),
  pp.~247 -- 260.
\newblock 35th Year Anniversary Issue.

\bibitem{teckentrup2013further}
{\sc A.~Teckentrup, R.~Scheichl, M.~Giles, and E.~Ullmann}, {\em Further
  analysis of multilevel monte carlo methods for elliptic pdes with random
  coefficients}, Numerische Mathematik, 125 (2013), pp.~569--600.

\bibitem{SoilOverview16}
{\sc H.~Vereecken, A.~Schnepf, J.~Hopmans, M.~Javaux, D.~Or, T.~Roose,
  J.~Vanderborght, M.~Young, W.~Amelung, M.~Aitkenhead, S.~Allison,
  S.~Assouline, P.~Baveye, M.~Berli, N.~Br{\"u}ggemann, P.~Finke, M.~Flury,
  T.~Gaiser, G.~Govers, T.~Ghezzehei, P.~Hallett, H.~Hendricks~Franssen,
  J.~Heppell, R.~Horn, J.~Huisman, D.~Jacques, F.~Jonard, S.~Kollet,
  F.~Lafolie, K.~Lamorski, D.~Leitner, A.~McBratney, B.~Minasny, C.~Montzka,
  W.~Nowak, Y.~Pachepsky, J.~Padarian, N.~Romano, K.~Roth, Y.~Rothfuss,
  E.~Rowe, A.~Schwen, J.~{\v S}im{\r u}nek, A.~Tiktak, J.~Van~Dam, S.~van~der
  Zee, H.~Vogel, J.~Vrugt, T.~W{\"o}hling, and I.~Young}, {\em Modeling soil
  processes: Review, key challenges, and new perspectives}, Vadose Zone
  Journal, 15 (2016), p.~vzj2015.09.0131.

\bibitem{ug4_ref2_2013}
{\sc A.~Vogel, S.~Reiter, M.~Rupp, A.~N{\"a}gel, and G.~Wittum}, {\em Ug 4: A
  novel flexible software system for simulating pde based models on high
  performance computers}, Computing and Visualization in Science, 16 (2013),
  pp.~165--179.

\bibitem{Voss_Souza}
{\sc C.~Voss and W.~Souza}, {\em Variable density flow and solute transport
  simulation of regional aquifers containing a narrow freshwater-saltwater
  transition zone}, Water Resources Research, 23 (1987), pp.~1851--1866.

\bibitem{Dongbin}
{\sc D.~Xiu}, {\em Fast numerical methods for stochastic computations: A
  review}, Commun. Comput. Phys., 5, No. 2-4 (2009), pp.~242--272.

\bibitem{xiuKarniadakis02a}
{\sc D.~Xiu and G.~E. Karniadakis}, {\em The {W}iener-{A}skey polynomial chaos
  for stochastic differential equations}, SIAM J. Sci. Comput., 24 (2002),
  pp.~619--644.

\end{thebibliography}

\end{document}